\documentclass{amsart}

\usepackage{amsmath,amsfonts, amssymb}									
\usepackage{graphicx}

\newtheorem{theorem}{Theorem}[section]
\newtheorem*{theorem A}{Theorem A}
\newtheorem*{theorem B}{N\"olker's Theorem}
\newtheorem{lemma}{Lemma}[section]

\newtheorem{proposition}{Proposition}[section]

\newtheorem{corollary}{Corollary}[section]

\theoremstyle{remark}
\newtheorem{remark}{Remark}[section]
\theoremstyle{remark}

\theoremstyle{definition}

\newtheorem{definition}{Definition}[section]
\newtheorem{example}{Example}[section]
\numberwithin{equation}{section}
\def\({\left ( }
\def\){\right )}
\def\<{\left < }
\def\>{\right >}

\def\h{\mathbb H}
\def\s{\mathbb S}
\def\r{\mathbb R}

\def\u{\mathbb U}

\def\e{\hbox{\bf E}}
\def\t{\hbox{\bf T}}
\def\n{\hbox{\bf N}}
\def\b{\hbox{\bf B}}

\def\c{\hbox{\bf c}}
 \def\T{ \mathcal T}

\begin{document}

\vspace{2cm}

\title[Differential Geometry in Lorentz-Minkowski space]{Differential Geometry of Curves and Surfaces in Lorentz-Minkowski space}

\author{Rafael L\'opez}
\address{Departamento de Geometr\'{\i}a y Topolog\'{\i}a\\
    Universidad de Granada-SPAIN}
\email{rcamino@ugr.es}

\thanks{Partially supported by a MEC-FEDER
 grant no. MTM2011-22547 and
Junta de Andaluc\'{\i}a grant no. P09-FQM-5088.}


\date{}


\begin{abstract}
We review part of the classical theory of curves and surfaces in $3$-dimensional Lorentz-Minkowski space. We focus in spacelike surfaces with constant mean curvature pointing the differences and similarities with the Euclidean space.
\end{abstract}
\maketitle

\tableofcontents
\addtocontents{toc}{\protect\setcounter{tocdepth}{1}}

The title of this work is motivated by the book of M. P. do Carmo, Differential Geometry of Curves and Surfaces (\cite{dc}), and its origin was a mini-course given by the author in September 2008 in  the Instituto de Matem\'atica e Estat\'{\i}stica (IME-USP) of the University of Sao Paulo, Brazil. The main purpose is  how to approach to the study of curves and surfaces in Lorentz-Minkowski space when one has basic concepts of curves and surfaces in Euclidean space. With this in mind, we will consider some of the topics that appear in \cite{dc} and we ask what type of differences and similarities are between both scenarios. Originally, these notes were posted in the ArXiv server \cite{lo1} such as I wrote in a first and quick version. I realize that the paper has been cited in many articles depicted the number of typographic mistakes of the original version. In occasion of this special volume dedicated to the memory of Prof. Franki Dillen by the International Electronic Journal of Geometry, I have the opportunity to improve the original draft. Furthermore, this article is a tribute to the work of Prof. Dillen and its influence, specially in the field of submanifolds in Lorentz-Minkowski space. In my case, I point out the article \cite{dk} on ruled surfaces in Minkowski space and the set of articles that appeared in the series of volumes `Geometry and Topology of Submanifolds' by World Scientific which took place in the nineties and where  Prof. Dillen had a high activity.

There is not a textbook with a systematic  study of curves and surfaces in Lorentz-Minkowski space such as it occurs in the Euclidean space. Some of the topics of this paper  can be found in some books and articles. For curves, it is usually cited \cite{wa}, which it is a doctoral thesis and not easily available. Very recently there is some space for some type of curves in \cite[p. 33]{ku}. For surfaces, I refer \cite[Part I]{gr} and \cite[Ch. 3]{we}, although here the focus lies on timelike surfaces.  A general reference including many topics in semi-Riemannian geometry is the classical book of O'Neill \cite{on}.   We begin in section 1 with an introduction to the metric space of Lorentz-Minkowski space $\e_1^3$ with an interest on the isometries of this space. The second section develops the Frenet equations for curves in $\e_1^3$.  In section 3 we study  surfaces in $\e_1^3$ with special attention on spacelike surfaces. We give the notion of mean curvature and Gauss curvature and we show  many examples of surfaces.   In section 4 we consider spacelike surfaces with constant mean curvature and finally, section 5  emphasizes the connection between this class of surfaces and the theory of elliptic equations.

\section{The Lorentz-Minkowski space $\e_1^3$}\label{section1}

\subsection{Basic definitions}

Let $\r^3$  denote the real vector space with its usual vector structure. Denote by  $B_u=\{E_1,E_2,E_3\}$ the canonical basis of $\r^3$, that is,
$$E_1=(1,0,0),\ \ E_2=(0,1,0),\ \ E_3=(0,0,1).$$
We denote $(x,y,z)$ the coordinates of a vector with respect to $B_u$. We also consider in $\r^3$ its affine structure, and we will say ``horizontal" or ``vertical" in its usual sense.

\begin{definition} The Lorentz-Minkowski space   is the metric space $\e_1^3=(\r^3,\langle,\rangle)$ where the metric $\langle,\rangle$  is
$$\langle u,v\rangle=u_1v_1+u_2 v_2-u_3 v_3,\ \ u=(u_1,u_2,u_3), v=(v_1,v_2,v_3),$$
which is called the Lorentzian metric.
\end{definition}
We also use the terminology   Minkowski space  and   Minkowski metric to refer the space and the metric, respectively. The Lorentzian metric  is a non-degenerate metric of index  $1$.
 The vector space $\r^3$ also supports the Euclidean metric, which will be denoted by  $\langle,\rangle_{e}$. We write the $3$-dimensional Euclidean space as  $\e^3=(\r^3,\langle,\rangle_{e})$ to distinguish from the Lorentz-Minkowski space.

\begin{definition}
A vector  $v\in\e_1^3$ is said
\begin{enumerate}
\item spacelike if  $\langle v,v\rangle>0$ or $v=0$,
\item timelike if  $\langle v,v\rangle<0$ and
\item  lightlike if  $\langle v,v\rangle=0$ and  $v\not=0$.
\end{enumerate}
\end{definition}

The \emph{light-cone} of $\e_1^3$ is the set of all lightlike vectors  of $\e_1^3$:
$$\mathcal{C}=\{(x,y,z)\in\e_1^3:x^2+y^2-z^2=0\}-\{(0,0,0)\}.$$
 The set of timelike vectors  is
$$\mathcal{T}=\{(x,y,z)\in\e_1^3:x^2+y^2-z^2<0\}.$$
We observe that both $\mathcal{T}$ and $\mathcal{C}$ have two connected components.

Given $U\subset\r^3$ a vector subspace, we consider the induced metric $\langle,\rangle_{U}$:
$$\langle u,v\rangle_{U}=\langle u,v\rangle,\ \ u,v\in U.$$
The metric on $U$ classifies in one of the next three types:
\begin{enumerate}
\item The metric is positive definite and $U$ is called spacelike.
\item The metric has index $1$ and $U$ is said timelike.
\item The metric is   degenerate and $U$ is called lightlike.
\end{enumerate}

 The  causal character of a vector or a subspace is the property to be spacelike, timelike or lightlike. In what follows, we give  some characterizations and properties of the causality of a subspace of $\e_1^3$.

\begin{proposition}\label{pr-2-21}
Let $U\subset\e_1^3$ be a vector subspace.
\begin{enumerate}
\item    $\mbox{dim}(U^\bot)=3-\mbox{dim}(U)$.
\item    $(U^\bot)^\bot=U$.
\item If $U$ is non-degenerate, then $U^\bot$ is a non-degenerate subspace.
\item  $U$ is timelike (resp. spacelike, lightlike) if and only if  $U^\bot$ is spacelike (resp. timelike, lightlike).
\item  If $v$ is timelike or spacelike, then $\e_1^3=\mbox{Span}\{v\}\oplus \mbox{Span}\{v\}^\bot$.
\end{enumerate}
\end{proposition}
Comparing with  Euclidean space $\e^3$, the existence of timelike and lightlike vectors give some `strange' properties, as the following:

\begin{proposition}\label{pr-1-1} \begin{enumerate}
\item Two lightlike vectors $u,v\in\e_1^3$  are  linearly dependent if and only  $\langle u,v\rangle=0$.
\item If $u$ and $v$ are two timelike or lightlike vectors with  $\langle u,v\rangle=0$, then they are lightlike vectors.
\item If $u$ and $v$ are two timelike vectors, then $º\langle u,v\rangle\not=0$.
\item If $U$ is a lightlike subspace, then  $\mbox{dim}(U\cap U^\bot)=1$.
\end{enumerate}
\end{proposition}

\begin{proposition}\label{1lightlike} Let $P\subset\e_1^3$ be a vector plane. The following statements are equivalent:
\begin{enumerate}
\item $P$ is a timelike subspace.
\item $P$ contains two  linearly independent lightlike vectors.
\item $P$ contains a timelike vector.
\end{enumerate}
\end{proposition}

  We now characterize lightlike subspaces.

\begin{proposition} Let $U$ be a vector subspace of  $\e_1^3$. The following statements are equivalent:
\begin{enumerate}
\item $U$ is a lightlike subspace.
\item $U$ contains a lightlike vector but not a timelike one.
\item $U\cap \mathcal{ C}=L-\{0\}$, and $\mbox{dim }L=1$.
\end{enumerate}
\end{proposition}

From the Euclidean viewpoint, the next result is useful.

\begin{proposition}\label{pr-n} Let $P\subset\e_1^3$ be a vector plane. Denote by  $\vec{n}_e$ an orthogonal vector with respect to the Euclidean metric. Then $P$ is a spacelike (resp. timelike, lightlike) plane if and only if $\vec{n}_e$ is a timelike (resp. spacelike, lightlike) vector.
\end{proposition}

\begin{proof} If $P$ writes as $P=\{(x,y,z)\in\r^3:ax+by+cz=0\}$, then  $\vec{n}_e$ is proportional to the vector $(a,b,c)$. We may also write $P$ as
$$P=\{(x,y,z)\in\r^3:ax+by- (-c)z=0\}=\mbox{Span}\{(a,b,-c)\}^\bot.$$
The causal character of  $(a,b,-c)$ is the same than  $\vec{n}_e$ and Prop. \ref{pr-2-21} proves the result.
\end{proof}

We define the norm (or modulus) of a vector.

\begin{definition} Given $u\in \e_1^3$,  the norm  of  $u$ is $|u|=\sqrt{|\langle u,u\rangle|}$. The vector $u$ is called unitary if its norm is  $1$.
\end{definition}

\begin{proposition} If  $P=\mbox{Span}\{v\}^\bot$ is a spacelike plane, then
$$|v|_{e}\geq |v|.$$
\end{proposition}
\begin{proof}It suffices if $|v|=1$. Assume $\vec{n}_e=(a,b,c)$, with $a^2+b^2+c^2=1$. Then
$$v=\pm\frac{(a,b,-c)}{\sqrt{c^2-a^2-b^2}}.$$
The Euclidean norm $|v|_e$ is $$|v|_{e}^2=\frac{a^2+b^2+c^2}{c^2-a^2-b^2}=\frac{1}{c^2-a^2-b^2}\geq 1$$
because $c^2-a^2-b^2=1-2(a^2+b^2)\leq 1$.
\end{proof}
This result justifies why when one draws a unit orthogonal vector to a spacelike plane, the Euclidean size is greater than $1$.

\subsection{Timelike vectors}
If $u$ is a timelike vector,   the  \emph{timelike cone}   of  $u$ is
$$C(u)=\{v\in\T:\langle u,v\rangle<0\}.$$
This set is non-empty since $u\in C(u)$. Moreover, if $v$ is other timelike vector, and using $\langle u,v\rangle\not =0$ (Prop. \ref{pr-1-1}), then $\langle u,v\rangle<0$ or $\langle u,v\rangle>0$. This means that
$\mathcal{ T}$ is the disjoint union $\mathcal{T}\doteq C(u)\cup C(-u)$, with $C(u)\cap C(-u)=\emptyset$.
Some properties of  timelike cones are:

\begin{proposition} \begin{enumerate}
\item Two timelike vectors $u$ and $v$ lie in the same timelike cone if and only if  $\langle u,v\rangle<0$.
\item $u\in C(v)$ if and only if $C(u)=C(v)$.
\item The timelike cones are convex sets.
\end{enumerate}
\end{proposition}

\begin{remark} The existence of timelike cones occurs because ${\mathcal T}$ has two components. For lightlike vectors there is a similar situation since ${\mathcal C}$ has two components, namely, ${\mathcal C}^+=\{p\in {\mathcal C}: z>0\}$ and ${\mathcal C}^{-}=\{p\in {\mathcal C}: z<0\}$. Given  two linearly independent  vectors $u,v\in {\mathcal C}$, then $\langle u,v\rangle\not=0$ by Prop. \ref{pr-1-1}. In this case $\langle u,v\rangle<0$ if and only if both vectors are in the same component of ${\mathcal C}$.

Up an order, a basis $\{e_1,e_2,e_3\}$ of $\e_1^3$ is called a null basis (or null frame) if $e_1$ is a unit spacelike vector, $e_2,e_3$ are lightlike vectors in $\mbox{Span}\{e_1\}^\bot$ such that $\langle e_2,e_3\rangle=-1$. In particular, $e_2$ and $e_3$  belong to the same component of ${\mathcal C}$.
\end{remark}

A difference that we find between  $\e^3$ and $\e_1^3$ refers to the Cauchy-Schwarz inequality. Recall that if   $u,v\in\e^3$,  the Cauchy-Schwarz inequality asserts $|\langle u,v\rangle|\leq |u||v|$ and the equality holds if and only if $u, v$ are proportional.

In Minkowski space, and for timelike vectors, there exists a `reverse' inequality called  the {\it backwards Cauchy-Schwarz inequality} (\cite[p. 144]{on}) .

\begin{theorem} Let $u,v\in\e_1^3$  two timelike vectors.  Then
$$|\langle u,v\rangle|\geq |u| |v|$$
and the equality holds if and only if $u$ and $v$ are  proportional.
In the case that both vectors lie in the same timelike cone, there exists a unique number $\varphi\geq 0$ such that
\begin{equation}\label{angleh}
\langle u,v\rangle=-|u| |v|\cosh\varphi.
\end{equation}
The number $\varphi$ is called  the hyperbolic angle   between $u$ and $v$.
\end{theorem}

\begin{proof}
Consider two linearly independent timelike vectors  $u$ and $v$. Then $U=\mbox{Span}\{u,v\}$ is a timelike plane. By Prop. \ref{1lightlike} the equation on $\lambda$
$$\langle u+\lambda v,u+\lambda v\rangle= \langle u,u\rangle+2\lambda \langle u,v\rangle+\lambda^2\langle v,v\rangle=0$$
 has a least two solutions.  In particular, the discriminant of the quadratic equation must be positive, that is,
$$\langle u,v\rangle^2>\langle u,u\rangle\langle v,v\rangle.$$
This shows the inequality in the case that $u$ and $v$ are linearly independent. On the other hand, if they are proportional, then we obtain directly the
equality.

For the second part of the theorem, we write
\begin{equation}\label{cs}
\frac{\langle u,v\rangle^2}{(|u| |v|)^2}\geq 1.
\end{equation}
If $u$ and $v$ lie in the same timelike cone, then $\langle u,v\rangle<0$ and the expression  \eqref{cs}  implies
$$\frac{-\langle u,v\rangle}{|u| |v|}\geq 1.$$
As the hyperbolic cosine function   $\cosh:[0,\infty)\rightarrow [1,\infty)$ is one-to-one, there exists a unique number $\varphi\in [0,\infty)$ such that
$$\cosh\varphi=\frac{-\langle u,v\rangle}{|u| |v|}.$$
 \end{proof}

 After the definition of the angle between two vectors that lie in the same timelike cone, we ask how to define the angle between two any vectors $u,v\in\e_1^3$. Assume that $u,v$ are linearly independent and that $u,v$ are not lightlike. The angle is defined depending on the plane $P$ determined by $u$ and $v$. The induced metric on $P$ can be Riemannian, Lorentzian or degenerate.
 \begin{enumerate}
 \item If the plane is Riemannian, then the definition of the angle between both (spacelike) vectors is the usual as in Euclidean space.
\item  If the plane is Lorentzian, then it is isometric to the Lorentz-Minkowski plane $\e_1^2$ and  an isometry does not change the definition of angle. We have defined the angle for two timelike vectors in the same timelike cone.  It suffices to consider that $u$ and $v$ are unitary. The set $\u_1^2$ of unit vectors of $\e_1^2$ has four components, namely,
 $$\h^1_{+}=\{(x,y)\in\e_1^2:x^2-y^2=-1, y>0\}, \ \h^1_{-}=\{(x,y)\in\e_1^2:x^2-y^2=-1, y<0\}$$
 $$\s_1^{1+}=\{(x,y)\in\e_1^2:x^2-y^2=1, x>0\}, \ \s_1^{1-}=\{(x,y)\in\e_1^2:x^2-y^2=1, x<0\}.$$
The vectors in $\h^1_{+}\cup\h_{-}^1$ are timelike and in $\s_1^{1+}\cup \s_1^{1-}$ are spacelike.

 \begin{remark}\label{note} We point out that by changing $(x,y)$ by $(y,x)$, the plane $\e_1^2$ changes by $\r^2$ equipped with the metric $-(dy)^2+(dx)^2$. Then a spacelike vector (resp. timelike) of $\e_1^2$ converts to a timelike (resp. spacelike) vector of the new metric space.
\end{remark}

Consider two unit spacelike vectors $u,v\in \u_1^2$ and,  \emph{in addition}, we assume that they lie in the same component of $\u_1^2$, that is, $u,v\in \s_1^{1+}$ or $u,v\in\s_1^{1-}$. By the above remark, we conclude then that $\langle u,v\rangle\geq 1$.
\begin{definition} Let $u,v\in\e_1^2$ be two non-zero  spacelike vectors such that $u/|u|$ and $v/|v|$ lie in the same component of $\u_1^2$. Then the angle $\angle(u,v)=\varphi$ is the unique number $\varphi\in[0,\infty) $ such that
\begin{equation}\label{angles}
\cosh\varphi=\frac{\langle u,v\rangle}{|u||v|}.
\end{equation}
\end{definition}

We do not define the angle between two unit spacelike (or timelike) vectors of $\e_1^2$ that do not belong to the same component of $\u_1^2$, neither the angle between a spacelike and timelike vectors.  See more justifications in Th. \ref{te-angle}.
 \item Finally, a third case appears if the plane containing both vectors is lightlike. Necessarily, $u$ and $v$ are not timelike. Here we do not define the angle between two (spacelike) vectors.
     \end{enumerate}

We give the definition of  {\it timelike orientation}. First, we recall the notion of orientation in any vector space. For this,   in the set of all ordered basis of $\r^3$, we consider the equivalence relation $\thicksim_o$ given by $ B \thicksim_o B'$ if the change of basis matrix  has positive determinant. There exist exactly two equivalence classes,  called orientations of  $\r^3$. Fix a basis $B$.  Given other basis $B'$,  we say that   $B'$ is positively oriented if  $B'\thicksim_o B$; on the contrary, we say that $B'$ is negatively oriented. Moreover, the choice of the ordered pair $(\r^3,[B])$ reads saying that $\r^3$ is oriented (by $B$),

In Minkowski space $\e_1^3$ and since the background space is  $\r^3$, the notion of orientation is the same.  The timelike orientation that we  introduce is a metric concept because we use the Lorentzian metric
 $\langle,\rangle$ and thus,  there is not relation with the above notion.

In $\e_1^3$ we consider the set $\mathcal{B}$ of all ordered orthonormal basis where if  $B=\{e_1,e_2,e_3\}\in \mathcal{B}$, then $e_3$ is a timelike vector. If $B=\{e_1,e_2,e_3\}$ and $B'=\{e_1',e_2',e_3'\}$ are two basis, we define the equivalence relation  $\sim$ by
$$B\sim B'\mbox{ if $e_3$ and $e_3'$ lies in the same timelike cone},$$
that is, if  $\langle e_3,e_3'\rangle<0$. The equivalence relation $\sim$ determines two equivalence classes, which are called  \emph{ timelike orientations}. Moreover, each class determines a unique timelike cone  which is defined by the third vector $e_3$ of $B$. Conversely, given a timelike cone, there exists a unique timelike orientation in such way that
any  basis $B$ that belongs to this orientation has the last vector $e_3$ lies in such timelike cone.

We say that  $\e_1^3$ is \emph{timelike oriented} if we fix a timelike orientation, that is, we consider an ordered pair $(\e_1^3,[B])$ for some $B$.

\begin{definition} Let $E_3=(0,0,1)$. Given a timelike vector  $v$, we say that  $v$ is future-directed (resp. past-directed) if
 $v\in C(E_3)$, that is, if $\langle v,E_3\rangle<0$ (resp. $v\in C(-E_3)$, or $\langle v,E_3\rangle>0$).
\end{definition}
In coordinates,  $v=(v_1,v_2,v_3)$ is future-directed if  $v_3>0$. Thus if we fix the timelike cone $C(E_3)$ then we have associated a timelike orientation. Then we say that an orthonormal basis $B=\{e_1,e_2,e_3\}$ is future-directed if $e_3$ is future-directed, or equivalently, if $e_3\in C(E_3)$.

We end this introduction with the definition of the  vector product.

\begin{definition} If $u,v\in\e_1^3$,  the Lorentzian vector product of $u$ and $v$ is to the unique vector denoted by  $u\times v$ that
satisfies
\begin{equation}\label{1pro}
\langle u\times v,w\rangle=\mbox{det }(u,v,w),
\end{equation}
where $\mbox{det}(u,v,w)$ is the determinant of the matrix obtained by placing by columns the coordinates of the three vectors
   $u$, $v$ and $w$ with respect to $B_u$.
\end{definition}
The bilinearity of the metric assures the existence and uniqueness of this vector $u\times v$. By taking $w$ in \eqref{1pro} each one of the vectors $E_i$ of $B_u$, we obtain the expression of $u\times v$ in coordinates with respect to $B_u$:
$$u\times v=\Bigg|\begin{array}{ccc}i&j&-k\\u_1&u_2&u_3\\v_1&v_2&v_3\end{array}\Bigg|.$$
Thus, if we denote by $u\times_e v$ the Euclidean vector product, then $u\times v$ is the reflection of  $u\times_{e}v$ with respect to the plane  of equation $z=0$. Let us observe that if $u$ and $v$ are two non-degenerate vectors, then $B=\{u,v,u\times v\}$ is a basis of $\e_1^3$. However, and in contrast to the Euclidean space,  the causal character of $u$ and $v$ determines  if the basis is  or is not positively oriented. Exactly, if $u,v$ are spacelike, then $u\times v$ is timelike and $B$ is negatively oriented because $det(u,v,u\times v)=\langle u\times v,u\times v\rangle<0$. If $u$ and $v$ have different causal character, then $B$ is positively oriented.

\subsection{Isometries of  $\e_1^3$}

We  consider $O_1(3)$ the set of all vector isometries  of $\e_1^3$. The matrix expression $A$ of an isometry with respect to an orthonormal basis satisfies $A^tGA=G$, where
$$G=\left(\begin{array}{ccc}1&0&0\\ 0&0&1\\ 0&0&-1\end{array}\right).$$
In other terms, we express $O_1(3)$ as the set of matrices
$$O_1(3)=\{A\in Gl(3,\r): A^tGA=G\}.$$
In particular,  $\mbox{det}(A)=\pm 1$. This means that $O_1(3)$ has at least two connected components.
Denote by $SO_1(3)$ the set of isometries with determinant $1$. The set $SO_1(3)$ is called   the \emph{special Lorentz group} and it is related with the notion of orientation of $\r^3$. Exactly, given an orientation $B\in \mathcal{B}$,  $B'\in \mathcal{B}$  is positively oriented if the matrix $A$ of change of basis belongs to $SO_1(3)$.

We define the \emph{ortocrone group}  by
$$O_1^+(3)=\{A\in O_1(3): \mbox{$A$ preserves the timelike orientation}\}.$$
We say that $A$ preserves the timelike orientation if $A$ carries a future-directed basis $B$ in other future-directed basis. The set $O_1^+(3)$ is a group with two components, one of them is $O_1^+(3)\cap SO_1(3)$. This proves that $O_1(3)$ has exactly four components. This contrast with the   isometries $O(3)$ of  Euclidean space $\e^3$, which  has exactly two connected components, being one of them, the special orthogonal group $SO(3)$.

The \emph{special Lorentz ortocrone} group  is the set $O_1^{++}(3)=SO_1(3)\cap O_1^+(3)$. This set is a group and $I\in O_1^{++}(3)$. From a topological viewpoint,   $O_1^{++}(3)$ is not a compact set, in contrast to $SO(3)\subset O(3)$, which is compact.

\begin{theorem}\label{iso} The connected components of   $O_1(3)$ are
$O_1^{++}(3)$ and
\begin{eqnarray*}
O_1^{+-}(3)&=&\{A\in SO_1(3): a_{33}<0\}\\
O_1^{-+}(3)&=&\{A\in O_1^+(3):\mbox{det}(A)=-1\}\\
O_1^{--}(3)&=&\{A\in O_1(3):{\mbox det}(A)=-1, a_{33}<0\}
\end{eqnarray*}
\end{theorem}
If we denote by  $T_1$ and  $T_2$ the isometries of $\e_1^3$ defined by
 $$T_1=\left(\begin{array}{ccc}
 1&0&0\\ 0&-1&0\\ 0&0&1\end{array}\right),\ \ T_2=\left(\begin{array}{ccc}
 1&0&0\\ 0&1&0\\ 0&0&-1\end{array}\right),$$
  then the three last components that appear in Th. \ref{iso} correspond, respectively, with
$T_2\cdot T_1\cdot O_1^{++}(3)$, $T_1\cdot O_1^{++}(3)$ and  $T_2\cdot O_1^{++}(3)$.

In order to clarify   why there appear  four connected components, we  compute the isometries in the two-dimensional case, that is, in  $\e_1^2$.    Let $A=\left(\begin{array}{cc}a&b\\ c&d\end{array}\right)$. Then $A\in O_1(2)$ if and only if $G=A^tGA$, where $G=\left(\begin{array}{cc}1&0\\ 0&-1\end{array}\right)$. This leads to the next three equations
\begin{equation}\label{eq1-12}
a^2-c^2=1,\ b^2-d^2=-1,\ ab-cd=0.\end{equation}
 If we compare with the isometries $O(2)$ of the Euclidean plane $\e^2$, the difference lies in the first two equations, because in $\e^2$ change by $a^2+b^2=1$ and $c^2+d^2=1$. The solution is $a=\cos\theta$ and $b=\sin\theta$.

However in \eqref{eq1-12} we have   $a^2-b^2=1$, which describes a hyperbola, in particular, there are two branches (two components).  The same occurs with the third equation and combining all the cases, we obtain the four components. Exactly, we have:
\begin{enumerate}
\item There exists $t$ such that $a=\cosh(t)$ and $c=\sinh(t)$. From  $b^2-d^2=-1$, it appears two cases again:
\begin{enumerate}
\item There exists $s$ such that $b=\sinh(s)$ and $d=\cosh(s)$. Using the third equation in \eqref{eq1-12}, we conclude   $s=t$.
\item There exists  $s$ such that  $b=\sinh(s)$ and $d=-\cosh(s)$. Now we have   $s=-t$.
\end{enumerate}
\item There exists $t$ such that $a=-\cosh(t)$ and $c=\sinh(t)$. Equation $b^2-d^2=-1$ in \eqref{eq1-12} yields two possibilities:
\begin{enumerate}
\item There exists  $s$ such that $b=\sinh(s)$ and $d=\cosh(s)$. The third equation of \eqref{eq1-12} concludes that  $s=-t$.
\item There exists  $s$ such that $b=\sinh(s)$ and $d=-\cosh(s)$. From $ab-cd=0$, we have  $s=t$.
\end{enumerate}
\end{enumerate}
As conclusion, we obtain  four kinds of isometries. In the same order that we have obtained them, they are the following:
$$\left(\begin{array}{cc}\cosh(t)&\sinh(t)\\ \sinh(t)&\cosh(t)\end{array}\right),\ \left(\begin{array}{cc}\cosh(t)&-\sinh(t)\\ \sinh(t)&-\cosh(t)\end{array}\right),$$
$$ \left(\begin{array}{cc}-\cosh(t)&-\sinh(t)\\ \sinh(t)&\cosh(t)\end{array}\right),\ \left(\begin{array}{cc}-\cosh(t)&\sinh(t)\\ \sinh(t)&-\cosh(t)\end{array}\right).$$
By using the notation as in Theorem  \ref{iso}, each one of the matrices that have appeared belong to
$O_1^{++}(2)$, $O_1^{--}(2)$, $O_1^{-+}(2)$ and $O_1^{+-}(2)$, respectively.

We end   the study of  isometries with the family of isometries that leave pointwise fixed a straight-line $L$. This kind of isometries are called  {\it boosts} of axis $L$ (\cite[p. 236]{on}) and they are the counterpart of the groups of rotations of $\e^3$. Depending on the causal character of  $L$, there are three types of such isometries.

\begin{enumerate}
\item The axis $L$ is timelike. Assume    $L=\mbox{Span}\{E_3\}$. Since the restriction of   the isometry to $L^{\bot}$ is an isometry in a vector plane which it is positive definite, the isometry is
$$A= \left(\begin{array}{ccc}
\cos\theta&-\sin\theta&0\\ \sin\theta&\cos\theta&0\\
0& 0&1\end{array}\right).$$
\item The axis $L$ is spacelike. Let  $L=\mbox{Span}\{E_1\}$. Now  $L^\bot$   is a Lorentzian plane and the restriction of the isometry to  $L^{\bot}$ belongs to $O_1^{++}(2)$. Then the isometry is
$$A=\left(\begin{array}{ccc}1&0&0\\
0&\cosh\theta&\sinh\theta\\
0&\sinh\theta&\cosh\theta\end{array}\right).$$
\item The axis $L$ is lightlike. Suppose    $L=\mbox{Span}\{E_2+E_3\}$. A straightforward computations gives
$$A= \left(
\begin{array}{ccc}
1&\theta &-\theta\\
-\theta &1-\frac{\theta^2}{2}&\frac{\theta^2}{2}\\
-\theta &-\frac{\theta^2}{2}&1+\frac{\theta^2}{2}
\end{array}
\right).$$
\end{enumerate}

The boosts allow to define a circle in $\e_1^3$.  In Euclidean space $\e^3$, there are different ways to define a circle. A first assumption is that the circle is included in a plane and that the curve is complete.  We have the next possibilities:
a) the set of points equidistant   from a given point; b) a curve with constant curvature; c) the orbit of a point under a group of rotations of $\e^3$.

In Lorentz-Minkowski space $\e_1^3$, we follow the last approximation but replacing rotations by boosts.  Let $L$ be a fixed straight-line of $\e_1^3$ and let  $G_L=\{\phi_\theta: \theta\in\r\}$  be the group of boost that leave pointwise fixed  $L$. A \emph{circle} is the orbit  $\{\phi_\theta( p_0):\phi_{\theta}\in G_L\}$ of a point $p_0\not\in L$, $p_0=(x_0,y_0,z_0)$. We distinguish three cases depending on the causal character of $L$.  After an isometry of $\e_1^3$, we have:
\begin{enumerate}
\item The axis $L$ is timelike. Consider  $L=\mbox{Span}\{E_3\}$.  Then the set $\{\phi_\theta(p_0):\theta\in\r\}$ is the Euclidean  circle that lies in the plane  of equation $z=z_0$ and radius $\sqrt{x_0^2+y_0^2}$.
\item The axis $L$ is spacelike. We take    $L=\mbox{Span}\{E_1\}$. Here we suppose $y_0^2-z_0^2\not=0$ since on the contrary, we get a straight-line. Then the orbit of  $p_0$ is a branch of the hyperbola  $y^2-z^2=y_0^2-z_0^2$  in the plane of equation  $x=x_0$. Depending if $y_0^2-z_0^2>0$ or $y_0^2-z_0^2<0$, we will have  four possibilities.
\item The axis $L$ is  lightlike. Assume that  $L=\mbox{Span}\{E_2+E_3\}$ and consider the plane $L^\bot=\mbox{Span}\{E_1,E_2+E_3\}$. The orbit of a point $p_0=(x_0,y_0,z_0)\not\in  L^{\bot}$ is a plane curve included in the plane $y-z=y_0-z_0$. If  $X=x_0+\theta(y_0-z_0)$ and $Y=y_0-x_0\theta-(y_0-z_0)\theta^2/2$, then the orbit of $p_0$ satisfies the equation
$$Y=\frac{-X^2+2y_0(y_0-z_0)+x_0^2}{2(z_0-y_0)}.$$
This means that the circle $\{\phi_\theta(p_0):\theta\in\r\}$ is a parabola.
\end{enumerate}

We point out that the  orbits are Euclidean circles, hyperbolas and parabolas only in the case that the axis of the group of boosts is one of the above ones. In general, they are affine ellipse, hyperbola or parabola, depending on the case. For example, we consider the rotations with respect to the timelike line $L=\mbox{Span}\{(0,1,2)\}$. Then $L^{\bot}=\mbox{Span}\{E_1,(0,2,1)/\sqrt{3}\}$. If $p=(1,0,0)\in L^{\bot}$, then
$$\phi_\theta(p)=\cos\theta E_1+\sin\theta\frac{1}{\sqrt{3}}(0,2,1),$$
which is an affine ellipse  included in $L^{\bot}$

\section{Curves in Minkowski space}\label{section2}

In this section we develop the theory of the Frenet trihedron for curves in  $\e_1^3$.  A (smooth) curve is a differentiable map $\alpha:I\subset\r\rightarrow\e_1^3$  where $I$ is an open interval. We also say that $\alpha$ is a parametrized curve. A curve is said to be regular if $\alpha'(t)\not=0$ for all $t\in I$. Here we do not use that $\e_1^3$ is a metric space, but that the codomain is $\r^3$, that is, a $3$-dimensional manifold. In other words, a regular curve is an immersion between  the (one-dimensional) manifold $I$ and the (three-dimensional) manifold $\r^3$.

\subsection{The local theory of curves}

Let $\alpha:I\rightarrow\e_1^3$ a regular curve. If $t\in I$, the tangent space $T_tI$ identifies with $\r$ and the differential map
 $(d\alpha)_t:T_t I\equiv\r\rightarrow T_{\alpha(t)}\e_1^3\equiv \r^3$  is
$$(d\alpha)_t(s)=\frac{d}{du}{\Big|}_{u=0}\alpha(t+su)=s\cdot\alpha'(t).$$
Thus the linear map $(d\alpha)_t$ is  a homothety from $\r$ to $\r^3$ given by $t\longmapsto s\cdot \alpha'(t)$. Here we identify $(d\alpha)_t$ by $\alpha'(t)$, or in other words, if $\partial/\partial t$ is the unit tangent vector on $T_tI$, then
$$(d\alpha)_t(\frac{\partial}{\partial t})=\alpha'(t).$$

We now endow $\r^3$ with the Lorentzian metric $\langle,\rangle$. On $I$ we consider the induced metric of $\e_1^3$ by the map $\alpha$ which converts
$$\alpha:(I,\alpha^{*}\langle,\rangle)\rightarrow \e_1^3=(\r^3,\langle,\rangle)$$
in an isometric immersion. The pullback metric $\alpha^{*}\langle,\rangle$ is  now
$$\alpha^*\langle,\rangle_t(m,n)=\langle (d\alpha)_t(m),(d\alpha)_t(n)\rangle=mn\langle\alpha'(t),\alpha'(t)\rangle,\ m,n\in\r,$$
or if we take the basis $\{\partial/\partial t\}$ in $T_tI$,
$$\alpha^{*}\langle,\rangle_t(\frac{\partial}{\partial t},\frac{\partial}{\partial t})=\langle\alpha'(t),\alpha'(t)\rangle.$$
In order to classify the manifold $(I,\alpha^{*}\langle,\rangle)$ and since $I$ is a one-dimensional manifold,  we need to know the sign of $\langle\alpha'(t),\alpha'(t)\rangle$. Thus
\begin{enumerate}
\item If $\langle\alpha'(t),\alpha'(t)\rangle>0$, $(I,\alpha^{*}\langle,\rangle)$ is a Riemannian manifold.
\item If $\langle\alpha'(t),\alpha'(t)\rangle<0$, $(I,\alpha^{*}\langle,\rangle)$ is a Lorentzian manifold, that is, the induced metric is non-degenerate with index $1$.
\item If $\langle\alpha'(t),\alpha'(t)\rangle=0$, $(I,\alpha^{*}\langle,\rangle)$ is a degenerate manifold.
\end{enumerate}
This classification justifies the following definition.

\begin{definition}\label{de-21} A curve $\alpha$  in  $\e_1^3$ is said spacelike (resp. timelike, lightlike) at  $t$ if $\alpha'(t)$ is a spacelike (resp. timelike, lightlike) vector. The curve \ $\alpha$ is   spacelike (resp. timelike, lightlike) if it is spacelike (resp. timelike, lightlike) for all $t\in I$.
\end{definition}
In particular, a timelike or a lightlike curve is regular. We point out that a curve in $\e_1^3$ may not be of one of the above types.   For example,   we consider the curve
$$\alpha:\r\rightarrow\e_1^3,\ \ \alpha(t)=(\cosh(t),\frac{t^2}{2},\sinh(t)).$$
Since $\alpha'(t)=(\sinh(t),t,\cosh(t)$, $\alpha$ is a regular curve. As $\langle\alpha'(t),\alpha'(t)\rangle=t^2-1$, then the curve is spacelike in  $(-\infty,-1)\cup(1,\infty)$, timelike in the interval $(-1,1)$ and lightlike in $\{-1,1\}$. Observe that $|\alpha'(\pm 1)|=0$, but $\alpha$ is regular at $t=\pm 1$.

However the spacelike (or timelike) condition is an open property, that is, if $\alpha$ is spacelike (or timelike) at  $t_0\in I$, there exists an interval  $(t_0-\delta,t_0+\delta)$ around $t_0$ where $\alpha$ is spacelike (or timelike):  if at $t_0\in I$
we have $\langle\alpha'(t_0),\alpha'(t_0)\rangle\not=0$, the continuity assures the existence of an interval around $t_0$ where  $\langle\alpha'(t),\alpha'(t)\rangle$ has the same sign than at $t=t_0$.

\begin{example}\label{curves1}  Consider  plane curves, that is, curves included in an affine plane of $\r^3$ and we study its causal character.
\begin{enumerate}
\item The straight-line  $\alpha(t)=p+t v$, $p,v\in\r^3$, $v\not=0$. This curve has the same causal character than the vector $v$.
\item The circle  $\alpha(t)=r(\cos t,\sin t,0)$ is a spacelike curve included in the spacelike plane of equation $z=0$.
\item The hyperbola  $\alpha(t)=r(0,\sinh t,\cosh t)$ is a spacelike curve in the timelike plane of equation $x=0$.
\item The hyperbola $\alpha(t)=r(0,\cosh t,\sinh t)$ is a timelike curve in the timelike plane of equation $x=0$.
\item The parabola  $\alpha(t)=(t,t^2,t^2)$ is a spacelike curve in the lightlike plane of equation $y-z=0$.
\end{enumerate}
\end{example}

\begin{example}
Consider spatial curves.
\begin{enumerate}
\item The helix  $\alpha(t)=(r\cos t,r\sin t, ht)$, $h\not=0$, of radius $r>0$ and pitch $2\pi h$. This curve is included in the cylinder of equation $x^2+y^2=r^2$. If $r^2>h^2$ (resp. $r^2-h^2<0$, $r^2=h^2$), $\alpha$ is a spacelike (resp. timelike, lightlike) curve.
\item Let $\alpha(t)=(ht,r\sinh t,r\cosh t)$, $h\not=0$, $r>0$. This curve is spacelike and included in the hyperbolic cylinder of equation $y^2-z^2=-r^2$.
\item Let $\alpha(t)=(ht,r\cosh t,r\sinh t)$, $h\not=0$, $r>0$. Here if $h^2-r^2>0$ (resp. $<0$, $=0$), the curve is spacelike (resp. timelike, lightlike). Moreover, the curve $\alpha$ is included in the hyperbolic cylinder of equation $y^2-z^2=r^2$.
    \end{enumerate}
\end{example}

It is well known that a regular curve in Euclidean space $\e^3$ is locally the graph (on a coordinate axis of $\r^3$) of two differentiable functions defined on a coordinate axis of $\r^3$. This is a consequence of the regularity of the curve and the inverse function theorem and it is not depend on the metric. If the curve is included in $\e_1^3$, the   causal character of the curve informs what is the axis where the above two functions are defined.

\begin{proposition}\label{pr-twof} Let $\alpha:I\rightarrow\e_1^3$ be a timelike (resp. lightlike) curve and $t_0\in I$. Then there exists $\epsilon>0$ and smooth functions $f,g:J\subset\r\rightarrow\r$ such that $t=\phi(s)$ and
$\beta(s)=\alpha(\phi(s))=(f(s),g(s),s)$.
\end{proposition}
\begin{proof}
Writing $\alpha(t)=(x(t),y(t),z(t))$, we know that $x'(t)^2+y'(t)^2-z'(t)^2\leq 0$. Then $z'(t_0)\not=0$. By the inverse function theorem, there exist $\delta,\epsilon>0$ such that $z:(t_0-\delta,t_0+\delta)\rightarrow (z(t_0)-\epsilon,z(t_0)+\epsilon)$ is a diffeomorphism. Denote $J=(z(t_0)-\epsilon,z(t_0)+\epsilon)$ and $\phi=z^{-1}$. Then the curve $\beta=\alpha\circ\phi$ satisfies
$$\alpha(\phi(s))=\beta(s)=((x\circ\phi)(s),(y\circ\phi)(s),s).$$
Take $f=x\circ \phi$ and $g=y\circ\phi$.
\end{proof}

In Euclidean plane $\e^2$ there is a rich theory of closed curves involving classical topics, as for example, the isoperimetric inequality, the four vertex theorem or the theorem of turning tangent (see \cite[sect. 1.7]{dc}). We see how the causal character of a curve in $\e_1^3$ imposes restrictions on planar closed curves.

A closed curve $\alpha:\r\rightarrow\e_1^3$ is a   parametrized curve  that is  periodic. If the curve is regular, there exists a minimum value $T>0$ such that $\alpha(t+T)=\alpha(t)$. In particular, the trace of $\alpha$ is a compact set.

\begin{theorem}  Let $\alpha$ be a closed regular curve in  $\e_1^3$ included in a plane  $P$.
 If $\alpha$ is spacelike, then $P$ is a spacelike plane.
\end{theorem}

\begin{proof} We distinguish cases according the causal character of $P$.
\begin{enumerate}
\item The plane $P$ is timelike. After a rigid motion of $\e_1^3$, we assume that $P$ is the plane of equation $x=0$. Then $\alpha(t)=(0,y(t),z(t))$. Since the function $y:\r\rightarrow\r$ is periodic, it attains a maximum at some point $t_0$. Then  $y'(t_0)=0$ and  so $\alpha'(t_0)=(0,0,z'(t_0))$. As $\alpha$ is regular, $z'(t_0)\not=0$ and this implies that $\alpha$ is timelike at $t=t_0$, a contradiction with the spacelike property of $\alpha$.
    \item The plane $P$ is lightlike. Suppose that $P$ is the plane of equation $y-z=0$. Then  $\alpha(t)=(x(t),y(t),y(t))$. Let $t_0$ be the maximum of the function
 $x(t)$. This implies $x'(t_0)=0$ and so, $\alpha'(t_0)=(0,y'(t_0),y'(t_0))$. Again, $y'(t_0)\not=0$ by regularity, but his implies that $\alpha'(t_0)$ is lightlike, a contradiction. \end{enumerate}
 From the above discussion, we conclude that $P$ is necessarily a spacelike plane.
\end{proof}
Therefore, and after an isometry, a plane closed spacelike curve is a closed curve in an Euclidean plane $\e^2$. This means that the theory of plane closed spacelike curves is the same than in Euclidean plane.

With the same arguments, we have:

\begin{theorem}
There are not closed curves in $\e_1^3$ that are timelike or lightlike.
\end{theorem}

\begin{proof} By contradiction, assume that the curve is closed. Using the same notation and because the function $z=z(t)$ is periodic,   there exists $t=t_0$ such that $z'(t_0)=0$. Then
$$\langle\alpha'(t_0),\alpha'(t_0)\rangle=x'(t_0)^2+y'(t_0)\geq 0.$$
This is a contradiction if  $\alpha$ is timelike. If $\alpha$ is lightlike, then $x'(t_0)=y'(t_0)=0$ and so, $\alpha'(t_0)=0$. In particular, $\alpha$ is not regular at $t=t_0$, a contradiction.
\end{proof}

Again, we conclude that there is no a theory of closed timelike or lightlike curves.

In Euclidean space, a regular curve $\alpha$ can be parameterized by the arc-length, that is, $\langle\alpha'(s),\alpha'(s)\rangle =1$ for all $s$. The same result holds for spacelike and timelike curves of $\e_1^3$.

\begin{proposition} Let $\alpha:I\rightarrow\e_1^3$ be a  spacelike or timelike curve. Given  $t_0\in I$, there is  $\delta,\epsilon>0$ and a diffeomorphism  $\phi:(-\epsilon,\epsilon)\rightarrow (t_0-\delta,t_0+\delta)$ such that the curve
$\beta:(-\epsilon,\epsilon)\rightarrow\e_1^3$ given by $\beta=\alpha\circ\phi$ satisfies $|\beta'(s)|=1$ for all $s\in (-\epsilon,\epsilon)$.
\end{proposition}

\begin{proof} We do the proof for timelike curves.  Define the function
$$S:I\rightarrow\r,\ \ S(t)=\int_{t_0}^t|\alpha'(u)|\ du.$$
 Since  $S'(t_0)>0$, the function $S$ is a local diffeomorphism around
$t=t_0$.  Because $S(t_0)=0$,  there exist $\delta,\epsilon>0$ such that  $S:(t_0-\delta,t_0+\delta)\rightarrow (-\epsilon,\epsilon)$ is a diffeomorphism. The map that we are looking for is $\phi=S^{-1}$.
\end{proof}

For a lightlike curve,   there is not sense reparametrize by the arc-length. However a differentiation of
$\langle\alpha'(t),\alpha'(t)\rangle=0$ gives  $\langle\alpha''(t),\alpha'(t)\rangle=0$. By Prop. \ref{pr-2-21}, $\mbox{Span}\{\alpha'(t)\}^\bot$ is a lightlike plane. We distinguish the next cases:
\begin{enumerate}
\item If $\alpha''(t)$ is lightlike, then $\alpha''(t)$ is proportional to $\alpha'(t)$ by Prop. \ref{pr-1-1}. If this holds for all $t$, then an easy integration yields
    $$\alpha(t)=e^t a+b,\ \ a,b\in\r^3, \ \langle a,a\rangle=0.$$
    This means that $\alpha$ is a parametrization of a (lightlike) straight-line.
    \item If $\alpha''(t)$ is spacelike, then we can parametrize $\alpha$ to get  $|\alpha''(t)|=1$. This is given in the next result.
    \end{enumerate}

\begin{lemma} Let $\alpha:I\rightarrow\e_1^3$ be a lightlike curve such that the trace of $\alpha$ is not a straight-line. There exists a reparametrization of $\alpha$ given by  $\beta(s)=\alpha(\phi(s))$   such that $|\beta''(s)|=1$. We say that $\alpha$ is   pseudo-parametrized by arc length.
\end{lemma}
\begin{proof}
We write  $\beta(s)=\alpha(\phi(s))$. Then
$$\beta''(s)=\phi''(s)\alpha'(t)+\phi'(s)^2\alpha''(t)\Rightarrow \langle\beta''(s),\beta''(s)\rangle=\phi'(s)^4|\alpha''(t)|^2.$$
It suffices by defining $\phi$ as the solution of the differential equation
$$\phi'(s)=\frac{1}{\sqrt{|\alpha''(\phi(s))|}}.$$
\end{proof}

\begin{remark} If $\alpha=\alpha(t)$ is a regular curve and $\beta=\alpha\circ\phi$ is a reparametrization of $\alpha$, the causal character of $\alpha$ and $\beta$ coincides.
\end{remark}

\subsection{Curvature and torsion. Frenet equations}
We  want to assign a basis of $\e_1^3$ for each point of a regular curve $\alpha(s)$ whose variation describes the geometry of the curve. This will be given by the Frenet trihedron $\{\t(s),\n(s),\b(s)\}$. In Euclidean space, the Frenet frame is a positively oriented  orthonormal basis, with $\b=\t\times\n$.

We assume that the curve is parameterized by the arc length or the pseudo arc length. Recall that the vector $\t(s)$ is the velocity of $\alpha$. In Minkowski space there appear some problems.
\begin{enumerate}
\item If the curve is lightlike, $\t(s)$ is a lightlike vector and so, $\{\t,\n,\b\}$ is not an orthonormal basis. In this situation, we will use the concept of null frame.
    \item Assume that $\{\t,\n,\b\}$ is an orthonormal basis of $\e_1^3$. The binormal vector $\b$ will be always defined as $\t\times\n$. Now the basis $\{\t,\n,\b\}$ is not necessarily positively oriented, as for example, if $\t,\n$ are spacelike vectors.
        \item It would be desirable that in the case that $\{\t,\n,\b\}$ is an orthonormal basis, this basis is future-directed. This can not assure a priori. Even in the case that $\alpha$ is a timelike curve, $\alpha'(s)=\t(s)$ could not be future-directed.
    \end{enumerate}

The simplest example of a curve is a straight-line. If $p\in\e_1^3$ and $v\not=0$ , the straight-line through the point  $p$ in the direction $v$ is parametrized by  $\alpha(s)=p+sv$. Then $\alpha''(s)=0$. In such a case,  we say that   the curvature   is  $0$.

Conversely, if  $\alpha$ is a regular curve that satisfies $\alpha''(s)=0$ for any  $s$, an integration  gives
 $\alpha(s)=p+sv$, for some $p,v\in\e_1^3$, $v\not=0$. This means that $\alpha$ parametrizes a straight-line through the point
$p$ along the direction given by  $v$. Let us observe that given a straight-line (as a set of $\e_1^3$), there are other parametrizations. For example,  $\alpha(s)=(s^3+s,0,0)$ is a parametrization of the straight-line $\mbox{Span}\{E_1\}$ where $\alpha''(s)\not=0$.

Consider $\alpha:I\rightarrow\e_1^3$ a regular curve  parametrized by arc length or by the pseudo arc length depending on the case. We call $$\t(s)=\alpha'(s)$$ the
 tangent vector at $s$. Since $\langle \t(s),\t(s)\rangle$ is constant, indeed, $1$, $-1$ or $0$, by differentiating with respect to $s$, we have   $\langle \t(s),\t'(s)\rangle=0$ and $\t'(s)$ is orthogonal to $\t(s)$. We shall restrict to curves such that $\t'(s)\not=0$ for all $s$ and that $\t'(s)$ is not proportional to $\t(s)$ for each $s$. This avoids that the curve is a straight-line.

We distinguish three cases depending on the causal character of $\t(s)$.

\begin{center}{\bf The curve is timelike.}\end{center}
As $\t(s)$ is a timelike vector,   Prop. \ref{pr-2-21} asserts that $\t'(s)$ is a spacelike vector. Then  $\t'(s)\not=0$ is a spacelike vector  linearly independent with  $\t(s)$.
We define the curvature of $\alpha$ at $s$ as
$$\kappa(s)=|\t'(s)|.$$
The normal vector $\n(s)$ is defined by
$$\n(s)=\frac{\t'(s)}{\kappa(s)}\Rightarrow \t'(s)=\kappa(s)\n(s).$$
Moreover $\kappa(s)=\langle \t'(s),\n(s)\rangle$. We define the binormal vector $\b(s)$ as
$$\b(s)=\t(s)\times\n(s).$$
 The vector $\b(s)$ is unitary and spacelike. For each  $s$, $\{\t(s),\n(s),\b(s)\}$ is an orthonormal basis of $\e_1^3$ which is called   the Frenet trihedron of  $\alpha$ at $s$. The basis $\{\t,\n,\b\}$ is positively oriented because $\mbox{det}(\t,\n,\b)=\langle\t\times\n,\b\rangle=\langle\b,\b\rangle=1>0$.

 We define the torsion $\tau$ of  $\alpha$ at $s$ as
$$\tau(s)=\langle\n'(s),\b(s)\rangle.$$
By differentiation each one of the vector functions of the Frenet trihedron and writing in coordinates  with the same Frenet basis, we obtain the  Frenet equations (or Frenet formula), namely,

\begin{equation}\label{frenet1}
\left(\begin{array}{c}
 \t' \\   \n' \\ \b'
\end{array} \right)=
 \left(\begin{array}{ccc}
  0 & \kappa & 0 \\
 \kappa& 0 & \tau \\
 0 & -\tau & 0
\end{array}\right)\left(\begin{array}{c} \t \\ \n \\ \b
\end{array} \right).
\end{equation}

\begin{center}{\bf The curve is spacelike.}\end{center}

 Since $\t'(s)$ is orthogonal to the spacelike vector $\t(s)$, $\t'(s)$ may be spacelike, timelike or lightlike  by Prop. \ref{pr-2-21}. We analyse the three cases.
 \begin{enumerate}
 \item The vector   $\t'(s)$ is spacelike. Again, we write the curvature  $\kappa(s)=|\t'(s)|$, $\n(s)=\t'(s)/\kappa(s)$ and
  $\b(s)=\t(s)\times\n(s)$. The vectors $\n$ and $\b$ are called the normal vector and the binormal vector respectively. Here $\b(s)$ is a timelike vector.  The Frenet equations are
\begin{equation}\label{frenet2} \left(\begin{array}{c}
     \t' \\
     \n' \\
     \b'
   \end{array} \right)=
   \left(\begin{array}{ccc}
             0 & \kappa & 0 \\
             -\kappa& 0 & \tau \\
             0 & \tau & 0
   \end{array}\right)\left(\begin{array}{c} \t\\ \n \\ \b   \end{array} \right).
   \end{equation}
 The torsion   of $\alpha$ is $\tau=-\langle\n',\b\rangle$. Here the basis $\{\t,\n,\b\}$ is negatively oriented because $\mbox{det}(\t,\n,\b)=
 \langle\t\times\n,\b\rangle=\langle\b,\b\rangle=-1<0$.

\item The vector  $\t'(s)$ is timelike. The curvature is
$$\kappa(s)=| \t'(s)|=\sqrt{-\langle \t'(s),\t'(s)\rangle}$$
and the normal vector is  $\n(s)=\t'(s)/\kappa(s)$. The binormal vector is  $\b(s)=\t(s)\times\n(s)$, which is a spacelike vector. The Frenet equations are
\begin{equation}\label{frenet3} \left(\begin{array}{c}
     \t' \\
     \n' \\
     \b'
   \end{array} \right)=
   \left(\begin{array}{ccc}
             0 & \kappa & 0 \\
             \kappa& 0 & \tau \\
             0 & \tau & 0
   \end{array}\right)\left(\begin{array}{c} \t\\ \n \\ \b   \end{array} \right).
   \end{equation}
The torsion of $\alpha$ is  $\tau=\langle\n',\b\rangle$. The Frenet basis is now positively oriented.

\item The vector $\t'(s)$ is lightlike for all $s$. We define the normal vector as $\n(s)=\t'(s)$, which is linearly independent with $\t(s)$. Let $\b(s)$ be the unique lightlike vector such that $\langle \n(s),\b(s)\rangle=-1$ and it is
    orthogonal to  $\t$. The vector  $\b(s)$ is the binormal vector of $\alpha$ at  $s$.  The Frenet equations are
\begin{equation}\label{frenet4}\left(\begin{array}{c}
     \t' \\
     \n' \\
     \b'
   \end{array} \right)=
   \left(\begin{array}{ccc}
             0 & 1 & 0 \\
             0& \tau &0 \\
             1 & 0&-\tau
 \end{array}\right)\left(\begin{array}{c} \t\\ \n \\ \b   \end{array} \right).
 \end{equation}
 The function $\tau$ is called the {\it pseudo-torsion}   of  $\alpha$ and it is obtained by $\tau=-\langle\n',\b\rangle$. There is not a definition of the curvature of $\alpha$.   Moreover, $\{\t,\n,\b\}$ is not an orthonormal basis of $\e_1^3$ since $\n$ and $\b$ are lightlike. Let us observe that $\{\t,\n,\b\}$ is a null frame and that a priori we do not know if it is positively oriented.
 \end{enumerate}

\begin{remark} We know that $\t'(s)$ is spacelike, timelike or lightlike, but it may occur that this does not hold in all the interval $I$, that is, the causal character of $\t'(s)$ may be change in $I$. Thus, in the above discussion of cases, we have assumed that the causal character of $\alpha''(s)$ is the same in $I$. As an example, for $s\in(1,\infty)$, let
$$\alpha(s)=\left(\cos(s)+s\sin(s),\sin(s)-s\cos(s),\frac12\left(s\sqrt{s^2-1}-\log(s+\sqrt{s^2-1})\right)\right).$$
Then
$$\t(s)=(s\cos(s),s\sin(s),\sqrt{s^2-1})$$
$$\t'(s)=\left(\cos(s)-s\sin(s),\sin(s)+s\cos(s),\frac{s}{\sqrt{s^2-1}}\right).$$
The causal character of $\t'(s)$ is given by the sign of $s^4-s^2-1$ since
$$\langle\alpha''(s),\alpha''(s)\rangle=\frac{s^4-s^2-1}{s^2-1}.$$
Thus $\t'(s)$ is spacelike if $s>\sqrt{1+\sqrt{5}}/\sqrt{2}$ and timelike if $1<s<\sqrt{1+\sqrt{5}}/\sqrt{2}$.
In both cases, the curvature and the torsion are
$$\kappa(s)=\sqrt{\frac{|s^4-s^2-1|}{s^2-1}},\ \tau(s)=\frac{s^6-2s^4-2s^2+2}{(s^4-s^2-1)\sqrt{s^2-1}}.$$
\end{remark}

\begin{center}{\bf The curve is lightlike.}\end{center}

Let $\alpha$ be a lightlike curve parametrized by the pseudo arc length. The tangent vector is $\t(s)=\alpha'(s)$. Define the normal vector as $\n(s)=\t'(s)$,  which it is a unit spacelike vector. The binormal vector is the unique lightlike vector is
orthogonal to $\n(s)$ such that  $\langle \t(s),\b(s)\rangle=-1$. Thus $\{\t,\n,\b\}$ is a null frame of $\e_1^3$. The Frenet equations are:
\begin{equation}\label{frenet5} \left(\begin{array}{c}
     \t' \\
     \n' \\
     \b'
   \end{array} \right)=
   \left(\begin{array}{ccc}
             0 & 1 & 0 \\
             \tau&0  &1 \\
             0 & \tau&0
 \end{array}\right)\left(\begin{array}{c} \t\\ \n \\ \b   \end{array} \right).
 \end{equation}
 The {\it pseudo-torsion} of $\alpha$ is $\tau=-\langle\n',\b\rangle$. As in the case that $\alpha$ is spacelike with $\t'$ lightlike, we do not define the curvature of the curve. We point out the reader that for lightlike curves there exists a variety of possibilities of concepts  where not all authors coincide. This is essentially due to different the possibilities in the choice of the Frenet frame. We refer \cite{bo,fgl,gr,in}.

\begin{example}\label{ex-21}
\begin{enumerate}
\item Let $\alpha(s)=r(\cos(s/r),\sin(s/r),0)$. Then
$$\t(s)=\left(-\sin(\frac{s}{r}),\cos(\frac{s}{r}),0\right),\ \ \t'(s)=\frac{1}{r}\left(-\cos(\frac{s}{r}),-\sin(\frac{s}{r}),0\right).$$
Then $\kappa=1/r$, and
$$\n(s)=\left(-\cos(\frac{s}{r}),-\sin(\frac{s}{r}),0\right),\ \ \b(s)=(0,0,-1).$$
As $\b'=0$, $\tau=0$. This basis is not positively oriented nor future directed.
\item Let $\alpha(s)= r(0,\sinh(s/r),\cosh(s/r))$. Then
$$\t(s)=\left(0,\cosh(\frac{s}{r}),\sinh(\frac{s}{r})\right),\ \ \t'(s)=\frac{1}{r}\left(0,\sinh(\frac{s}{r}),\cosh(\frac{s}{r})\right).$$
Then $\kappa=1/r$. Moreover,
$$\n(s)=\left(0,\sinh(\frac{s}{r}),\cosh(\frac{s}{r})\right),\ \  \b(s)=(1,0,0).$$
Here $\tau=0$.
\item Let $\alpha(s)= r(0,\cosh(s/r),\sinh(s/r))$. Then
$$\t(s)=\left(0,\sinh(\frac{s}{r}),\cosh(\frac{s}{r})\right),\ \ \t'(s)=\frac{1}{r}\left(0,\cosh(\frac{s}{r}),\sinh(\frac{s}{r})\right).$$
Here $\kappa=1/r$. Moreover
$$\n(s)=(0,\cosh(\frac{s}{r}),\sinh(\frac{s}{r})),\ \ \b(s)=(-1,0,0).$$
Again, $\tau=0$.
\item Let $\alpha(s)=(hs/\sqrt{r^2-h^2},r\cosh(s/\sqrt{r^2-h^2}),r\sinh(s/\sqrt{r^2-h^2}))$, where $r^2-h^2>0$. Then
    $$\t(s)=\frac{1}{\sqrt{r^2-h^2}}\left(h,r\sinh(\frac{s}{\sqrt{r^2-h^2}}),r\cosh(\frac{s}{\sqrt{r^2-h^2}})\right).$$
    This vector is timelike and future-directed. We have
    $$\t'(s)=\frac{r}{r^2-h^2}\left(0,\cosh(\frac{s}{\sqrt{r^2-h^2}}),\sinh(\frac{s}{\sqrt{r^2-h^2}})\right)$$
    $$\kappa=\frac{r}{r^2-h^2}.$$
Hence
$$\n(s)=\left(0,\cosh(\frac{s}{\sqrt{r^2-h^2}}),\sinh(\frac{s}{\sqrt{r^2-h^2}})\right)$$
$$\b(s)=\frac{1}{\sqrt{r^2-h^2}}\left(-r,-h\sinh(\frac{s}{\sqrt{r^2-h^2}}),-h\cosh(\frac{s}{\sqrt{r^2-h^2}})\right).$$
The torsion is $\tau=h/(r^2-h^2)$.
\item Let $\alpha(s)=r(s/r,(s/r)^2,(s/r)^2)$. This curve is spacelike with
$$\t(s)=\left(1,\frac{2s}{r},\frac{2s}{r}\right),\ \ \t'(s)=\n(s)=\left(0,\frac{2}{r},\frac{2}{r}\right).$$
Then $\b(s)=(0,-r/4,r/4)$ and $\tau=0$.
\item Consider a curve $\alpha$ constructed by the boosts about the lightlike axis $\mbox{Span}\{(0,1,1)\}$. Take the orbit $\beta$ of the point $(0,1,-1)$. Then $\beta(s)=(2s,1-s^2,-1-s^2)$. Hence $\beta'(s)=(2,-2s,-2s)$ and $\beta$ is spacelike. As $|\beta'(s)|=2$, we reparametrize by the pseudo  arc length changing $s$ by $s/2$. Thus let $\alpha(s)=(s,1-s^2/4,-1-s^2/4)$. Then
$\t(s)=(1,-s/2,-s/2)$ and $\t'(s)=(0,-1/2,-1/2)$. This means that the curve is spacelike with lightlike normal vector. Thus
$$\n(s)=\left(0,-\frac{1}{2},-\frac{1}{2}\right),\ \ \b(s)=\left(s,1-\frac{s^2}{4},-1-\frac{s^2}{4}\right).$$
Then $\tau=0$. We point out that $\alpha$ is included in the plane of equation $y-z=2$.
\item Consider the lightlike curve
$$\alpha(s)=\frac{1}{r^2}\left(\cosh(rs), r s,\sinh(rs)\right).$$
Then
$$\t(s)=\frac{1}{r}\left(\sinh(rs),1,\cosh(rs)\right),\
\n(s)=\t'(s)=(\cosh(rs),0,\sinh(rs)).$$
Hence that $\alpha$ is pseudo arc length. Then
$$\b(s)=\frac{r}{2}\left(\sinh(rs),-1,\cosh(rs)\right),\ \n'(s)=r(\sinh(rs),0,\cosh(rs)).$$
The pseudo-torsion is   $\tau=-r^2/2$.
\item Define
$$\alpha(s)=\frac{1}{r^2}\left(\cos(rs), \sin(rs),r s\right).$$
Then
$$\t(s)=\frac{1}{r}\left(-\sin(rs), \cos(rs),1\right),\ \t'(s)=\left(
-\cos(rs),-\sin(rs),0\right).$$
This curve is lightlike and parametrized by the pseudo arc length. Here
$$\b(s)=\frac{r}{2}\left(\sin(rs),-\cos(rs),1\right).$$
We deduce that $\tau=-r^2/2$.
\end{enumerate}

\end{example}

Timelike curves and spacelike curves with spacelike or timelike normal vector are called {\it Frenet curves}. In this case, the Frenet equations write in a unified way. If $\langle\t,\t\rangle=\epsilon$ and $\langle\n,\n\rangle=\delta$, then
 $$ \left(\begin{array}{c}
     \t' \\
     \n' \\
     \b'
   \end{array} \right)=
   \left(\begin{array}{ccc}
             0 & \kappa & 0 \\
            -\delta \kappa& 0 & \tau \\
             0 & \epsilon \tau & 0
   \end{array}\right)\left(\begin{array}{c} \t\\ \n \\ \b   \end{array} \right). $$
   We point out that for these curves, the curvature $\kappa$ is defined as the function such that $\t'$ is proportional to $\n$. Similarly, the torsion $\tau$ is defined as the third coordinate of $\n'$ with respect to the Frenet basis.

For  spacelike curves with lightlike normal vector or  lightlike curves, the Frenet equations write as follows: let $\langle\t,\t\rangle=\epsilon$, $\langle\n,\n\rangle=\delta$ where $\epsilon,\delta\in\{0,1\}$ and $\epsilon\not=\delta$. Then
$$ \left(\begin{array}{c}
     \t' \\
     \n' \\
     \b'
   \end{array} \right)=
   \left(\begin{array}{ccc}
             0 & 1 & 0 \\
            \delta  \tau&\epsilon\tau  &\delta \\
             \epsilon &\delta \tau&-\epsilon\tau
 \end{array}\right)\left(\begin{array}{c} \t\\ \n \\ \b   \end{array} \right). $$

The torsion is
$$\tau(s)=-\epsilon\delta\langle\n'(s),\b(s)\rangle.$$
In Euclidean space, the torsion measures how far is a curve to be planar in such way that a curve is planar if and only if its torsion is zero. In Minkowski space, we obtain the same result for Frenet curves.

\begin{theorem}\label{torsion}
Let $\alpha:I\rightarrow\e_1^3$ be a Frenet curve parametrized by the arc length. Then $\alpha$ is included in an affine plane if and only if the torsion vanishes.
\end{theorem}
The proof is the same and we omit it. However, there are more curves to consider.

\begin{theorem} Let $\alpha$ be a  spacelike curve with lightlike vector normal or a lightlike curve.
\begin{enumerate}
\item If the pseudo-torsion is zero, then the curve is included in a plane.
\item If a lightlike curve  is included in a plane, then it is a straight-line. There are spacelike plane curves with lightlike vector normal with non-zero pseudo-torsion.
\end{enumerate}
\end{theorem}

\begin{proof} For the first part, we do the arguments assuming that $\alpha$ is a  spacelike curve with lightlike vector normal (analogous if $\alpha$ is lightlike). If $\tau=0$, then $\n'=0$ and so, $\n(s)=v\in\e_1^3$ for all $s$. Let $s_0\in I$ and define the function $f(s)=\langle\alpha(s)-\alpha(s_0),v\rangle$. Then $f(s_0)=0$ and $f'(s)=\langle\t(s),v\rangle=0$. This implies that $f(s)$ is a constant point,  proving that the curve is planar.

Assume that $\alpha$ is a lightlike curve included in a plane. Then this plane must be timelike or lightlike. In the first case, there are only two linearly independent lightlike directions, and in the second one, there is only one. Thus $\t(s)$ is proportional to a fixed direction, proving that the curve is a straight-line.
\end{proof}

The reverse of the above theorem is false for spacelike curves with lightlike vector normal. An example is the following.   Let $\alpha(s)=(s,s^3/3,s^3/3)$, $s>0$, which is included in the plane $y-z=0$. To compute the Frenet frame, we have
$\t(s)=(1,s^2,s^2)$ and $\n(s)=\t'(s)=(0,2s,2s)$. As $\n'(s)=(0,2,2)$ and $\b(s)=(s/2,-1/(4s),1/(4s))$,  the pseudo-torsion is $\tau(s)=-\langle \n'(s),\b(s)\rangle=1/s$.

The rest of this subsection is devoted to prove the invariance of $\kappa$ and $\tau$ by rigid motions of $\e_1^3$ and  the theorem of existence and uniqueness. As in \cite[Ex. 6, p. 23]{dc}, we have:

\begin{theorem}\label{rigi} For a Frenet curve, the curvature is invariant under a rigid motion and the torsion is invariant up a sign depending if the motion is direct or inverse. In the case that the curve is spacelike with lightlike normal vector or that it is lightlike, then the pseudo-torsion is invariant by rigid motions.
\end{theorem}
\begin{proof}
The proof is the same than in Euclidean case. Let $Mx=Ax+b$ be a rigid motion, $A\in O_1(3)$, $b\in\e^3_1$, and  $\beta=M\circ\alpha$. Suppose that $\alpha$ is Frenet curve. Then the relation between the Frenet basis is
$$\t_\beta=A\t_\alpha,\ \n_\beta=A\n_\alpha,\ \b_\beta=\pm A\b_\alpha$$
depending of the sign of $\mbox{det}(A)$. This proves that $\kappa_\beta=\kappa_\alpha$ and $\tau_\beta=\pm\tau_\alpha$.

Assume now that $\alpha$ is a spacelike curve with lightlike normal vector. Then $\t_\beta=A\t_\alpha$ and
$\t'_\beta=A\t'_\alpha$ is lightlike, and so,  $\n_\beta=\t'_\beta=A\t'_\alpha=A\n_\alpha$. The vector $A\b_\alpha$ is a lightlike vector orthogonal to $\t_\beta$ and $\langle A\b_\alpha,\n_\beta\rangle=\langle A\b_\alpha,A\n_\alpha\rangle=-1$. Thus $\b_\beta=A\b_\alpha$. This implies
$$\tau_\beta=-\langle\n'_\beta,\b_\beta\rangle=-\langle A\n'_\alpha,A\b_\alpha\rangle=\tau_\alpha.$$

Consider now that $\alpha$ is lightlike. Then $\beta=M\circ\alpha$ is a lightlike curve parametrized by the pseudo arc length. Again $\t_\beta=A\t_\alpha$, $\n_\beta=A\n_\alpha$ and $\b_\beta=A\b_\alpha$, proving $\tau_\beta=\tau_\alpha$.
\end{proof}

In Euclidean space, the theorem of existence and uniqueness asserts that given two functions $\kappa>0$, $\tau$, there exists a unique curve, up a rigid motion, with curvature $\kappa$ and torsion $\tau$. In Minkowski space, the result of existence is the same but the uniqueness is not true by the causal character of the curve. For example, the curves $\alpha(s)=(\cos(s),\sin(s),0)$ and $\beta(s)=(0,\cosh(s),\sinh(s))$ have $\kappa=1$ and $\tau=0$ (see Example \ref{curves1}). However, there is not a rigid motion carrying $\alpha$ into $\beta$ because $\alpha$ is spacelike and $\beta$ is timelike. Even if both curves have the same causal character, we have to pay attention on the causal character of the other vectors of the trihedron of Frenet. For example, the curve $\gamma(s)=(0,\sinh(s),\cosh(s))$ has $\kappa=1$, $\tau=0$, but there is not a rigid motion between $\alpha$ and $\gamma$. Let us observe that for  $\alpha$, $\t$ and $\n$ are spacelike but $\gamma$ is a spacelike curve with timelike normal vector.

In Minkowski space, the result of existence is the same, although now the initial conditions will impose the causal character of the curve $\alpha$.  We analyse the different cases. First we consider Frenet curves.
 \begin{theorem}
 Let $\kappa(s)>0$ and $\tau(s)$, $s\in I$, two differentiable maps. Then there exists three different  regular parametrized curves $\alpha:I\rightarrow\e_1^3$, $\alpha=\alpha(s)$, with curvature $\kappa$ and torsion $\tau$.
\end{theorem}
By different we mean that there exist not a rigid motion carrying one in another one.
\begin{proof}
Let $s_0\in I$ and let $\{e_1,e_2,e_3\}$ be an orthonormal basis, which it will be the initial conditions of an ODE system. Depending on the causal character of the vectors $e_i$ we will obtain the different cases.
First assume that $e_1$ is timelike and that the basis is positively oriented. In such a case, we solve the next ODE system of $9$ equations
\begin{eqnarray*}
\t'(s)&=&\kappa(s)\n(s)\\
\n'(s)&=&\kappa(s)\t(s)+\tau(s)\b(s)\\
\b'(s)&=&-\tau(s)\n(s)
\end{eqnarray*}
with initial conditions
\begin{eqnarray*}\t(s_0)&=&e_1\\
\n(s_0)&=&e_2\\
\b(s_0)&=&e_3
\end{eqnarray*}
Let $\{\t,\n,\b\}$ be the unique solution and define
\begin{equation}\label{alfa2}
\alpha(s)=\int_{s_0}^s\t(u)\ du.
\end{equation}
We prove that this curve is timelike with curvature $\kappa$ and torsion $\tau$.
We first show that $\{\t(s),\n(s),\b(s)\}$ is an orthonormal basis  with the same causal properties that the initial basis $\{e_1,e_2,e_3\}$. Consider the ODE system:
\begin{eqnarray*}
\langle\t,\t\rangle'&=&2\kappa \langle\t,\n\rangle\\
\langle\n,\n\rangle'&=& 2\kappa\langle\t,\n\rangle+2\tau\langle\b,\n\rangle\\
\langle\b,\b\rangle'&=&-2\tau\langle \b,\n\rangle\\
\langle\t,\n\rangle'&=& \kappa\langle\n,\n\rangle+\kappa\langle\t,\t\rangle+\tau\langle\t,\b\rangle\\
\langle\t,\b\rangle'&=&\kappa\langle\n,\b\rangle-\tau\langle\n,\t\rangle\\
\langle\n,\b\rangle'&=&\kappa\langle\t,\b\rangle+\tau\langle\b,\b\rangle-\tau\langle\n,\n\rangle
\end{eqnarray*}
with initial conditions at $s=s_0$ given by  $(-1,1,1,0,0,0)$. On the other hand, the functions
$$f_1=-1,\ f_2=1,\ f_3=1,\ f_4=0,\ f_5=0,\ f_6=0$$
satisfy the same ODE system and initial conditions. By uniqueness,
$$-\langle\t,\t\rangle=\langle\n,\n\rangle=\langle\b,\b\rangle=1,\ \langle\t,\n\rangle=\langle\t,\b\rangle=\langle\n,\b\rangle=0.$$
This implies that $\{\t,\n,\b\}$ is an orthonormal basis of $\e_1^3$, where $\t$ is timelike. From \eqref{alfa2},  $\alpha'(s)=\t(s)$, and so $\alpha$ is a timelike curve parametrized by arc length. Now we are in conditions to prove that the curvature and torsion of $\alpha$ is $\kappa$ and $\tau$, respectively. The steps to follow are similar as in Euclidean space and we omit the details.

We can say a bit more. If the vector $e_1$ is future directed, then $\alpha$ is future-directed: since $\t(s)$ and $E_3$ are unit timelike vectors, then $\langle\t(s),E_3\rangle\geq 1$ in $I$ or $\langle\t(s),E_3\rangle\leq -1$ in $I$. As $\langle \t(s_0),E_3\rangle\leq -1$,    the same occurs for all $s\in I$  by connectedness.

If  we want to obtain  a spacelike curve with spacelike normal vector and curvature $\kappa$ and torsion $\tau$, consider the  initial conditions
$$\t(s_0)=e_1,\ \ \n(s_0)=e_2,\ \ \b(s_0)=e_3,$$
where $\{e_1,e_2,e_3\}$ is a  negatively oriented orthonormal basis and $e_3$ is timelike. The  ODE system that we solve is \eqref{frenet2}. Finally, if we are looking for a spacelike curve with timelike normal vector, the initial condition is a  positively oriented orthonormal basis  $\{e_1,e_2,e_3\}$, where $e_2$ is timelike and  the ODE system is \eqref{frenet3}.
\end{proof}

We analyse the remaining two cases. We prove the existence of  a spacelike curve with lightlike normal vector or a lightlike curve. Let $\tau:I\rightarrow \r$ be a smooth function and we ask for curves with the above causal character with pseudo-torsion $\tau$. The situation is similar than  Frenet curves and the solution depends on the initial conditions.

\begin{theorem} Let $\tau:I\rightarrow\r$ be a smooth function. Then there are a spacelike curve with lightlike normal vector and a lightlike curve with pseudo-torsion $\tau$.
\end{theorem}

\begin{proof} Let $\{e_1,e_2,e_3\}$ be a null frame of $\e_1^3$   such that $e_1$ is spacelike. We pose the ODE system \eqref{frenet4} with initial conditions
$$\t(s_0)=e_1,\ \n(s_0)=e_2,\ \b(s_0)=e_3.$$
Let $\{\t,\n,\b\}$ be the unique solution and define
\begin{equation}\label{ali}\alpha(s)=\int_{s_0}^s\t(u)\ du.
\end{equation}
We prove that $\alpha$ is a spacelike curve with lightlike normal vector. First we consider  the next ODE system of $6$ equations:
\begin{eqnarray*}
\langle\t,\t\rangle'&=&2  \langle\t,\n\rangle\\
\langle\n,\n\rangle'&=& 2\tau\langle\n,\n\rangle\\
\langle\b,\b\rangle'&=&2\langle\t,\b\rangle-2\tau\langle \b,\b\rangle\\
\langle\t,\n\rangle'&=&\langle\n,\n\rangle+\tau\langle\t,\n\rangle\\
\langle\t,\b\rangle'&=&\langle\n,\b\rangle+\langle\t,\t\rangle-\tau\langle\ \t,\b\rangle\\
\langle\n,\b\rangle'&=&\langle\t,\n\rangle
\end{eqnarray*}
with initial conditions $(1,0,0,0,0,-1)$. Since the functions $(1,0,0,0,0,-1)$ also are a set of solutions, by uniqueness, we have obtained the right solution. Thus $\{\t,\n,\b\}$ is a basis of $\e_1^3$ that satisfies the same properties than $\{e_1,e_2,e_3\}$.

From \eqref{ali}, $\alpha'(s)=\t(s)$ and it follows that  $\alpha$ is a spacelike curve. The arguments to end the proof are standard and we omit them.

If we want to obtain a lightlike curve with pseudo-torsion $\tau$, then we solve the ODE \eqref{frenet5} and we change the initial condition by
$$\t(s_0)=e_1,\ \n(s_0)=e_2,\ \b(s_0)=e_3,$$
where $\{e_1,e_2,e_3\}$ is a null frame with $e_2$ a unit spacelike vector.
\end{proof}

Once we have established the existence, the uniqueness does not hold as we have remarked previously, but it holds if the causal character  of the Frenet frame agree for both curves. To abbreviate the statements, we give the next definition.

\begin{definition} Let $\alpha,\beta:I\rightarrow\e_1^3$ be two curves parametrized by arc length or by the pseudo arc length.  We say that $\alpha$ and $\beta$ have the same causal character of the Frenet frame if $\t_\alpha$,  $\n_\alpha$ and  $\b_\alpha$  have the same causal character than $\t_\beta$,  $\n_\beta$ and $\b_\beta$, respectively.
\end{definition}

\begin{theorem} Let $\alpha,\beta:I\rightarrow\e_1^3$ be two regular curves that have the same causal character of the Frenet frame. If they have the same curvature and torsion, or they have same pseudo-torsion depending on the case, then there exists a rigid motion $M$ of $\e_1^3$ such that $\beta=M\circ\alpha$.
\end{theorem}

\begin{proof}  Let $s_0\in I$ and consider the isometry $A\in O_1(3)$ such that $A\t_\alpha(s_0)=\t_\beta(s_0)$, $A\n_\alpha(s_0)=\n_\beta(s_0)$ and $A\b_\alpha(s_0)=\b_\beta(s_0)$. If $b=\beta(s_0)-A\alpha(s_0)$, define the rigid motion $Mx=Ax+b$. We know by Th. \ref{rigi} that the curve $\gamma=M\circ\alpha$ satisfies the same ODE system of the Frenet equations than $\beta$. As the initial conditions coincide, then by uniqueness of ODE system, $\beta=\gamma$, proving the result. See \cite[p. 310]{dc}.
\end{proof}

Similarly as in Euclidean space, one can find the formula for the curvature and torsion function in the case that the curve is not parametrized by arc length. We only focus  on Frenet curves. Let $\alpha:I\rightarrow\e_1^3$ be a regular curve and  $\beta=\alpha\circ\phi$ be any parametrization by arc length. We define
$$\kappa_\alpha(t)=\kappa_\beta\circ\phi^{-1},\ \tau_\alpha=\tau_\beta\circ\phi^{-1}.$$
Assume that $\beta$ is a Frenet curvature.  The definition does not depend on the reparametrization, except perhaps a sign for the torsion. The proof is similar as in \cite[p. 25]{dc}. Then
 $$\kappa_{\alpha}(t)=\frac{|\alpha'(t)\times\alpha''(t)|}{|\alpha'(t)|^3},\ \
 \tau_\alpha(t)=-\epsilon\delta\frac{\mbox{det}(\alpha'(t),\alpha''(t),\alpha'''(t))}{  |\alpha'(t)\times\alpha''(t)|^2}.$$

We point out that given a regular curve $\alpha:I\rightarrow\e_1^3$ not necessarily parametrized by the arc length, we do not know a priori what is the causal character of the Frenet frame of its  parametrization by the arc length  $\beta=\alpha\circ\phi$. For example, consider the curve $\alpha(s)=(s^2,\sinh(s^2),\cosh(s^2))$, $s>0$. Then $\alpha'(s)$ is a spacelike and
$\alpha''(s)=2(1,\cosh(s^2)+s\sinh(s^2), \sinh(s^2)+s\cosh(s^2))$. Thus
$$\alpha''(s)\mbox { is}\left\{\begin{array}{ll}
\mbox{spacelike} & s\in (0,\sqrt{2})\\
\mbox{lightlike} & s =\sqrt{2}\\
\mbox{timelike} & s>\sqrt{2}\\
\end{array}\right.$$
However the parametrization by the arc length is $\beta(s)=(s/\sqrt{2},\sinh(s/\sqrt{2}),\cosh(s/\sqrt{2}))$, which it is spacelike.

\subsection{Curves in Lorentz-Minkowski plane}
We study plane curves in Minkowski space $\e_1^3$ giving a sign to the curvature $\kappa$.   A problem appears in a first moment showing a difference with the Euclidean context. We have two options. First, consider the two dimensional case of Lorentz-Minkowski space, that is, the Lorentz-Minkowski plane $\e_1^2$. The second possibility is to consider a curve of $\e_1^3$ included in an affine plane. In the latter case, there are three possibilities depending if the plane is    spacelike, timelike of lightlike. If the plane is spacelike, the theory corresponds to curves in a Riemannian surface. In this case, the plane is isometric to the Euclidean plane $\e^2$ and thus the theory is known; if  the plane is timelike, then it is isometric to $\e_1^2$, and we are in the first option; the case that  the plane is lightlike  is new and not covered by the Euclidean plane $\e^2$ or the Lorentzian plane $\e_1^2$.

We consider the first option. Denote $\e_1^2=(\r^2,(dx)^2-(dy)^2)$ the Lorentz-Minkowski plane.   We define the Frenet dihedron in such way that the curvature has a sign. Let $\alpha:I\rightarrow\e_1^2$ be a curve  parametrized by arc length. Define the tangent vector as
$$\t(s)=\alpha'(s).$$
We discard lightlike curves because in $\e_1^2$ there are two linearly independent directions of lightlike vectors. Thus $\t(s)$ would be proportional to a given direction, obtaining that the curve is a straight-line. In what follows, we suppose that $\alpha$ is spacelike or timelike.

The vector $\t'(s)$ is orthogonal to $\t(s)$. This means that $\t(s)$ and $\n(s)$ will have different causal character.

A new difference with Euclidean setting appears now. In $\e^2$, the unit normal $\n_e(s)$ is chosen so $\{\t(s),\n_e(s)\}$ is a positively oriented basis. Now in $\e_1^2$ we will choose again the Frenet frame as a positively oriented basis but the order of the vectors $\t$ and $\n$ is chosen under the condition the first vector is spacelike and the second one is timelike. In other words,  the cases are:
  \begin{enumerate}
  \item The curve is spacelike. Then define the normal vector $\n(s)$ such that $\{\t(s),\n(s)\}$ is positively oriented.
   \item The curve is timelike. Then define the normal vector $\n(s)$ such that $\{\n(s),\t(s)\}$ is positively oriented.
\end{enumerate}
Let $\langle\t,\t\rangle=\epsilon\in \{1,-1\}$ depending if the curve is spacelike or timelike. Then $\langle\n,\n\rangle=-\epsilon$. We define the curvature of $\alpha$ as the function $\kappa(s)$ such that
$$\t'(s)=\kappa(s)\n(s).$$
Thus
$$\kappa(s)=-\epsilon\langle\t'(s),\n(s)\rangle.$$
The Frenet equations are
$$\t'(s)=\kappa(s)\n(s)$$
$$\n'(s)=\kappa(s)\t(s).$$
As we expect, we observe that both equations coincide with the first two equations in \eqref{frenet1} and \eqref{frenet2}, neglecting the coordinate with respect to the binormal vector.

\begin{example}\label{ej-ab}
\begin{enumerate}
\item The set  $A=\{(x,y)\in\r^2: x^2-y^2=-r^2\}$ has two components
$$A^+=\{(x,y)\in A: y>0\},\ \ A^{-}=\{(x,y)\in A: y<0\}$$
   which parametrize as spacelike curves. For $A^+$, let $\alpha(s)=(r\sinh(s/r),r\cosh(s/r))$. Then
$$\t(s)=\left(\cosh(\frac{s}{r}),\sinh(\frac{s}{r})\right),\ \ \n(s)=\left(\sinh(\frac{s}{r}),\cosh(\frac{s}{r})\right).$$
As $\t'(s)=(1/r)(\sinh(s/r),\cosh(s/r))$, then $\kappa=1/r$.

For $A^-$, let $\beta(s)=(r\sinh(s/r),-\cosh(s/r))$. Then
$$\t(s)=(\cosh(s/r),-\sinh(s/r)),\ \n(s)=(-\sinh(s/r),\cosh(s/r)).$$
Hence we deduce $\kappa=-1/r$.
\item The set  $B=\{(x,y)\in\r^2: x^2-y^2=r^2\}$ has two components again, namely, $B^+=\{(x,y)\in B: x>0\}$ and $B^{-}=\{(x,y)\in B: x<0\}$. A parametrization of $B^+$ is
$\alpha(s)=(r\cosh(s/r),\sinh(s/r))$. Then
$$\t(s)=\left(\sinh(\frac{s}{r}),\cosh(\frac{s}{r})\right),\ \ \n(s)=\left(\cosh(\frac{s}{r}),\sinh(\frac{s}{r})\right).$$
Here   $\kappa(s)=1/r$.

   For $B^-$, let  $\beta(s)=(-r\cosh(s/r),\sinh(s/r))$. Then
   $$\t(s)=(-\sinh(\frac{s}{r}),\cosh(\frac{s}{r})),\ \ \n(s)=(\cosh(\frac{s}{r}),-\sinh(\frac{s}{r})).$$
   Thus $\kappa=-1/r$.
  \end{enumerate}
\end{example}

The existence and uniqueness result holds here and the proof is analogous that for Frenet curves in  the $3$-dimensional case. We point out two remarks. For the existence, and in Euclidean plane, it is possible to obtain a parametrization of the curve in terms of integrals of the curvature. Given a differentiable function  $\kappa$, let
\begin{equation}\label{teta}
\theta(s)=\int_{s_0}^s\kappa(t)\ dt.
\end{equation}
Define two  curves $\alpha$ and $\beta$ with curvature $\kappa$, where $\alpha$ is spacelike and $\beta$ is timelike:
$$\alpha(s)=\left(\int_{s_0}^s\cosh\theta(t)\ dt,\int_{s_0}^ s\sinh\theta(t)\ dt\right)$$
$$\beta(s)=\left(\int_{s_0}^s\sinh\theta(t)\ dt,\int_{s_0}^ s\cosh\theta(t)\ dt\right).$$
We extend to the Lorentzian space, the Euclidean result that asserts that the curvature of a plane curve is the variation of the angle between the tangent vector with a fix direction.

\begin{theorem}\label{2angulo} Let $\alpha:I\rightarrow\e_1^2$ be a timelike curve parametrized by arc length. Suppose that there exists  a unit timelike vector   $v\in\e_1^2$  such that $\t(s)$ and $v$ lies in the same timelike cone for all $s$. If $\theta$ is the angle between the tangent vector of $\alpha$ and $v$, then
$$\kappa(s)=\pm \theta'(s).$$
\end{theorem}

\begin{proof} We know that $-\cosh(\theta(s))=\langle\t(s),v\rangle$. By differentiating,
$$-\theta'(s)\sinh\theta(s)=\kappa(s)\langle\n(s),v\rangle.$$
As $v=-\langle v,\t(s)\rangle\t(s)+\langle v,\n(s)\rangle\n(s)$, then
$$-1=-\langle v,\t(s)\rangle^2+\langle v,\n(s)\rangle^2=-\cosh(\phi(s))^2+\langle v,\n(s)\rangle^2.$$
Then  $\langle v,\n(s)\rangle=\pm\sinh(\theta(s))$. Thus  $\theta'(s)=\pm\kappa(s)$.
\end{proof}
The same result holds for a spacelike curve $\alpha$, where $\t(s)$ makes constant angle with a unit spacelike vector in the same component of $\u_1^2$, that is, both belong to $\s_1^{1+}$ or $\s_1^{1-}$.

The second remark is about the uniqueness. Since the curvature has a sign, the curvature is only preserved by {\it direct} rigid motions. Exactly, if $\alpha:I\rightarrow\e_1^2$ is a spacelike or timelike curve and $Mx=Ax+b$ is a rigid motion of $\e_1^2$ with $A\in SO_1(2)$, then the relation between the Frenet frames of $\alpha$ and $\beta=M\circ\alpha$ is $\t_\beta=A\t_\alpha$ and $\n_\beta=A\n_\alpha$. Thus
$$\t_{\beta}'(s)= A\t_{\alpha}'(s)=A(\kappa_\alpha(s)\n_\alpha(s))=\kappa_\alpha(s)A\n_\alpha(s)=\kappa_\alpha(s)\n_\beta(s).$$
This proves $\kappa_\beta=\kappa_\alpha$.

Let $\alpha,\beta:I\rightarrow\e_1^2$ be two spacelike curves parametrized by arc length $s\in I$ and with the same curvature $\kappa$. We do a different proof that $\alpha$ and $\beta$ differ of a rigid motion without the use of the theory of ODE. We follow \cite[p. 20]{dc}  (the same arguments hold for a timelike curve by Rem. \ref{note}). Fix $s_0\in I$. The Frenet dihedrons $\{\t_\alpha(s_0),\n_\alpha(s_0)\}$ and $\{\t_\beta(s_0),\n_\beta(s_0)\}$ are two positively oriented basis.  Then  there exists an isometry $A\in SO_1(2)$ such that
\begin{eqnarray*}
A(\t_\alpha(s_0))&=&\t_\beta(s_0)\\
A(\n_\alpha(s_0))&=&\n_\beta(s_0).
\end{eqnarray*}
Let  $b=\beta(s_0)-A\alpha(s_0)$ and define the rigid motion  $Mx=Ax+b$.
Then the curve $\gamma=M\circ\alpha$ satisfies   $\kappa_\gamma=\kappa$ because $\mbox{det}(A)=1$. Denote $\{\t_\gamma,\n_\gamma\}$ the Frenet frame of  $\gamma$. Observe that the Frenet frames of $\gamma$ and $\beta$ coincide at $s=s_0$. Define
$$f:I\rightarrow\r,\ f(s)=| \t_\gamma(s)-\t_\beta(s)|^2-|\n_\gamma(s)-\n_\beta(s)|^2.$$
A differentiation of $f$, together the Frenet equations gives $f'(s)=0$ for all $s$, that is, $f$ is a constant function. As   $f(s_0)=0$, then $f(s)=0$ for all $s$. Now there is a great difference with the Euclidean case, because we can not assert that $|\t_\gamma-\t_\beta|=|\n_\gamma-\n_\beta|=0$ because the metric is not positive definite. However, expanding the equation $f(s)=0$, we obtain
$$\langle\n_\beta,\n_\gamma\rangle+1=\langle\t_\beta,\t_\gamma\rangle-1.$$
Observe that $\t_\beta(s)$ and $\t_\gamma(s)$ lie in the same component of unit spacelike vectors, that is, both lie in $\s_1^{1+}$ or both in $\s_1^{1-}$. This is because by connectedness and Rem. \ref{note}, $\langle\t_\beta(s),\t_\gamma(s)\rangle\geq 1$ or   $\langle\t_\beta(s),\t_\gamma(s)\rangle\leq -1$. As at $s=s_0$, $\t_\beta(s_0)=\t_\gamma(s_0)$, then  $\langle\t_\beta(s),\t_\gamma(s)\rangle\geq 1$.   For the vectors $\n_\beta(s)$ and $\n_\gamma(s)$ occurs the same. Using that $\n_\beta$ and $\n_\gamma$ lie in the same timelike cone,
$$0\geq \langle\n_\beta,\n_\gamma\rangle+1=\langle\t_\beta,\t_\gamma\rangle-1\geq 0.$$
This proves that $\langle\t_\beta,\t_\gamma\rangle= 1$ and so $\t_\beta=\t_\gamma$. Then $A\alpha'(s)=\beta'(s)$ for all $s$. By integrating, there exists $c\in\e_1^2$ such that
$\beta(s)=A\alpha(s)+c$.  Evaluating at $s=s_0$, we obtain $c=b$ and we conclude that $\beta=\gamma=M\circ\alpha$.

\begin{remark} The author has not been able to extend this proof in the $3$-dimensional case.
\end{remark}

We finish describing the curves in $\e_1^2$ of constant curvature. Assume that the curvature $\kappa$ is a constant $a\not=0$. Then
$$\theta(s)=\int_{s_0}^s a\ dt=as+b,\ b\in\r.$$
From \eqref{teta}, the next curves have curvature $a$:
\begin{enumerate}
\item The spacelike curve $$\alpha(s)=\frac{1}{a}\left(\sinh(as+b),\cosh(as+b)\right).$$
\item The timelike curve $$\beta(s)=\frac{1}{a}\left(\cosh(as+b),\sinh(as+b)\right).$$
 \end{enumerate}
 From the Euclidean viewpoint, both curves are Euclidean hyperbolas.

 As usual, a non-degenerate plane curve with constant curvature is called a circle. Recall the word `circle' appeared in  section \ref{section1} as the orbit of a point by the motion of a group of  boosts of $\e_1^3$. We now relate the notion of circle with plane curves with constant curvature.

 \begin{theorem} Let $\alpha:I\rightarrow\e_1^3$ be a Frenet curve included in a plane of $\e_1^3$. Then $\alpha$ is a circle if and only if the curvature is a non-zero constant and the torsion is $0$.
 \end{theorem}

 \begin{proof}
 We know by Th. \ref{torsion} that a Frenet curve included in a plane has $\tau=0$. If $\alpha$ is a Frenet curve, $\alpha$ can not be included in a lightlike plane since $\t$ or $\n$ is a timelike vector. Therefore, if $\alpha$ is a circle,  and  after a rigid motion, $\alpha$ is an Euclidean circle in the plane $z=0$ or a hyperbola en the plane $x=0$. By Ex. \ref{ej-ab}, these curves have constant curvature.

 The reverse statement is immediate because, after a rigid motion,  a planar Frenet curve is included in the plane $z=0$ or in the plane $x=0$. In both cases we know that a curve with constant curvature is a circle or a hyperbola, which both are invariant by a uniparametric group of boosts with timelike or spacelike axis.
 \end{proof}
By the way, in a spacelike or timelike plane, a circle is the set of equidistant points from a fix point $p_0$. If the plane is spacelike, we suppose that it is the $xy$-plane. Then the set of equidistant points from $p_0$ is an Euclidean circle, which it is a circle in Minkowski space. If the plane is timelike, we suppose that it is the $yz$-plane. Then the equidistant points from $p_0$ satisfies the equation $(y-y_0)^2-(z-z_0)^2=r^2$ or $(y-y_0)^2-(z-z_0)^2=-r^2$, which are circles in $\e_1^3$.

 Finally we focus on  spacelike curves included in a lightlike plane, in particular, the normal vector is lightlike.
  We have seen in Ex. \ref{ex-21}  circles obtained by boosts about lightlike axis which have not constant torsion. We study this type of curves  with non-zero constant pseudo-torsion. After a rigid motion, we suppose that the lightlike plane is the plane of equation $y-z=0$.

  \begin{theorem} Let $P$ be the lightlike plane  of equation $y-z=0$. The only spacelike curves  in $P$ with constant pseudo-torsion $\lambda\not=0$ are, up a change of parameter,
  $$\alpha(s)=\left(s+d,\frac{a}{\lambda^2}e^{\lambda s}+bs+c,\frac{a}{\lambda^2}e^{\lambda s}+bs+c\right),\ \ a,b,c,d\in\r.$$
 \end{theorem}

 \begin{proof} Let $\alpha(s)=(x(s),y(s),y(s))$. As $\alpha$ is parametrized by the arc length, then $x'(s)=\pm 1$. Up a constant and up a change of parameter, $x(s)=s$. Now
 $$\n(s)=\t'(s)=(0,y'',y''),\ \ \b(s)=\left(\frac{y'}{y''},\frac{-1+y'^2}{2y''},\frac{1+y'^2}{2y''}\right).$$
 Observe that $y''\not=0$ because on the contrary, $y(s)=as+b$, $a,b\in\r$, showing that $\alpha$ is straight-line. The computation of the pseudo-torsion gives $\tau=-\langle\n',\b\rangle=y'''/y''$, with $y''\not=0$. Since $y'''/y''=\lambda$, by solving this differential equation, we obtain the explicit parametrization of the curve claimed in the statement of theorem.
 \end{proof}

\subsection{Helices  in $\e_1^3$}

In Euclidean space, a helix is a curve whose tangent straight-lines make a constant angle with a fixed direction. This direction is called the axis of the helix. A result due to Lancret shows that a curve is a helix if and only if $\tau/\kappa$ is a constant function. For example, plane curves are helices. A helix with constant curvature and torsion is called a  cylindrical helix.

We try to extend this notion to the Lorentz-Minkowski space. The problem appears when we speak of the angle between two vectors because the angle is not defined for all couple of vectors of $\e_1^3$.  For example, this is the case if the curve  is lightlike. In other cases, it may also occur problems when the axis is timelike (or spacelike) and the curve is spacelike (or timelike). Even in the case  that the straight-line and the curve are timelike, the directions of the axis ant the tangent vector may not be in the same timelike cone.   For these reasons, we extend the notion of helix in $\e_1^3$ as follows:

\begin{definition}
A helix $\alpha:I\rightarrow\e_1^3$   is a regular curve  parametrized by arc length (or by the pseudo arc length if $\alpha$ is lightlike) such that  there exists a vector $v\in \e_1^3$ with the property that the function  $\langle \t(s),v\rangle$ is  constant. Any line parallel this direction $v$ is called the axis of the helix.
\end{definition}

In particular, a straight-line and a plane curve are helices. In what follows, we discard both situations.

For Frenet curves, there holds the equivalence of a helix in terms of the constancy of $\tau/\kappa$  and the proof follows the same steps as in Euclidean space.

\begin{theorem} Let $\alpha:I\rightarrow\e_1^3$ be a Frenet curve. Then $\alpha$ is a helix if and only if $\tau/\kappa$ is constant.
\end{theorem}

We analyse  what   happens in the other cases. Let $\alpha$ be spacelike curve with lightlike normal vector. Suppose $\langle\t(s),v\rangle=a\in\r$. By differentiating $\langle\t(s),v\rangle$ we have $\langle\n(s),v\rangle=0$. Since $\n(s)$ is lightlike, there exists  a function $b=b(s)$ such that  $v=a\t(s)+b(s)\n(s)$. Differentiating again, and using the Frenet equations, we have $(b'+b\tau +a)\n(s)=0$. Thus, $b'+b\tau +a=0$. This says that any spacelike curve with lightlike normal vector is a helix since the vector $v$ is {\it any} vector of type $v=a\t(s)+b(s)\n(s)$, where $a\in\r$ and $b$ satisfies the above ODE.  Remark that the Frenet equations implies that $v$ does not depend on $s$, that is, $v$ is a fix vector.

Assume now that $\alpha$ is a lightlike curve such that the function $\langle\t(s),v\rangle $ is a constant $a\in\r$, $a\not=0$, for some vector $v\in\e_1^3$. Then $\langle\n(s),v\rangle=0$. As $\t(s)$ is lightlike, there exists a function $b=b(s)$ such that $v=b(s)\t(s)-a\b(s)$. By differentiating and using the Frenet equations,
$$0=b'\t+(b-a\tau)\n.$$
Then $b$ is a constant function and $\tau$ is constant. The reverse always holds, that is, any lightlike curve with constant torsion is a helix. For this, we take any $a\not=0$ and consider $b=a\tau$. Then the vector $v=b\t-a\b$, which does not depend on $s$, is an axis of the helix.

\begin{theorem} A spacelike curve with lightlike normal vector is a helix. A lightlike curve is a helix if and only if its torsion is constant.
\end{theorem}

\subsection{Angle between two vectors}

We revisit the notion of angle in $\e_1^2$ given in section 1.  In Euclidean plane $\e^2$, let $u,v\in\s^1$ be two unit vectors. It is not difficult to prove that the angle  $\angle(u,v)$ between $u$ and $v$ is the arc length of the shortest curve in the circle $\s^1$ between the points $u$ and $v$.

We now consider two unit vectors in $\e_1^2$ and we assume that they are both timelike or both spacelike. Recall that in Lorentz-Minkowski plane $\e_1^2$ we have defined in \eqref{angleh} and \eqref{angles} the angle between two unit timelike   (resp. spacelike)   vectors that lie in the same component of $\u_1^2$. The cases are:
\begin{enumerate}
\item If $u,v$ are timelike, then $u,v\in \h_{+}^1$ or $u,v\in\h_{-}^1$.
\item If $u,v$ and spacelike, then $u,v\in\s_1^{1+}$ or $u,v\in\s_1^{1-}$.
\end{enumerate}
As in Euclidean space, we prove that $\angle(u,v)$ is the length of piece of circle joining both points.

\begin{theorem}\label{te-angle} Let $u,v$ be two unit vectors $\e_1^2$ in the same component of $\u_1^2$. The angle $\angle(u,v)$ is the length of the arc of $\u_1^2$ joining $u$ and $v$.
\end{theorem}
\begin{proof}
Take $u,v$ two unit timelike vectors  in the same timelike cone, and without loss of generality, $u,v\in\h^1_{+}$. Suppose $u=(\sinh(\varphi),\cosh(\varphi))$, $v=(\sinh(\psi),\cosh(\psi))$ with $\varphi\leq\psi$. The arc of $\h^1$ between $u$ and $v$ is $\alpha(s)=(\sinh(s),\cosh(s))$, $s\in[\varphi,\psi]$, and the length between $u$ and $v$ is
$$L_{\varphi}^{\psi}(\alpha)=\int_\varphi^{\psi}|\alpha'(s)|\ ds=\int_\varphi^{\psi} 1\ ds=\psi-\varphi.$$
If we compute the angle between both vectors,
$$-\langle u,v\rangle=-\sinh(\varphi)\sinh(\psi)+\cosh(\varphi)\cosh(\psi)=\cosh(\psi-\varphi)$$
and so $\angle(u,v)=\psi-\varphi$ because $\psi-\varphi\geq 0$.

Consider now two unit spacelike vectors $u,v$ that belong to the same connected component  of $\u_1^2$. Assume  $u,v\in\s^{1+}_1$ with $u=(\cosh(\varphi),\sinh(\varphi))$, $v=(\cosh(\psi),\sinh(\psi))$ and $\varphi\leq \psi$. Then the arc $\alpha$ of $\s_1^{1+}$ joining $u$ and $v$ is $\alpha(s)=(\cosh(s),\sinh(s))$, $s\in[\varphi,\psi]$. The length of $\alpha$ between $u$ and $v$ is
$$L_{\varphi}^{\psi}(\alpha)=\int_\varphi^{\psi}|\alpha'(s)|\ ds=\int_\varphi^{\psi} 1\ ds=\psi-\varphi.$$
 The angle $\angle(u,v)=\theta$ satisfies $\langle u,v\rangle=\cosh(\theta)$. Since
 $$\langle u,v\rangle=\cosh(\varphi)\cosh(\psi)-\sinh(\varphi)\sinh(\psi)=\cosh(\psi-\varphi),$$
and $\psi-\varphi\geq 0$,  the angle $\theta$ is $\psi-\varphi$.
 \end{proof}

We point out  that other authors (e.g. \cite{npv}) have defined the angle slightly different. Let $u,v\in\h_{+}^1$ (or $u,v\in\h^1_{-}$). Then there exists $\varphi\in [0,\infty)$ that that $R_\varphi(u)=v$, where
 $$R_\varphi=\left(\begin{array}{cc}\cosh(\varphi)&\sinh(\varphi)\\ \sinh(\varphi)&\cosh(\varphi)\end{array}\right)\in O_1^{++}(2).$$
 Then $\angle (u,v)=\varphi$.

 For unit spacelike vectors   $u,v\in\s_1^{1+}$ (or $u,v\in\s_1^{1-}$), a reflection $F$ with respect to the line  $y=x$ (included in the the lightcone of $\e_1^2$) gives $F(u), F(v)\in\h^1_{+}$ or $F(u),F(v)\in \h^1_{-}$. Then the angle $\angle(u,v)$ is defined as $\angle(F(u),F(v))$. In both cases, the definition agrees with the one given here. In \cite{npv} the authors {\it also} define the angle between unit timelike and  spacelike vectors in different components of $\u_1^2$ carrying both vectors to $\h_{+}^1$ by the use of  successive reflections across the lightlike cone of $\e_1^2$. However in the new cases, there is not an arc of $\u_1^2$ joining both vectors.


\section{Surfaces in Minkowski space}

We introduce the notion of spacelike and timelike surface and we will define the mean curvature and the Gaussian curvature for this kind of surfaces. Next, we will compute these curvatures by using parametrizations and, finally, we will characterize umbilical and isoparametric surfaces of $\e_1^3$. The development of this chapter is similar to the Euclidean space, even in the local formula of the curvature. However, we will see how the causal character  imposes restrictions, as for example, the surfaces can not be closed and that the Weingarten map for timelike surfaces may not be diagonalizable.

\subsection{Spacelike and timelike surfaces in  $\e_1^3$}

Let $M$ be a smooth, connected  surface possibly with non-empty   boundary $\partial M$.  Let  $x:M\rightarrow \e_1^3$  be an  immersion, that is, a differentiable map such that its differential map  $dx_p:T_p M\rightarrow \r^3$ is injective. We identify the tangent plane $T_p M$ with  $(dx)_p(T_p M)$. We consider the pullback metric $x^*(\langle,\rangle_p)$, that is,
$$x^*(\langle,\rangle_p)(u,v)=\langle dx_p(u),dx_p(v)\rangle,\ u,v\in T_pM,$$
so $x:(M,x^{*}\langle,\rangle)\rightarrow(\e_1^3,\langle,\rangle)$ is an isometric immersion. The metric $x^{*}\langle,\rangle$ can be positive definite, a metric with index $1$ or a degenerate metric.

\begin{definition} Let $M$ be a surface. An immersion $x:M\rightarrow\e^3$ is called spacelike (resp. timelike, lightlike) if all tangent planes $(T_pM,x^{*}(\langle,\rangle_p))$ are spacelike (resp. timelike, lightlike). A non-degenerate surface is a spacelike or timelike surface.
\end{definition}

As the curves of $\e_1^3$, given an immersed surface in $\e_1^3$, the causal character may change in different points of the same surface. This means that a surface is not necessarily classified in one of the above types. For example, in the sphere $\s^2=\{(x,y,z)\in\r^3: x^2+y^2+z^2=1\}$,  the region $\{(x,y,z)\in\s^2: |z|<1/\sqrt{2}\}$ is timelike,  $\{(x,y,z)\in\s^2: |z|>1/\sqrt{2}\}$ is spacelike  and $\{(x,y,z)\in\s^2: |z|=1/\sqrt{2}\}$ is lightlike. Similarly, we point out that the spacelike and timelike conditions are open properties.

For a spacelike (resp. timelike) surface $M$ and $p\in M$, we have the decomposition    $\e_1^3=T_p M\oplus (T_p M)^\bot$, where $ (T_p M)^\bot$ is a timelike (resp. spacelike) subspace of dimension $1$. A Gauss map is a differentiable map $N:M\rightarrow\e_1^3$ such that $|N(p)|=1$ and $N(p)\in  (T_p M)^\bot$ for all $p\in M$. Let us recall that a surface is orientable if there is a family of coordinate charts where the change of parameters has positive Jacobian. For a non-degenerate surface this is equivalent to the existence of a Gauss map, called also an orientation of $M$. Recall that locally a surface is a graph of a function and thus, it is locally orientable.

 The causal character of an immersion imposes conditions on the surface $M$. For example, we have:

\begin{proposition} Let $M$ be a compact surface and let $x:M\rightarrow\e_1^3$ be a spacelike, timelike or lightlike immersion.  Then  $\partial M\not=\emptyset$.
\end{proposition}

\begin{proof} Assume $\partial M=\emptyset$. Consider that the immersion is spacelike  (resp. timelike or lightlike). Let $a\in\e^3$ be a spacelike (resp. timelike) vector. Since $M$ is compact, let $p_0\in M$ be the  minimum of the function $f(p)=\langle x(p),a\rangle$. As $\partial M=\emptyset$, then $p_0$ is a critical point of the function $f$ and so, $\langle (dx)_{p_0}(v),a\rangle$, $\forall v\in T_{p_0}M$. Then $a\in (T_{p_0}M)^\bot$, a contradiction because $(T_{p_0}M)^\bot$ is timelike (resp. spacelike or lightlike).
 \end{proof}

This result discards the study in Minkowski space $\e_1^3$ of the class of closed surfaces (compact without boundary), which plays an important role in Euclidean space.

\begin{proposition}\label{co-closed} Let $x:M\rightarrow\e_1^3$ be a spacelike immersion of a surface $M$ and let the projection map  $\pi:M\rightarrow\r^2$, $\pi(x,y,z)=(x,y)$. \begin{enumerate}
\item The projection $\pi$ is a local diffeomorphism.
\item Assume that $M$ is compact and that
$x_{|\partial M}$ is a diffeomorphism between  $\partial M$ and a plane, closed, simple curve. Then  $x(M)$ is a graph on the planar domain determined by $x(\partial M)$.
\end{enumerate}
\end{proposition}

\begin{proof}
\begin{enumerate}
\item The map  $\pi:M\rightarrow\r^2$ satisfies     $$|(d\pi)_p(u)|^2=|(u_1,u_2,0)|^2=u_1+u_2^2>u_1^2+u_2^2-u_3^2=|u|^2.$$
This means that  $d\pi_p$ is an isomorphism and $\pi$ is a local diffeomorphism.
\item Let $\Omega$ be the planar domain that encloses $x(\partial M)$. We know that  $\pi:M\rightarrow\r^2$ is a local diffeomorphism. We claim that  $\pi(M)\subset \Omega$. On the contrary, let $q\in\partial\pi(M)\setminus \overline{\Omega}$  and let   $x(p)=q$. Then $p\not\in \partial M$ because
$x(\partial M)=\partial \Omega$, and so, $p$ is an interior point.  But $x(p)\in\partial\pi(M)$ and then, the tangent plane at  $p$ must be vertical: contradiction. This shows the claim.

Thus $\pi:M\rightarrow \Omega$ is a local diffeomorphism. In particular, it is a covering map. Since $\Omega$ is simply connected,
 $\pi\circ x$ is a diffeomorphism, which means that  $x(M)$ is a graph on  $\Omega$.
 \end{enumerate}
\end{proof}

We study the causal character of surfaces of $\e_1^3$.
\begin{example}
\begin{enumerate}
\item A  plane $P=p_0+\mbox{Span}\{v\}^\bot$. The causal character of $P$ coincides with the one of  $v$. If $v$ is a unit     timelike or spacelike vector, then a Gauss maps is $N(p)=v$.
\item  A \emph{hyperbolic plane} of center $p_0\in\e_1^3$ and radius $r>0$ is the surface
$$\h^{2}(r;p_0)=\{p\in\e_1^3:\langle p-p_0,p-p_0\rangle=-r^2, \langle p-p_0,E_3\rangle<0\}.$$
Here $E_3=(0,0,1)$. We observe that the set $\{p\in\e_1^3:\langle p-p_0,p-p_0\rangle=-r^2\}$ has exactly two connected components and that the condition $\langle p-p_0,E_3\rangle <0$ chooses one of them. This component, when $p_0$  is the origin of $\r^3$ and $r=1$, is denoted  by $\h^2$, that is,
$$\h^2=\{(x,y,z)\in\e_1^3: x^2+y^2-z^2=-1, z>0\}.$$
From the Euclidean viewpoint, this surface is one component of a hyperboloid of two sheets. We will justify in Ex. \ref{ex-complete} why $\h^2(r;p_0)$ is called a hyperbolic plane. A hyperbolic plane   is a spacelike surface. Indeed, if $v\in T_p\h^2(r;p_0)$ and $\alpha=\alpha(s)\subset \h^2(r;p_0)$ is a curve that represents $v$, then $\langle \alpha(s)-p_0,\alpha(s)-p_0\rangle=-r^2$. By differentiating with respect to $s$ and letting $s=0$, we obtain $\langle v,p-p_0\rangle=0$. This means that $T_p M=\mbox{Span}\{p-p_0\}^\bot$. As $p-p_0$ is a timelike vector, then $M$ is a spacelike surface. Moreover,  $N(p)=(p-p_0)/r$ is a Gauss map. Since   $\langle N,E_3\rangle<0$, $N$ is future directed. See Fig.  \ref{3fig1}.

\item The {\it pseudosphere} of  center $p_0$ and radius $r$  is the surface
$$\s^{2}_1(r;p_0)=\{p\in\e_1^3:\langle p-p_0,p-p_0\rangle=r^2\}.$$
The tangent plane at $p$ is  $T_pM=\mbox{Span}\{p-p_0\}^\bot$ and $N(p)=(p-p_0)/r$. This vector is spacelike vector and so, the surface is timelike. If $p_0$ is the origin and $r=1$, the surface is also called the {\it De Sitter space} and we denote by $\s^2_1$. Then
 $$\s_1^2=\{(x,y,z)\in\e_1^3: x^2+y^2-z^2=1\}.$$
See Fig.  \ref{3fig1}. From an Euclidean viewpoint, this surface is a ruled hyperboloid, also called, a hyperboloid of one sheet.
\item The \emph{lightlike cone} of center $p_0$ is
 $$\mathcal{C}(p_0)=\{p\in\e_1^3:\langle p-p_0,p-p_0\rangle=0\}-\{p_0\}.$$
Here $T_p\mathcal{C}(p_0)=\mbox{Span}\{p-p_0\}^\bot$. The surface is lightlike. If $p_0$ is the origin of $\r^3$, then $\mathcal{C}(p_0)$ is the lightlike cone $\mathcal{C}$ of $\e_1^3$.
\begin{figure}[hbtp]
\begin{center}\includegraphics[width=.3\textwidth]{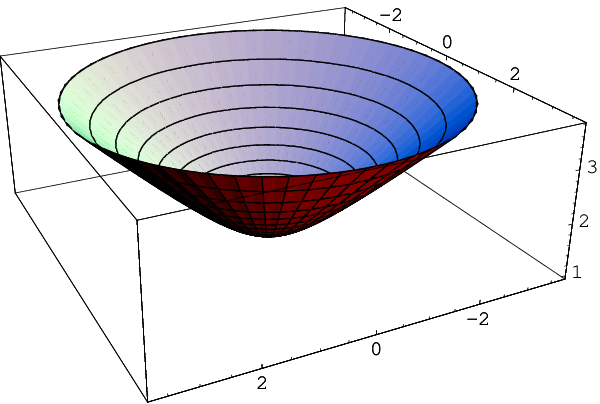}\ \ \includegraphics[width=.3\textwidth]{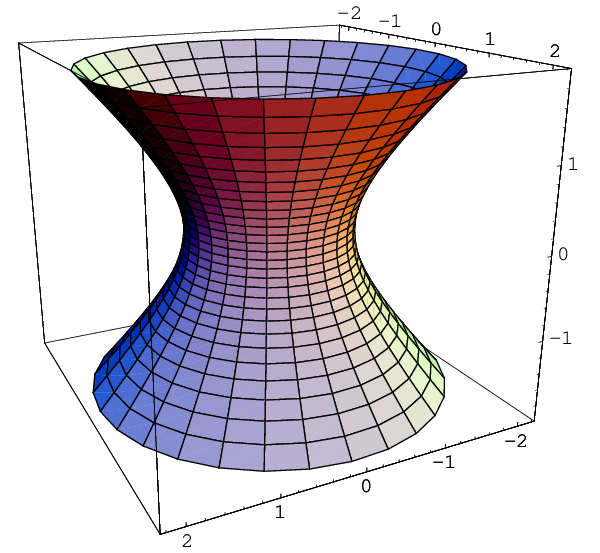}\ \ \includegraphics[width=.3\textwidth]{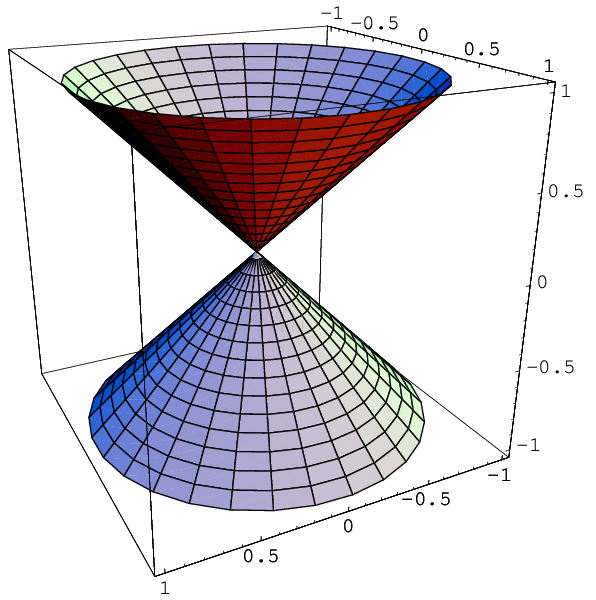}\end{center}
\caption{The hyperbolic plane, the pseudosphere and the lightlike cone of center the origin of $\r^3$.\label{3fig1}}
\end{figure}

\item  Let $f:\Omega\subset\r^2\rightarrow\r$ be a smooth function defined on a domain  $\Omega\subset\r^2$.  Let the graph of $f$ defined by
\begin{equation}\label{planegraph}
    S=\mbox{graph}(f)=\{(x,y,f(x,y)):(x,y)\in\Omega\}.
    \end{equation}
    Consider $S$ as the image of the immersion
 $$\Psi:\Omega\rightarrow\e_1^3,\ \ \Psi(x,y)=(x,y,f(x,y)).$$
 As $\Psi_x=(1,0,f_x)$ and $\Psi_y=(0,1,f_y)$, the matrix of the induced metric with respect to $\{\Psi_x,\Psi_y\}$  is
 $$\left(\begin{array}{cc}
 1-f_x^2 &-f_x f_y\\ -f_x f_y&1-f_y^2\end{array}\right)$$
and the determinant is
\begin{equation}\label{grad}
1-f_x^2-f_y^2=1-|\nabla f|^2.
\end{equation}
Thus the immersion is spacelike if $|\nabla f|^2<1$,  timelike if $|\nabla f|^2>1$ and lightlike if $|\nabla f|^2=1$.

Here we observe a difference with the Euclidean setting (compare with \cite[p. 58]{dc}). Given a function $f$, one can consider the graph of $f$ on the $xy$-plane, such as it is done in \eqref{planegraph} but also on the $yz$-plane or on the $xz$-plane. In each case, the causal character changes for the {\it same} function $f$. For example, if the surface is a graph on the timelike plane of equation $x=0$ given by  $Q=\{(f(y,z),y,z): (y,z)\in\Omega\}$, the matrix of the metric is now
 $$\left(\begin{array}{cc}
 1+f_y^2 &f_y f_z\\ f_y f_z&f_z^2-1\end{array}\right).$$
 The determinant is $-f_y^2+f_z^2-1$, which is different than $1-|\nabla f|^2$ as in \eqref{grad} and the sign determines the causal character of the surface. Thus the same function $f$ may give surfaces with different causal character. For example, if   $\Omega=\r^2$ and $f(x,y)=0$, then $S$ in \eqref{planegraph} is a spacelike (horizontal) plane but $Q$ is a timelike (vertical) plane.
\item Surface given by a regular value. Let $f:\Omega\subset\r^3\rightarrow\r$ be a differentiable function and $a\in\r$ a regular value of $f$, that is, $(df)_p$ is regular for all $p\in f^{-1}(\{a\})$. It is well known that $S=f^{-1}(\{a\})$ is a surface. In particular, the gradient $\nabla f$ computed with the metric $\langle,\rangle$ is an orthogonal vector to $S$ and this gives the causal character of $S$. Indeed, for any $v\in T_pS$, if $\alpha=\alpha(s)$ is a curve representing $v$, $f(\alpha(s))=0$ and so, $(df)_p(v)=0$. This writes as
    $\langle(\nabla f)_p,v\rangle=0$. It is not difficult to prove that   $\nabla f$ is the reflection of  the Euclidean gradient $\nabla_e f$   with respect to the plane of equation $z=0$.

    Consider the next example. Let $f(x,y,z)=x^2+y^2-z^2$. Then $(df)_p(v)=2(xv_1+yv_2-zv_3)$, with $p=(x,y,z)$ and $v=(v_1,v_2,v_3)$. Then $p$ is a critical point only if $p=(0,0,0)$. In such a case, $f(0,0,0)=0$. This means that for all $a\not=0$, $S_a=f^{-1}(\{a\})$ is a surface. Here
    $$\nabla f=2(x,y,z),\ \ \nabla_e f=2(x,y,-z).$$
    Consider $a\in\{-1,1\}$. Then
    $$\langle \nabla f,\nabla f\rangle=4(x^2+y^2-z^2)=4f(x,y,z)=4a.$$
    Therefore, if $a<0$ (resp. $a>0$)  the surface is spacelike (resp. timelike). In fact, $S_{-1}$ has two (spacelike) components, one of them is $\h^2$ and $S_1$ is the pseudosphere $\s_1^2$.
\end{enumerate}
\end{example}

Given a surface $M$, we know that an immersion $x:M\rightarrow\r^3$ is locally a graph on one of the three coordinate planes. Remark that this result is not metric. If we endow on $\r^3$ the Lorentzian metric and  for a spacelike or timelike surface, we can precise with respect to what coordinate plane is a graph.

\begin{proposition} A spacelike (resp. timelike) surface is locally the graph of a function defined in the plane of equation  $z=0$ (resp. $x=0$ or $y=0$).
\end{proposition}
\begin{proof}
We consider a local parametrization of the surface
$$\Psi=\Psi(u,v)=(x(u,v),y(u,v),z(u,v)).$$
Since the vector  $\Psi_u\times \Psi_v$ is orthogonal to
 $\Psi_u$ and  $\Psi_v$, then $\Psi_u\times \Psi_v$  is a timelike (resp. spacelike) vector. Thus, its third coordinate does not vanish (resp. the first or the second coordinate). This coordinate is
 $$-\left|\begin{array}{cc}x_u&y_v\\ x_v&y_v\end{array}\right|\ \ (\mbox{resp.}
 \left|\begin{array}{cc}y_u&z_u\\ y_v&z_v\end{array}\right|, \mbox{ or } \left|\begin{array}{cc}x_u&z_u\\ x_v&z_v\end{array}\right|)
 .$$
The implicit function theorem  asserts that around a point of the surface, the map
$$\phi:(u,v)\longmapsto (x(u,v),y(u,v))$$
(resp.
$$(u,v)\longmapsto (y(u,v),z(u,v)),\ \ \ \  (u,v)\longmapsto (x(u,v),z(u,v)))$$
is a diffeomorphism. We reparametrize the immersion  by $\Psi\circ\phi^{-1}$. Then locally the surface is  the graph of the function
$z\circ\phi^{-1}$ (resp. $x\circ\phi^{-1}$ or $y\circ\phi^{-1}$).
\end{proof}

We end this subsection showing that a spacelike surface is orientable.

\begin{theorem}\label{t-or} Let $M$ be a surface and  let  $x:M\rightarrow\e_1^3$ be a spacelike immersion. Then  $M$ is  orientable.
\end{theorem}

\begin{proof} We know that for each connected coordinate open $\Omega\subset M$ there exists two unit orthogonal vectors $N$ and $-N$.  Since $N(p)$ and $E_3$ are both timelike, we have   $\langle N(p),E_3\rangle\geq 1$ or $\langle N(p),E_3\rangle \leq -1$. On $\Omega$ we choose $N$ such that $\langle N,E_3\rangle\leq -1$. This allows us to define a globally normal vector field on $M$ and $M$ is an orientable surface.
\end{proof}
This result imposes an a priori condition to an abstract surface $M$ (a $2$-dimensional manifold)  to be isometrically immersed in $\e_1^3$ as a spacelike surface: $M$ must be orientable.

Depending if the surface is spacelike or timelike, the codomain of the Gauss map is a hyperbolic plane or a pseudosphere.
Exactly,
\begin{enumerate}
\item If the immersion is spacelike, the Gauss map pointing to the future is a map $N:M\rightarrow\h^2$.
\item If the immersion is timelike, the (local) Gauss map writes as $N:M\rightarrow\s_1^2$.
\end{enumerate}
This means that in Lorentz-Minkowski space $\e_1^3$ the hyperbolic plane $\h^2$ and the pseudosphere $\s_1^2$ play the  role of a sphere $\s^2$ of $\e^3$.

\subsection{The mean curvature of a spacelike surface}\label{sec-122}
Let $x:M\rightarrow \e_1^3$ be a spacelike or timelike immersion of a surface $M$. Denote $\mathfrak{X}(M)$ the space of tangent vector fields to $M$ and denote by $\nabla^0$ the Levi-Civita connection of $\e_1^3$. If $X,Y\in  \mathfrak{X}(M)$, we have the decomposition
$$\nabla^0_XY=(\nabla^0_XY)^\top+(\nabla^0_XY)^\bot,$$
where $\top$ and $\bot$ denote the tangent part and the normal part with respect to $M$ of $\nabla^0_XY$, respectively. Denote $\nabla$ the induced connection on $M$ by the immersion $x$, that is,
$$\nabla_XY=(\nabla^0_XY)^\top$$
and we define  the second fundamental form of $x$ as
 the tensorial, symmetric map
 $$\sigma:\mathfrak{X}(M)\times\mathfrak{X}(M)\rightarrow {\mathfrak X}(M)^\bot,\ \ \sigma(X,Y)=(\nabla_X^0 Y)^\bot.$$
The expression of the Gauss formula is
\begin{equation}\label{3-gf}
\nabla_X^0 Y=\nabla_X Y+\sigma(X,Y),\ \ X,Y\in \mathfrak{X}(M).
\end{equation}
 Consider now $\xi$  a normal vector field to  $x$ and let $A_\xi (X)$ be the tangent component of
 $-\nabla_X^0 \xi$,
$$A_\xi(X)=-(\nabla_X^0 \xi)^\top.$$
We have from \eqref{3-gf}
\begin{equation}\label{3-gf2}
\langle A_\xi(X),Y\rangle=\langle\sigma(X,Y),\xi\rangle.
\end{equation}
Because $\sigma$ is symmetric,  \eqref{3-gf2} implies
\begin{equation}\label{self}
\langle A_\xi(X),Y\rangle=\langle X,A_\xi(Y)\rangle.
\end{equation}
This means that $A_\xi$  is  self-adjoint with respect to the metric of $M$.

 Let $N$ be a (local) unit normal vector field on $M$. We know that if the immersion is spacelike, the surface is always orientable by Th. \ref{t-or}. Denote
$$\langle N,N\rangle=\epsilon=\left\{\begin{array}{cl}-1 &\mbox{if $M$ is spacelike}\\
 1&\mbox{if $M$ is timelike.}
 \end{array}
 \right.$$
Take in the above formula $\xi=N$.  Because $\langle N,N\rangle$ is constant, we have $\langle\nabla_X^0N,N\rangle=0$. Then $\nabla^0_XN$ is tangent to $M$. Denote
\begin{equation}\label{3-wf}
-\nabla^0_XN=A_N(X)\ \mbox{(Weingarten formula)}
\end{equation}

\begin{definition} The Weingarten endomorphism at $p\in M$ is defined by
$$A_p:T_pM\rightarrow T_pM,\ A_p=(A_N(X))_p.$$
Moreover   \eqref{3-wf} gives
$$A_p(v)=-\nabla^0_v N=-(dN)_p(v),\ \ v\in T_pM.$$
 We will write $AX$ instead of $A_N(X)$.
\end{definition}
Since $\sigma(X,Y)$ is proportional to $N$,   from  \eqref{3-gf} and \eqref{3-gf2} we deduce
\begin{equation}\label{s-1}
\sigma(X,Y)=\epsilon\langle\sigma(X,Y),N\rangle N=\epsilon\langle AX,Y\rangle N.
\end{equation}
Now \eqref{3-gf} writes as
$$\nabla^0_XY=\nabla_XY+\epsilon\langle AX,Y\rangle N.$$

We define the mean curvature and the Gauss curvature. Firstly, recall how is defined in Euclidean space. In $\e^3$,  the Weingarten map is diagonalizable because it is a self-adjoint endomorphism with respect to a Riemannian metric. The principal curvatures are the eigenvalues of the Weingarten map, and hence, the Gauss curvature and the mean curvature is the product and the arithmetic average of the principal curvatures, respectively.

The identity \eqref{self} says that  the Weingarten map $A$ is self-adjoint  with respect to the induced metric $\langle,\rangle$. If the metric is Riemannian, then $A$ is diagonalizable but if the metric is Lorentzian, the map $A$ could be not diagonalizable.  In other words, for a spacelike surface the principal curvatures are well defined but in a timelike surface could not be defined. Thus we have to choose other approach to define the mean curvature in both types of surfaces and for this purpose, we will consider the trace of the second fundamental form.

\begin{definition}
Let $M$ be a surface and let $x:M\rightarrow\e_1^3$ be a non-degenerate immersion. The mean curvature vector field $\vec{H}$  is the vector
$$\vec{H}=\frac12\mbox{ trace}(\sigma),$$
where the trace is computed with  respect to the metric of the surface.
The mean curvature function $H$ is defined  by the relation $\vec{H}=HN$. Therefore
$$H=\epsilon\langle\vec{H},N\rangle.$$
\end{definition}
Observe that $\vec{H}$ is a vector field orthogonal to $M$, that is, $\vec{H}\in \mathfrak{X}(M)^\bot$.
We write $\vec{H}$ and $H$ in terms of a local tangent basis. Let  $\{e_1,e_2\}$ be  an orthonormal local tangent vector fields on $M$ where $e_1$ is spacelike and $\langle e_2,e_2\rangle=-\epsilon$. Then \eqref{s-1} gives
\begin{eqnarray*}\vec{H}&=&\frac12\mbox{trace}(\sigma)=\frac12 (\sigma(e_1,e_1)-\epsilon \sigma(e_2,e_2))\\
&=&\frac{1}{2}\left(\epsilon\langle Ae_1,e_1\rangle-\langle Ae_2,e_2\rangle\right)N\\
&=&\frac{\epsilon}{2}\left(\langle Ae_1,e_1\rangle-\epsilon\langle Ae_2,e_2\rangle\right)N=\left(\frac{\epsilon}{2}\ \mbox{trace} A\right)N.
\end{eqnarray*}
On the other hand,
$$H=\epsilon\langle\vec{H},N\rangle=\frac{\epsilon}{2}\left(\langle Ae_1,e_1\rangle-\epsilon\langle Ae_2,e_2\rangle\right)=\frac{\epsilon}{2}\mbox{ trace}(A).$$
\begin{corollary} The mean curvature of a non-degenerate surface is
\begin{equation}\label{dh}
H=\frac{\epsilon}{2}\mbox{ trace}(A).
\end{equation}
\end{corollary}

We define the (intrinsic) Gauss curvature $K$ of the surface.  For a surface, $\rho=2K$, where $\rho$ is the scalar curvature. For this, we compute the curvature tensor of the surface (here we follow \cite{on}).

Denote by $R^0$ and $R$  the curvature tensors of $\e_1^3$ and $M$, respectively. Since  $R^0=0$, we can compute $R$.   Let $X,Y,Z\in\mathfrak{X}(M)$. We know that
$$R^0(X,Y)Z=\nabla_X^0\nabla^0_YZ-\nabla^0_Y\nabla^0_XZ-\nabla^0_{[X,Y]}Z.$$
 Also,  $\nabla_Y^0Z=\nabla_YZ+\sigma(Y,Z)$. Because $\sigma(Y,Z)=\epsilon\langle AY,Z\rangle N$, and using \eqref{3-gf}, we have
 \begin{eqnarray*}
 \nabla^0_X\nabla^0_YZ&=&\nabla^0_X(\nabla_YZ)+\nabla^0_X\sigma(Y,Z)\\
 &=&\nabla_X\nabla_YZ+\sigma(X,\nabla_YZ)-\epsilon\langle AY,Z\rangle AX+\epsilon\langle AY,Z\rangle N.
 \end{eqnarray*}
The tangent part on $M$ is $\nabla_X\nabla_YZ-\epsilon\langle AY,Z\rangle AX$. Similarly, we calculate $ \nabla^0_Y\nabla^0_XZ$ and
$\nabla^0_{[X,Y]}Z$ and considering the tangent parts. Using that $R^0=0$ and that $R(X,Y)Z=\nabla_X\nabla_YZ-\nabla_Y\nabla_XZ-\nabla_{[X,Y]}Z$, we conclude
\begin{eqnarray}\label{tensor}
R(X,Y)Z&=&-\epsilon\langle AY,Z\rangle AX+\epsilon\langle AX,Z\rangle AY\nonumber\\
&=&\epsilon (\langle AX,Z\rangle AY-\langle AY,Z\rangle AX).
\end{eqnarray}
Hence we compute the Ricci tensor  and the scalar curvature $\rho$. For the Ricci tensor, we obtain
\begin{eqnarray*}
\mbox{Ric}(X,Y)&=&\mbox{trace}\left(v\longmapsto R(X,v)Y\right)=
\langle R(X,e_1)Y,e_1)-\epsilon\langle R(X,e_2)Y,e_2\rangle\\
&=&\epsilon\left(\langle AX,Y\rangle(\langle Ae_1,e_1\rangle-\epsilon\langle Ae_2,e_2\rangle)\right)-\epsilon\langle AX,AY\rangle\\
&=&\epsilon\left(\mbox{trace}(A)\langle AX,Y\rangle-\langle AX,AY\rangle\right)=2H\langle AX,Y\rangle-\epsilon\langle AX,AY\rangle.
\end{eqnarray*}
Thus
\begin{eqnarray}\label{rho}
\rho&=&\mbox{trace (Ric)}=R(e_1,e_1)-\epsilon R(e_2,e_2)\nonumber\\
&=& 2H\left(\langle Ae_1,e_1\rangle-\epsilon\langle Ae_2,e_2\rangle\right)-\epsilon(\langle Ae_1,Ae_1\rangle-\epsilon\langle Ae_2,Ae_2\rangle)\nonumber\\
&=&\epsilon \left((\mbox{trace}(A)^2- \mbox{trace}(A^2)\right)=4\epsilon H^2-\epsilon\ \mbox{trace}(A^2)\nonumber\\
&=&2\epsilon\mbox{ det}(A).
\end{eqnarray}
The expression of this matrix $A$ in the basis $\{e_1,e_2\}$ is
$$A=\left(\begin{array}{cc}\langle A e_1,e_1\rangle&\langle Ae_2,e_1\rangle\\ -\epsilon\langle Ae_1,e_2\rangle&-\epsilon\langle Ae_2,e_2\rangle\end{array}\right).$$
As $\rho=2K$,  the Gauss curvature $K$ is
\begin{equation}\label{kk}
K=\epsilon\mbox{ det}(A)=\frac{\epsilon}{2}\left(4H^2-\mbox{trace}(A^2)\right).
\end{equation}

\begin{corollary} Consider the Weingarten map $A$ of a non-degenerate surface of $\e_1^3$. Then
\begin{equation}\label{dk}
K=\epsilon\ \mbox{det}(A).
\end{equation}
\end{corollary}

One can also compute $K$ by observing that in a $2$-dimensional manifold, the Gauss curvature coincides with the sectional curvature of the $2$-plane generated by $\{e_1,e_2\}$, that, is,  of the tangent plane. As a consequence of \eqref{tensor}, we obtain
\begin{eqnarray*}
K&=&\frac{\langle R(e_1,e_2)e_2,e_1\rangle}{\langle e_1,e_1\rangle\langle e_2,e_2\rangle-\langle e_1,e_2\rangle^2}\\
&=& \frac{\epsilon(\langle Ae_1,e_1\rangle\langle Ae_2,e_2\rangle-\langle Ae_1,e_2\rangle\langle Ae_2,e_1\rangle}{-\epsilon}\\
&=&-\left(\langle Ae_1,e_1\rangle\langle Ae_2,e_2\rangle-\langle Ae_1,e_2\rangle^2\right).
\end{eqnarray*}
It is immediate that this expression coincides with \eqref{kk}.

 Returning to the case that $A$ is diagonalizable, we have:

\begin{definition} Consider $x:M\rightarrow\e_1^3$ a non-degenerate immersion and $p\in M$. If the Weingarten map $A_p$ is diagonalizable, the eigenvalues of $A_p$ are called the principal curvature at $p$, and we denote by $\lambda_1(p)$ and $\lambda_2(p)$.
\end{definition}

From \eqref{dh} and \eqref{dk}, we have:

\begin{corollary}\label{co-hk} Assume that $A_p$ is diagonalizable in a non-degenerate surface of $\e_1^3$. Then
$$H(p)=\epsilon\ \frac{\lambda_1(p)+\lambda_2(p)}{2},\ \ K(p)=\epsilon\ \lambda_1(p)\lambda_2(p).$$
\end{corollary}

In  Euclidean space, once defined the principal curvatures, it is given the notion of an umbilic as a point where the two principal curvatures coincide (\cite[p. 147]{dc}). Thus in Lorentz-Minkowski space we can not adopt the definition of umbilic point in terms of the principal curvatures. See \cite[p. 105]{on}.

\begin{definition} Let  $x:M\rightarrow\e_1^3$ be a spacelike or timelike immersion. A point  $p\in M$ is called umbilic if there exists $\lambda(p)\in\r$ such that $$\langle\sigma(u,v),N(p)\rangle=\lambda(p)\langle u,v\rangle,\ \ u,v\in T_pM.$$
 A surface is called totally umbilical if all points are umbilic.
\end{definition}
Thus, an umbilic is a point where the second and the first fundamental forms are proportional. Also, it is equivalent to say that $\langle A_pu,v\rangle=\lambda(p)\langle u,v\rangle$. In particular,   and from \eqref{3-gf2}, $A_p$ must be diagonalizable because  $\langle Ae_1,e_2\rangle=0$. Thus we can say that $p$ is umbilical if and only if $\lambda_1(p)=\lambda_2(p)$.
In Euclidean space, it is well know the inequality $H^2-K\geq 0$ and the equality hold only in a umbilic. Now in $\e_1^3$ we have to assume that the Weingarten map is diagonalizable.

\begin{proposition}\label{pr-hk}Assume that $M$ is a non-degenerate surface of $\e_1^3$, $p\in M$ and $A_p$ is diagonalizable. Then
$$H(p)^2-\epsilon K(p)\geq 0$$
and the equality holds if and only if $p$ is umbilic.   In particular, in a timelike surface, if $H(p)^2-K(p)<0$, then   $p$ is not umbilic.
\end{proposition}
\begin{proof}
From the definition of $H$ and $K$, we have
$$0\leq\Big(\frac{\lambda_1-\lambda_2}{2}\Big)^2=\Big(\frac{\lambda_1+\lambda_2}{2}\Big)^2-
\lambda_1\lambda_2=H^2-\epsilon K.$$
Moreover the equality holds at a point $p$ if and only if $\lambda_1(p)=\lambda_2(p)$, that is, $p$ is an umbilic.
\end{proof}
The diagonalization of the Weingarten map depends on the existence of real roots of its characteristic polynomial $P(\lambda)$. A simple computation leads to $P(\lambda)=\lambda^2-2H\epsilon\lambda+\epsilon K$ and its discriminant is $D=4(H^2-\epsilon K)$. Thus:
\begin{enumerate}
\item If $H^2-\epsilon K>0$, there are two different real roots of $P(\lambda)$ and the  Weingarten map is diagonalizable.
 \item If $H^2-\epsilon K<0$, $A$ is not diagonalizable.
\item If $H^2-\epsilon K=0$, there is a double root of $P(\lambda)$. Then: a) if $\epsilon=-1$, the root $\lambda=-H$ is the eigenvalue of $A$ and the point is umbilic; b) if $\epsilon=1$, the matrix could be or not be diagonalizable.
\end{enumerate}
Finally,
$$|\sigma|^2=\sum_{i,j=1}^2\langle Ae_i,e_j\rangle^2 =4H^2-2\epsilon K,$$
and if $A_p$ is diagonalizable, $|\sigma|^2=\lambda_1^2+\lambda_2^2$.

We point out that there exist non-umbilical timelike surfaces such that $H^2-K=0$ on the surface. See examples \ref{scroll} and \ref{ex-ru} below.

After the definition, we compute $H$ and $K$ in some surfaces that  have previously appeared.

\begin{example}\label{ex-complete}The next surfaces are umbilical.
\begin{enumerate}
 \item (Plane) Consider a non-degenerate plane $P=p_0+\mbox{Span}\{v\}^\bot$, with $|v|=1$. Then $N=v$  and $dN=0$. Here $\lambda_1=\lambda_2=H=K=0$.
 \item (Hyperbolic plane) The unit normal vector pointing to the future of  $\h^{2}(r;p_0)$    is $N(p)=(p-p_0)/r$. Then $A=-I/r$ and
     $$\lambda_1=\lambda_2=-\frac{1}{r},\ H=\frac{1}{r},\ K=-\frac{1}{r^2}.$$
 Thus a hyperbolic plane has constant negative curvature. Here we collect the properties of the surface $\h^2(r;p_0)$. This surface is a simply-connected Riemannian $2$-manifold with constant negative curvature. Moreover, it is geodesically complete. We do a proof of this fact for   $\h^2$. Given $p\in\h^2$ and $v\in T_p\h^2$ with $|v|=1$, the curve  $\alpha(s)=\cosh(s)p+\sinh(s) v$ is a geodesic starting from $p$ with velocity $v$: $\alpha''(s)=\alpha(s)$ and its the tangent part vanishes  $\alpha''(s)^T=0$. Finally, the geodesic $\alpha$ is defined for all $s$, which means that the surface is (geodesically) complete. With this in mind, $\h^2(r;p_0)$ is a $2$-space form of negative curvature, called usually the hyperbolic plane.

\item (Pseudosphere) For  $\s^{2}_1(r;p_0)$, the Gauss map is  $N(p)=(p-p_0)/r$. Again $A=-I/r$. Thus
     $$\lambda_1=\lambda_2=-\frac{1}{r},\ H=-\frac{1}{r},\ K=\frac{1}{r^2}.$$
     Thus a pseudosphere has constant positive curvature. Let us observe that the fundamental group of the pseudosphere is ${\mathbb Z}$.
\end{enumerate}
\end{example}

\begin{example}\label{ex-cy} Let us define right circular cylinders in $\e_1^3$.  Let $L$ be a non-degenerate straight-line. The right circular cylinder $C(r;L)$ of axis $L$ and radius $r>0$ is the set of points   equidistant from $L$ a distance $r$.

It is usual in the literature to call \emph{Lorentzian  cylinder} if the axis is timelike and   \emph{hyperbolic cylinder} when the axis is spacelike \cite{kl,mp}.
A right circular cylinder can be also viewed as a ruled surface whose basis is a Lorentzian circle centered at $L$ and orthogonal to $L$ and whose rulings are parallel to $L$.

 We  write explicit equations of a right circular cylinder in $\e_1^3$. Consider that $L$ passes  through $p_0\in\e_1^3$ and with direction $\vec{a}$. Assume $\vec{a}$ is a unit vector, with $\langle\vec{a},\vec{a}\rangle=\delta$, $\delta\in\{1,-1\}$.
 If $p\in\e_1^3$, the straight-line orthogonal to $L$ passing through $p$ intersects $L$ at the point
  $\bar{p}=p_0+\lambda \vec{a}$, with $\lambda=\delta\langle p-p_0,a\rangle$. If $\vec{a}$ is timelike, then
the vector  $p-\bar{p}$ is spacelike, but if $\vec{a}$ is spacelike, then  $p-\bar{p}$ can be spacelike, timelike and lightlike. Since we consider positive distances,  we discard the case that $p-\bar{p}$ is lightlike.  As a conclusion
$$C(r;L)=p\in\e_1^3: |p-p_0|^2-\delta\langle p-p_0,\vec{a}\rangle^2=\pm r^2\}.$$
In order to know the causal character of a right circular cylinder, a unit vector orthogonal  is
$$N(p)=\frac{1}{r}\Big(p-p_0-\delta\langle p-p_0,\vec{a}\rangle\vec{a}\Big).$$
Hence a basis of the tangent plane at $p$ is $\{(p-p_0)\times\vec{a},\vec{a}\}$ and the Weingarten map is $Av=-(v-\delta\langle v,\vec{a}\rangle \vec{a})/r$. Thus
$$A((p-p_0)\times\vec{a})=-\frac{1}{r}(p-p_0)\times\vec{a},\ \ A(\vec{a})=0.$$
As a consequence, the Weingarten map is diagonalizable. Furthermore:
\begin{enumerate}
\item (Lorentzian cylinder). If the axis is timelike, then  $N$ is spacelike and the surface is timelike. The surface writes as
$$C(r;L)=\{p\in\e_1^3:  |p-p_0|^2+\langle p-p_0,\vec{a}\rangle^2=  r^2\}.$$
Here
$$\lambda_1=-\frac{1}{r},\ \lambda_2=0,\ H=-\frac{1}{2r},\ K=0.$$
\item (Hyperbolic cylinder). If the axis   is spacelike, $N$ is timelike or spacelike.  There are two circular cylinders:
\begin{enumerate}
\item If $N$ is timelike,
$$C_s(r;L)=\{p\in\e_1^3:  |p-p_0|^2-\langle p-p_0,\vec{a}\rangle^2=  -r^2\},$$
$$\lambda_1=-\frac{1}{r},\ \lambda_2=0,\ H=\frac{1}{2r},\ K=0.$$
\item If $N$ is spacelike,
$$C_t(r;L)=\{p\in\e_1^3:  |p-p_0|^2-\langle p-p_0,\vec{a}\rangle^2=  r^2\},$$
$$\lambda_1=-\frac{1}{r},\ \lambda_2=0,\ H=-\frac{1}{2r},\ K=0.$$
\end{enumerate}
\end{enumerate}
Let us take explicit vectors $\vec{a}$ (see Fig. \ref{fcylinder}).
\begin{enumerate}
\item If $\vec{a}=(0,0,1)$,   $C(r;L)=\{(x,y,z)\in\r^3: x^2+y^2=r^2\}$. This surface is a vertical Euclidean circular cylinder.
\item If $\vec{a}=(1,0,0)$, then
\begin{enumerate}
\item the spacelike hyperbolic cylinder is $C_s(r;L)=\{(x,y,z)\in\r^3:y^2-z^2=-r^2, z>0\}$;
 \item the timelike hyperbolic cylinder is $C_t(r;L)=\{(x,y,z)\in\r^3: y^2-z^2= r^2, y>0\}$.
 \end{enumerate}
 \end{enumerate}
\begin{figure}[hbtp]
\begin{center}\includegraphics[width=.3\textwidth]{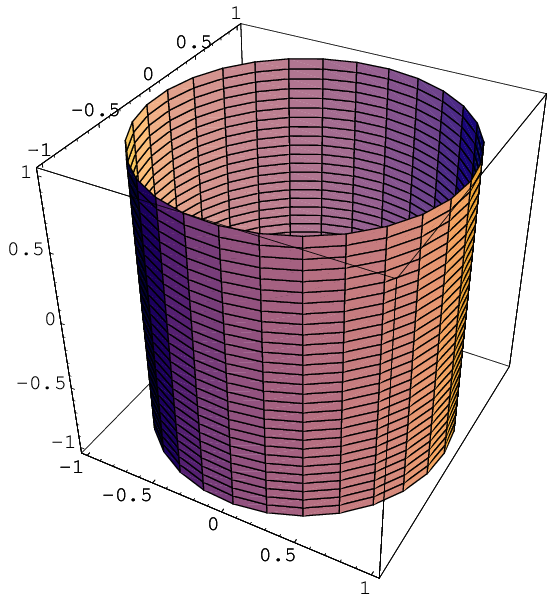}\ \ \includegraphics[width=.3\textwidth]{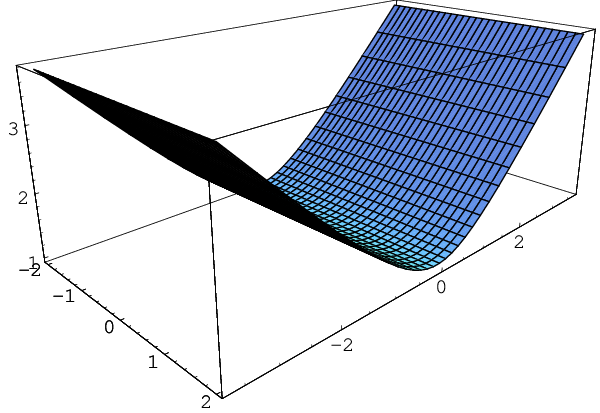}\ \ \includegraphics[width=.3\textwidth]{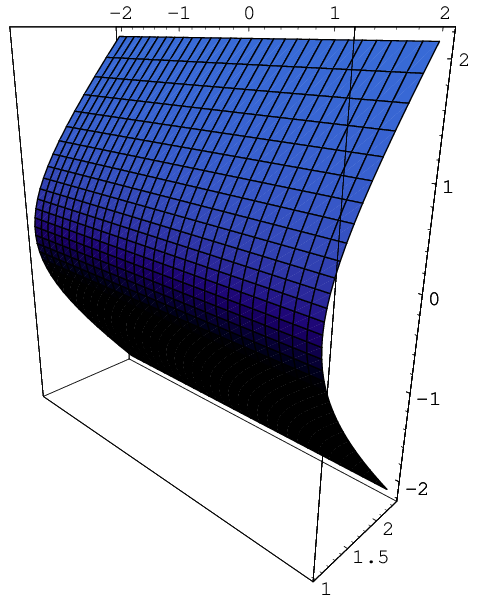}\end{center}
\caption{A Lorentzian cylinder, a spacelike hyperbolic cylinder and a timelike hyperbolic cylinder\label{fcylinder}}
\end{figure}
\end{example}

\subsection{Local computations of the curvature and examples}
We compute the curvatures of a non-degenerate surface by using a local parametrization. Here we follow the same ideas as in \cite{dc}. Consider a local parametrization
$$\Psi:U\subset\r^2\rightarrow\e_1^3,\ \ \Psi=\Psi(u,v),$$
of a (spacelike or timelike) immersion $x$. Let $B=\{\Psi_u,\Psi_v\}$ be a local basis of the tangent plane at each point of $\Psi(U)$.
With respect to  $B$, let  $\left(\begin{array}{cc}E&F\\F&G\end{array}\right)$ be the matricial expression of the first fundamental form, where
$$E=\langle \Psi_u,\Psi_u\rangle,\  F=\langle \Psi_u,\Psi_v\rangle, \ G=\langle \Psi_v,\Psi_v\rangle.$$
Denote $W=EG-F^2$. The surface is spacelike if $W>0$ and it is timelike if  $W<0$.  Take  the unit normal vector field
\begin{equation}\label{ene}
N=\frac{\Psi_u\times \Psi_v}{|\Psi_u\times \Psi_v|}.
\end{equation}
Again, we use the notation $\langle N,N\rangle=\epsilon$. Here
$$|\Psi_u\times\Psi_v|=\sqrt{-\epsilon(EG-F^2)}=\sqrt{-\epsilon W}.$$
Let  $\left(\begin{array}{cc}e&f\\f&f\end{array}\right)$ be the matricial expression of  $\sigma$ with respect to  $B$, that is,
\begin{eqnarray*}
e&=&-\langle A\Psi_u,\Psi_u\rangle=-\langle N_u,\Psi_u\rangle=\langle N,\Psi_{uu}\rangle\\
f&=&-\langle A\Psi_u,\Psi_v\rangle=-\langle N_v,\Psi_u \rangle=-\langle \Psi_v,N_u\rangle=\langle N,\Psi_{uv}\rangle\\
g&=&-\langle A\Psi_v,\Psi_v\rangle=-\langle N_v,\Psi_v \rangle=\langle N,\Psi_{vv}\rangle,
\end{eqnarray*}
where $A$ is the Weingarten map. Then
$$A=\left(\begin{array}{cc}E&F\\F&G\end{array}\right)^{-1} \left(\begin{array}{cc}e&f\\f&f\end{array}\right).$$
Therefore the mean curvature and the Gauss curvature are
\begin{equation}\label{hk}H=\epsilon\ \frac12\frac{eG-2fF+gE}{EG-F^2},\ \ K= \epsilon\ \frac{eg-f^2}{EG-F^2}.
\end{equation}
According to \eqref{ene}, we have
$$e=\langle N,\Psi_{uu}\rangle=\frac{\mbox{det}(\Psi_u,\Psi_v,\Psi_{uu})}{\sqrt{-\epsilon W}}.$$
$$f=\langle N,\Psi_{uv}\rangle=\frac{\mbox{det}(\Psi_u,\Psi_v,\Psi_{uv})}{\sqrt{-\epsilon W}}.$$
$$g=\langle N,\Psi_{vv}\rangle=\frac{\mbox{det}(\Psi_u,\Psi_v,\Psi_{vv})}{\sqrt{-\epsilon W}}.$$
Thus \eqref{hk} writes now as:
\begin{equation}\label{media2}
H=-\frac{\mbox{det}(\Psi_u,\Psi_v,\Psi_{uu})G-2\mbox{det}(\Psi_u,\Psi_v,\Psi_{uv})F+
 \mbox{det}(\Psi_u,\Psi_v,\Psi_{vv})E}{2(-\epsilon W)^{3/2}}.
\end{equation}
$$K=-\frac{\mbox{det}(\Psi_u,\Psi_v,\Psi_{uu})\mbox{det}(\Psi_u,\Psi_v,\Psi_{vv})-
\mbox{det}(\Psi_u,\Psi_v,\Psi_{uv})^2}{W^2}.$$

\begin{example} (Graph) Let $f\in C^2(\Omega)$ be a smooth function and consider the surface $M$ given by $z=f(x,y)$. Let $\Psi:\Omega\rightarrow\e_1^3$ denote the usual parametrization $\Psi(x,y)=(x,y,f(x,y))$. The coefficients of the first fundamental form are
$$E=1-f_x^2,\ F=-f_xf_y,\ G=1-f_y^2.$$
Thus $EG-F^2=1-f_x^2-f_y^2=1-|\nabla f|^2$. If the immersion is spacelike (resp. timelike) we have
 $|\nabla f|^2<1$ (resp. $>1$) on $\Omega$. The mean curvature $H$ satisfies
$$(1-f_y^2)f_{xx}+2 f_xf_y+(1-f_x^2)f_{yy}=-2H(-\epsilon(1-|\nabla f|^2))^{3/2}.$$
Similarly, the Gauss curvature $K$ is
$$K=-\frac{f_{xx}f_{yy}-f_{xy}^2}{(1-f_x^2-f_y^2)^2}.$$
\end{example}

In order to show the difference between the spacelike and timelike case, we study a family of timelike surfaces whose Weingarten endomorphism is not diagonalizable.

\begin{example}\label{scroll}
 Let  $\alpha:I\rightarrow\e_1^3$ be a lightlike curve and we denote by
$\{\t,\n,\b\}$ the Frenet trihedron. We consider the map  $$\Psi:I\times\r\rightarrow\e_1^3,\ \ \Psi(s,t)=\alpha(s)+t \b(s).$$
This surface is called a $B$-scroll and it was introduced by Graves in \cite[p. 374]{gr}. We compute the matrix of the Weingarten map   with respect to the basis   $\{\Psi_s,\Psi_t\}$. Since  $\Psi_s=\t+t\b'=\t+t\tau\n$ and $\Psi_t=\b$, then
$$\left(\begin{array}{cc}E&F\\ F&G\end{array}\right)=\left(\begin{array}{cc}t^2\tau^2&-1\\ -1 & 0\end{array}\right).$$
As the determinant is negative, the surface is timelike. Since
$$\Psi_{ss}=t\tau^2\t+(1+t\tau')\n+t\tau \b,\ \ \Psi_{st}=\tau\n,\ \ \Psi_{tt}=0,$$
the second fundamental form is
$$\left(\begin{array}{cc}e&f\\ f&g\end{array}\right)=\left(\begin{array}{cc}
-1-t(\tau'-t\tau^3)&-\tau\\ -\tau&0\end{array}\right).$$
The Weingarten endomorphism is now
$$A=\left(\begin{array}{cc}\tau&0\\ 1+t\tau'&\tau \end{array}\right).$$
This matrix is not diagonalizable. On the other hand, the mean curvature is $H=\tau$ and the Gauss curvature is  $K=\tau^2$.
Finally, let us see that this surface satisfies  $H^2-\epsilon K=H^2-K=0$ but it is not umbilical.
\end{example}
We give more examples of  surfaces.

\begin{example}(Helicoids) A helicoidal surface is a surface invariant by a uniparametric group of helicoidal motions of $\e_1^3$.   A particular case is a surface   of revolution that will be studied in Sect. \ref{section4} when $H$ is constant. Here our  interest is the computation of $H$ and $K$. We refer \cite{dk,ko,mp,wo}. The next surfaces are all minimal surfaces ($H=0$).
\begin{enumerate}
\item Helicoid of first kind. The parametrization is
$$\Psi(s,t)=(s\cos(t),s\sin(t),ht),\  s>h>0.$$
Since $W=s^2-h^2$, the surface is spacelike with $H=0$ and $K=h^2/(s^2-h^2)^2$.
\item Helicoid of second kind. Consider the surface
$$\Psi(s,t)=(ht,s\cosh(t),s\sinh(t)),\ h>0,\ s\in (h,\infty).$$
Then $W=h^2-s^2<0$ and the surface is timelike. Here we have $H=0$ and $K=h^2/(s^2-h^2)^2$. Since $H^2-K<0$, Proposition \ref{pr-hk} asserts that the surface has no umbilics.
\item Helicoid of third kind. The parametrization is
$$\Psi(s,t)=(ht,s\sinh(t),s\cosh(t)),\ h>0, s\in\r.$$
This surface is timelike with $H=0$ and $K=h^2/(s^2+h^2)^2$. The Weingarten map is not diagonalizable.
\item Cayley's surface. Its parametrization is
$$\Psi(s,t)=\left(st-ht+h\frac{t^3}{3},s+ht^2,st+ht+h\frac{t^3}{3}\right),\ h,s>0.$$
Now $W=-4hs$.  Then the surface is timelike with$H=0$ and $K=1/(4s^2)$.
\end{enumerate}
\end{example}

\begin{example}\label{ex-ru}  Ruled surfaces  is a class of surfaces of interest in Lorentz-Minkowski space. See \cite{dk}. We use some examples of ruled surface to compute $H$ and $K$. In all the next examples, the surface is timelike, the Weingarten endomorphism is not diagonalizable and $H^2-K=0$.
\begin{enumerate}
\item Consider the immersion
$$\Psi(s,t)=(s\cos(t),s\sin(t),s+ht),\ h>0.$$
This surface is timelike with $W=-h^2$. Here $H=1/h$ and $K=1/h^2$. Then $H^2-K=0$ and the Weingarten map is
    $$\left(\begin{array}{cc}\frac{1}{h}&-1\\0&\frac{1}{h}\end{array}\right).$$
\item Let $a\not=0$. Define the surface
$$\Psi(s,t)=(ht,(s+a)\cosh(t)+s\sinh(t),(s+a)\sinh(t)+s\cosh(t)).$$
The surface is timelike because $W=-a^2$. Here $H=K=0$ and the Weingarten map is
    $$\left(\begin{array}{cc}0&-\frac{h}{a}\\0&0\end{array}\right).$$
\item (Parabolic null cylinder) The parametrization is
$$\Psi(s,t)=\left(s+h(-t+\frac{t^3}{3}),ht^2,s+h(t+\frac{t^3}{3})\right),\ h>0.$$
The surface is timelike  with $W=-4h^2$ and $H=K=0$. The Weingarten map is
$$\left(\begin{array}{cc}0&1\\0&0\end{array}\right).$$
\end{enumerate}
\end{example}

We end this section with the description of all umbilical surfaces in $\e_1^3$.

\begin{theorem}\label{t-umbilical}
 The only totally umbilical surfaces in Minkowski space are a plane, a hyperbolic plane or a pseudosphere.
\end{theorem}

\begin{proof} The proof is step-by-step as in Euclidean space (\cite[p.147]{dc}). Consider a coordinate neighbourhood $\Omega\subset\r^2$ and let $\Psi=\Psi(u,v)$ be the corresponding parametrization. Since the Weingarten map is diagonalizable,  there is a   function  $f$ such that
\begin{eqnarray*}
(N\circ \Psi)_u&=&(f\circ \Psi) \Psi_u\\
(N\circ \Psi)_v&=&(f\circ \Psi) \Psi_v.
\end{eqnarray*}
As a consequence, $f$ is smooth.
A differentiation with respect to  $u$ and $v$ yields
$$(f\circ \Psi)_u \Psi_v+(f\circ \Psi) \Psi_{uv}=(f\circ \Psi)_v \Psi_u+(f\circ \Psi) \Psi_{uv}.$$
Thus $(f\circ \Psi)_u=(f\circ \Psi)_v$. This means that  $f$ is a constant function in $\Omega$, namely,
$f\circ \Psi=r$, $r\in\r$. Since the surface is connected, $f\circ \Psi=r$ on $M$.
\begin{enumerate}
\item If $r=0$, then $N_u=N_v$, that is,  $N$ is constant. This shows that the surface is a plane.
\item If $r\not=0$, then  $N_u=r \Psi_u$ and $N_v=r \Psi_v$. In $\Omega$ we define
$$h(u,v)=\Psi(u,v)-\frac{1}{r}(N\circ \Psi)(u,v).$$
It follows that $h_u=h_v=0$ and so,  $h$ is constant. Thus there exists  $p_0\in\e_1^3$ such that
$$\Psi(u,v)-\frac{1}{r}(N\circ \Psi)(u,v)=p_0.$$
Then
$$\langle \Psi-p_0,\Psi-p_0\rangle=\mp\frac{1}{r^2},$$
with the sign depending if the surface is spacelike or timelike, respectively. In all these cases, and according to the definition, the surface is included in a hyperbolic plane or a pseudosphere.

\end{enumerate}
\end{proof}

\subsection{Surfaces in Minkowski space  with   mean  curvature  and   Gauss curvature both constant}

A classical result in Euclidean space asserts that the only surfaces  with constant principal curvatures are open sets of umbilical surfaces (planes or spheres) or  right circular cylinders. As far as the author knows, the statement  appears as an exercise  in \cite[p. 263]{on1}) and a proof lies in \cite{mr}. We extend this result to Minkowski space following, if possible, the same proof steps.

For a   non-degenerate surface in $\e_1^3$ with diagonalizable Weingarten map $A$, the constancy of $H$ and $K$ is equivalent to say that the principal curvatures  $\lambda_1$ and $\lambda_2$ are constant. In such case, and with the corresponding changes, the proof in $\e_1^3$ is similar as in Euclidean space, concluding that the surface is  umbilical or a right circular cylinder.

When the Weingarten map is not diagonalizable (necessarily for timelike surfaces), then $H^2-K\leq 0$. An example of this situation appeared  in Ex. \ref{ex-ru}. Consider the surface $M$ parameterized by  $\Psi(s,\theta)=(s\cos(\theta),s\sin(\theta),s+h\theta)$, $h>0$. This surface is  helicoidal, that is, it is  invariant under a uniparametric group $G_{L,h}=\{\phi_\theta: \theta\in\r\}$ of motions of $\e_1^3$, whose axis $L$ is determined by the vector $(0,0,1)$  and a pitch $h$. Here $\phi_\theta$ is
\begin{equation}\label{htime}
\phi_\theta(x,y,z)=\left(\begin{array}{ccc}
\cos{(\theta)}&-\sin{(\theta)}&0\\ \sin{(\theta)}&\cos{(\theta)}&0\\
0& 0&1\end{array}\right)\left(\begin{array}{c}x\\y\\z\end{array}\right)+h\left(\begin{array}{c}0\\0\\ \theta\end{array}\right).
\end{equation}
See \cite{dk,mp}. If we consider the curve $\gamma(t)=(s,0,s)$, then  $M=G_{L,h}(\gamma)=\{\phi_\theta(\gamma(t)):s,\theta\in\r\}$. Recall that $M$ is a not umbilical timelike surface and it is also ruled with lightlike rulings. Moreover, $H$ and $K$ are constant with  $H=1/h$,  $K=1/h^2$ and $H^2-K=0$ on $M$.

Returning to the case that   $H$ and $K$ are both constant functions, we will show that necessarily  $H^2-K\geq 0$. Finally, a recent work by Clelland shows that a surface with $H^2-K=0$ and with non necessarily constant $H$ and $K$, is a ruled surface with both directrix and rulings of lightlike type \cite{cl}. In the case that $H$ and $K$ are constant, we characterize the parametrizations of these surfaces and we present explicit  examples.

\begin{theorem}\label{t-hk} Let $M$ be a non-degenerate surface in Minkowski space $\e_1^3$. If $H$ and $K$ are both constant, then $M$ is an umbilical surface,  a right circular cylinder or a  ruled surface with  directrix and rulings both lightlike. In the latter case, $H^2-K=0$.
\end{theorem}

In general, a hypersurface with constant coefficients
of the characteristic polynomial of the Weingarten map is called isoparametric. For surfaces, these coefficients are, up a constant, $H$ and $K$. In two-dimensional case, these surfaces have been locally
classified in Euclidean space and in Lorentzian space of arbitrary dimension. The case of one or two distinct (real) eigenvalues is treated in \cite[Th. 5.1]{aky} and  the case of a non-diagonalizable Weingarten map is treated in \cite[Th. 4.5]{ma}.

\begin{proof}
Since the result is local, we always consider a parametrized open set of the surface.
 The proof is divided in two cases depending if the Weingarten map is or is not diagonalizable.

\noindent \emph{Case I. The Weingarten endomorphism is diagonalizable.}

The proof follows similar steps as in Euclidean space.  This case  occurs  if the surface is spacelike of if it is timelike with   $H^2-K>0$. By the expressions of $H$ and $K$ in terms of the principal curvature, we deduce that $\lambda_i$ are also constant. If at some point $p\in M$, $\lambda_1(p)=\lambda_2(p)$, then $H(p)=\epsilon\lambda_1(p)$ and $K(p)=\epsilon\lambda_1(p)^2$. Because $H$ and $K$ are constant functions, then  $H^2-\epsilon K=0$ on $M$. This means that the surface is umbilical proving a part of the statement of Th.  \ref{t-hk}.

The other case is that there are not umbilical points. We show that the surface is a right circular cylinder. Since in an open of a non-umbilic point there are local coordinates by  lines of curvature, let $\Psi:U\subset\r^2\rightarrow M\subset\e_1^3$ be a parametrization of the surface with the property $F=f=0$ and so,
$$-N_u=\lambda_1 \Psi_u,\ \ -N_v=\lambda_2 \Psi_v.$$
If the surface is spacelike, both vectors $\Psi_u$ and $\Psi_v$ are spacelike. If $M$ is timelike, we assume, without loss of generality, that $\Psi_u$ is spacelike and $\Psi_v$ is timelike.
In these coordinates, the expressions of the second derivatives of the parametrizations in terms of $\Psi_u$, $\Psi_v$ and $N$ are:
\begin{eqnarray}\label{x1}
\Psi_{uu}&=&\frac{E_u}{2E}\Psi_u+eN.\nonumber\\
 \Psi_{uv}&=&\frac{E_v}{2E}\Psi_u+\frac{G_u}{2G}\Psi_v.\\
\Psi_{vv}&=&\frac{G_v}{2G}\Psi_v+gN.\nonumber
\end{eqnarray}
We differentiate the first equation with respect to $v$, the second one with respect to $u$ and taking into account that $\lambda_i$ are constants, we have $-N_{uv}=\lambda_1 \Psi_{uv}=\lambda_2 \Psi_{uv}$. Since $\lambda_1\not=\lambda_2$, we obtain $N_{uv}=\Psi_{uv}=0$ on the domain $U$. Thus
$$E_v=\langle \Psi_u,\Psi_u\rangle_v=2\langle \Psi_{uv},\Psi_u\rangle=0.$$
$$G_u=\langle \Psi_v,\Psi_v\rangle_u=2\langle \Psi_{uv},\Psi_v\rangle=0.$$
Then $E$ is a function depending only on $u$ and $G$ depends only on $v$. Similarly, and using $N_{uv}=\Psi_{uv}=0$ again, we obtain
$$e_v=-\langle N_u,\Psi_u\rangle_v=-\langle N_{uv},\Psi_u\rangle-\langle N_u,\Psi_{uv}\rangle=0.$$
$$g_u=-\langle N_v,\Psi_v\rangle_u=-\langle N_{uv},\Psi_v\rangle-\langle N_v,\Psi_{uv}\rangle=0.$$
We conclude that $e$ is a function depending only on $u$ and $g$ depends only on $v$. Hence, and using that $(\Psi_{uu})_v=(\Psi_{uv})_u$ and $(\Psi_{vv})_u=(\Psi_{uv})_v$, we have from (\ref{x1})
\begin{equation}\label{eg}
eN_v=0,\ \ gN_u=0.
\end{equation}
A first observation is that  both $N_u$ and $N_v$ do not vanish at some point $p$, because in such a case,   the two principal curvatures vanish at $p$ and it  would be an umbilic point. Also, both $e$ and $g$ can not vanish at some point because $\lambda_1=e/E $ and $\lambda_2=g/G$ and the point would be umbilical again. Without loss of generality, we suppose that $N_u\not=0$ in an open of the surface. From (\ref{eg}),   $g=0$ on $U$ and using the second equation of (\ref{eg}), $e\not= 0$ and $N_v=0$.

{\it Claim:} The vector field $\Psi_v/|\Psi_v|$ is constant on the domain $U$.

Proof of the claim. We show that $(\Psi_v/|\Psi_v|)_u=(\Psi_v/|\Psi_v|)_v=0$ on $U$. We begin with
$(\Psi_v/|\Psi_v|)_u$. Recall that the sign of $G$ is the same that the one of $-\epsilon$. Because $\Psi_{uv}=0$ and $|\Psi_v|=\sqrt{-\epsilon G}$ depends only on $v$, we have
$$\Big( \frac{\Psi_v}{|\Psi_v|}\Big)_u=\frac{\Psi_{uv}}{|\Psi_v|}-\Big(\frac{1}{\sqrt{-\epsilon G}}\Big)_u \Psi_v=0.$$
Now we prove that $(\Psi_v/|\Psi_v|)_v=0$ showing that this vector field is orthogonal to each one of the vectors $\{\Psi_u,\Psi_v,N\}$. Using $F=0$ and $\Psi_{uv}=0$, we have
$$\langle \Big( \frac{\Psi_v}{|\Psi_v|}\Big)_v,\Psi_u\rangle=\langle\frac{\Psi_v}{|\Psi_v|},\Psi_u\rangle_v-\langle \frac{\Psi_v}{|\Psi_v|},\Psi_{uv}\rangle=0.$$
Since $G/\sqrt{-\epsilon G}=-\epsilon\sqrt{-\epsilon G}$, we get
$$\langle \Big( \frac{\Psi_v}{|\Psi_v|}\Big)_v,\Psi_v\rangle=\langle\frac{\Psi_v}{|\Psi_v|},\Psi_v\rangle_v-\langle \frac{\Psi_v}{|\Psi_v|},\Psi_{vv}\rangle=(-\epsilon\sqrt{-\epsilon G})_v-\frac{G_v}{2\sqrt{-\epsilon G}}=0.$$
By taking into account that $N_v=0$,   we have
$$\langle \Big( \frac{\Psi_v}{|\Psi_v|}\Big)_v,N\rangle=-\langle \frac{\Psi_v}{|\Psi_v|},N_v\rangle=0.$$
This finishes the proof of the claim.

The claim assures the existence of a unit vector $\vec{a}$ such that
$$\vec{a}=\frac{\Psi_v}{|\Psi_v|}=\frac{\Psi_v}{\sqrt{-\epsilon G}}.$$
In particular,
\begin{eqnarray}
& & \langle \Psi_u,\vec{a}\rangle=\langle N,\vec{a}\rangle=0.\label{81}\\
& & \langle \Psi_v,\vec{a}\rangle=\langle \Psi_v,\frac{\Psi_v}{|\Psi_v|}\rangle=-\epsilon\sqrt{-\epsilon G}.\label{82}
\end{eqnarray}
We remark that
$$\langle \vec{a},\vec{a}\rangle=\langle\frac{\Psi_v}{\sqrt{-\epsilon G}},\frac{\Psi_v}{\sqrt{-\epsilon G}}\rangle=-\epsilon.$$
Define
$$h:U\rightarrow\e_1^3, \ \ h(u,v)=\Psi(u,v)+\epsilon\langle \Psi(u,v),\vec{a}\rangle\vec{a}+\frac{1}{\lambda_1}N(u,v).$$
We prove that $h$ is a constant function showing that the partial derivatives $h_u$ and $h_v$ vanish on $U$. Using that $N_u=-\lambda_1 \Psi_u$ and (\ref{81}), we have
$$h_u=\Psi_u+\frac{1}{\lambda_1}N_u=0.$$
By using (\ref{82}) and   $|\Psi_v|^2=-\epsilon G$, we have
$$h_v=\Psi_v+\epsilon\langle \Psi_v,\vec{a}\rangle\vec{a}=0.$$
As $h$ is constant, there exists $p_0\in\e_1^3$ such that
$$\Psi(u,v)+\epsilon\langle \Psi(u,v),\vec{a}\rangle\vec{a}+\frac{1}{\lambda_1}N(u,v)=p_0.$$
Hence, together (\ref{81}) and (\ref{82}),  we have
$$|\Psi(u,v)-p_0|^2+\epsilon\langle \Psi(u,v)-p_0,\vec{a}\rangle^2=\frac{\epsilon}{\lambda_1^2}.$$
By  Ex. \ref{ex-cy}, the surface lies included in the right circular cylinder  of axis $L$ and  radius $1/|\lambda_1|$, where $L$ is the straight-line passing $p_0$ with direction $\vec{a}$.

\noindent \emph{Case II. The Weingarten endomorphism is not diagonalizable.}

Now the surface is timelike. Assume that $H^2-K<0$ and we will arrive to a contradiction. We take null coordinates on $M$, that is,
let $\Psi:U\rightarrow M$ be a parametrization  such that $E=G=0$ (\cite{we}). From (\ref{hk}) and since $H$ and $K$ are constant,  there exists   $\lambda,\mu\in\r$ such that
$$H=\frac{f}{F}=\lambda,\ \ K=\frac{eg-f^2}{-F^2}=\mu.$$
In particular,
\begin{equation}\label{eg2}
eg=(\lambda^2-\mu)F^2.
\end{equation}
The matricial expression of the Weingarten map $A$ with respect to $\{\Psi_u,\Psi_v\}$ is
$$A=\left(\begin{array}{ll}0&F\\ F&0\end{array}\right)^{-1}\left(\begin{array}{ll}e&f\\ f&g\end{array}\right)=
\frac{1}{F}\left(\begin{array}{ll} f&e\\ g&f\end{array}\right).$$
Because $H^2-K<0$, then $eg\not=0$.   The second derivatives of $\Psi$ are
$$\Psi_{uu}=\frac{F_u}{F}\Psi_u+eN,\ \ \Psi_{uv}= fN,\ \ \Psi_{vv}=\frac{F_v}{F}\Psi_v+gN.$$
By combining these equations,   $(\Psi_{uu})_v=(\Psi_{uv})_u$, $(\Psi_{vv})_u=(\Psi_{uv})_v$ and together the expression of $A$ we get
\begin{eqnarray}
 \left(\frac{F_u}{F}\right)_v-\frac{e^2}{F}=-\frac{f^2}{F}, & & \left(\frac{F_v}{F}\right)_u-\frac{g^2}{F}=-\frac{f^2}{F}.\label{e1}\\
 \frac{ef}{F}& = &\frac{fg}{F}.\label{e2}\\
 \frac{F_u}{F}f+e_v=
 f_u, & &  \frac{F_v}{F}f+g_u=f_v.\label{e3}
\end{eqnarray}
Equation (\ref{e1}) gives  $e^2=g^2$, that is, $e=\pm g$.
We distinguish two possibilities:
\begin{enumerate}
\item Case $f=0$. Equation (\ref{e3}) yields $e_u=e_v=0$. Thus $e$ and $g$ are constant functions. From (\ref{eg2}), the coefficient $F$ is constant, and (\ref{e1})  concludes $e^2=f^2=0$: contradiction.
\item Case $f\not=0$. Equation (\ref{e2}) implies  $e=g$, and  (\ref{eg2}) that $F$ is constant. Equation (\ref{e1})  assures $e^2=f^2$, that is, $e=\pm f$.
In particular, $K=0$ and  $H^2-K\geq 0$: contradiction. In fact,   the matricial expression of the Weingarten map would be
$$\frac{e}{F}\left(\begin{array}{ll}
\pm 1&1\\1&\pm 1\end{array}\right),$$
which is diagonalizable with eigenvalues $0$ and $\pm e/F$.
\end{enumerate}
Suppose $H^2-K=0$ on the surface. From \cite{cl},  the surface is ruled and $M$ can parametrized in such way that the directrix and the rulings are both lightlike. This concludes the proof.
\end{proof}

We give examples of non-umbilical timelike surfaces with   $H$ and $K$ both constant functions and $H^2-K=0$.   The surface $M$   parametrizes as  $\Psi(s,t)=\alpha(s)+t \vec{w}(s)$, where $\alpha$ is a lightlike curve,  and the rulings $\vec{w}(s)$ are also lightlike (\cite{cl}). See also Ex. \ref{scroll}.
Now $W=-F^2$. Then
$$H=\frac{f}{F}=\frac{\mbox{det}(\Psi_u,\Psi_v,\Psi_{uv})}{F|F|}=\lambda,\ \ K=\frac{f^2}{F^2}=\frac{{\mbox{det}(\Psi_u,\Psi_v,\Psi_{uv})}^2}{F^4}=\mu,$$
for some constants $\lambda,\mu\in\r$. Parametrize $\alpha(s)=(x(s),y(s),z(s))$. Since $\alpha$ is lightlike,  Prop. \ref{pr-twof} asserts that the curve $\alpha$ is graph of two functions with respect to the $z$-axis: $\alpha(s)=(x(s),y(s),s)$. Moreover, $x'(s)^2+y'(s)^2-1=0$ and then,
\begin{equation}\label{xy}
x'(s)=\cos\psi(s),\ \ y'(s)=\sin\psi(s),
\end{equation}
for some differentiable function $\psi$. We distinguish different cases.

\begin{enumerate}
\item Suppose that $\vec{w}$ is a constant vector field. Then $f=0$ and so, $H=K=0$. This means that any cylinder with basis a  lightlike curve and with lightlike rulings is a surface with $H=K=0$. In fact the result is  more general: any cylindrical surface with lightlike rulings is a timelike surface with $H=K=0$. This is because from the parametrization $\Psi(s,t)=\alpha(s)+t\vec{w}$,  $\Psi_{t}=\vec{w}$ and
$\Psi_{st}=\Psi_{tt}=0$. Thus $G=f=g=0$.

\item Suppose   that  $\vec{w}(s)$ is a horizontal curve. Without loss of generality, assume that $w$ lies in the plane of equation $z=1$. Since $\vec{w}(s)$ is lightlike, $\vec{w}(s)=(\cos\theta(s),\sin\theta(s),1)$ for some differentiable function $\theta$. After a new reparametrization of the surface by $\bar{s}=\psi\circ\theta^{-1}(s)$, the parametrization of $M$ is $\Psi(s,t)=\alpha(s)+t(\cos{(s)},\sin{(s)},1)$, together the equations (\ref{xy}).  Now
$$f=1-\cos{(s-\psi)},\ \ F=-1+\cos{(s-\psi)}.$$
The condition $f=\lambda F|F|$ gives
$$1-\cos{(s-\psi)}=-\lambda(1-\cos{(s-\psi)})^2.$$
If $\lambda=0$, then $1-\cos{(s-\psi)}=0$ and so, $F=0$: contradiction. Thus  $\lambda\not=0$. Then $1=-\lambda(1-\cos{(s-\psi)})=0$. This means that $\psi(s)=s+c$ for some $c\in\r$,
$c\not= 2n\pi$, $n\in \mathbb{Z}$. The solution of (\ref{xy}) is
$$x(s)=\sin{(s+c)}+x_0,\ \ y(s)=-\cos{(s+c})+y_0,$$
with $x_0,y_0\in\r$ integration constants. After a horizontal translation, the surface parametrizes as
\begin{equation}\label{heli}
\Psi(s,t)=(t\cos{(s)}+\sin{(s+c)},-\cos{(s+c)}+t\sin{(s)},s+t).
\end{equation}
This surface appears in \cite[Sect. 5]{cl} and it is a helicoidal surface whose axis $L$ is the $z$-axis  and pitch $h=1$.  For this, let
$\phi_m\in G_{L,1}$, $m\in\r$, whose expression is given by \eqref{htime}. Then one can easily prove that
$\phi_m(\Psi(s,t))=\Psi(m+s,t)$. By \cite[Lemma 2.2]{mp}, the generating curve $\gamma(t)$ lies in the plane of equation $y=0$. In
the particular case  $c=\pi/2$, then $\Psi(s,t)=((1+t)\cos{(s)},(1+t)\sin{(s)},s+t)$, which it is a helicoidal surface with timelike axis and the generating curve is the lightlike curve $\gamma(t)=(t+1,0,t)$. Up a rigid motion of $\e_1^3$, this
  surface  appeared in Ex. \ref{ex-ru}.
\item The rest of surfaces are obtained when  $\vec{w}(s)$ is not a planar horizontal curve. We parametrize  $\vec{w}(s)=(a(s),b(s),s))$. Since $\vec{w}$ is lightlike, there exists $\theta(s)$ such that
$\vec{w}(s)=s(\cos\theta(s),\sin\theta(s),1)$. We suppose that $s>0$. Then
$$f=s^2\theta'(1-\cos(\psi-\theta))\ \ F=s(-1+\cos(\psi-\theta)),$$
and we have
\begin{equation}\label{general}
\theta'(1-\cos(\psi-\theta))=-\lambda(1-\cos(\psi-\theta))^2,\ \ \lambda\in\r.
\end{equation}
Since $1-\cos(\psi-\theta)=0$ implies  $F=0$, we conclude
$$\frac{\theta'}{1-\cos{(\psi-\theta)}}=-\lambda.$$
\end{enumerate}

\begin{theorem} The only non-umbilical timelike surfaces in $\e_1^3$ with $H$ and $K$ constant and $H^2=K=\lambda$  are:
\begin{enumerate}
\item Cylinders with lightlike  rulings. These surfaces have $H=K=0$.
\item Helicoidal surfaces parametrized by (\ref{heli}). These surfaces are helicoidal and satisfy $H,K\not=0$.
\item Ruled surfaces $X(s,t)=\alpha(s)+t s(\cos{\theta(s)},\sin{\theta(s)},1)$, where $\theta$ is a solution of (\ref{general}), $\alpha$ is given by (\ref{xy}) and   $\psi-\theta\not=2 k\pi$, $k\in\mathbb{Z}$. These examples satisfy $H^2=K=\lambda$.
\end{enumerate}
\end{theorem}
We end this section giving explicit examples of surfaces in the third case. Take  $\psi(s)=n\theta(s)$, with $n\in\mathbb{Z}$.
\begin{example}\label{ex1}
 \begin{enumerate}
 \item If $n=2$, the solution is, up a reparametrization,  $\theta(s)=-2\mbox{ arccot}{(c s)}$, $c\in\r$. From (\ref{xy}) we have
$$\alpha(s)=\left(s+\frac{4s}{1+c^2s^2}-4\frac{\arctan{(c s)}}{c},-\frac{4}{c+c^3 s^2}-2\frac{\log{(1+c^2s^2)}}{c} ,s\right).$$
\item If $n=3$, the solution of (\ref{general}) is $\theta(s)=-\mbox{arccot}{(2c s)}$ and the integration of (\ref{xy}) leads to
$$\alpha(s)=\left(\frac{5+4s^2}{2c\sqrt{1+4c^2s^2}},\frac{4s}{\sqrt{1+4c^2s^2}}-\frac{3}{2c}
\log\Big(2c s+\sqrt{1+4c^2s^2}\Big),s\right).$$
\end{enumerate}
\end{example}


\section{Spacelike surfaces with constant mean curvature}\label{section4}

In this section we study spacelike surfaces with constant mean curvature (cmc). We will obtain examples of such surfaces with some added geometric assumption, as for example, that the surface is rotational or  a graph.

As in Euclidean space, a cmc spacelike surface is a critical point of the area functional. Let $M$ be a compact surface and let $x:M\rightarrow\e_1^3$ be a spacelike immersion. In particular, $\partial M\not=\emptyset$. A variation of $x$ is a differentiable map $X:(-\epsilon,\epsilon)\times M\rightarrow\e_1^3$ with the following properties:
 \begin{enumerate}
 \item $X(0,p)=x(p)$ for all $p\in M$.
 \item The maps $X_t:M\rightarrow\e_1^3$,  $X_t(p)=X(t,p)$, are spacelike immersions for all $t$.
  \item $X_t$ pointwise fixes the boundary of $\partial M$, that is,  $X(t,p)=x(p)$ for $p\in\partial M$.
  \end{enumerate}
  We define the variational vector field of the variation $X$ as
$$\xi(p)=\frac{\partial X}{\partial t}(0,p).$$

The area functional is
$$A(t)=\int_M dA_t,$$
where $dA_t$ is the area element on $M$ induced by the metric $X_t^{*}(\langle,\rangle)$. Similarly, the volume functional is
$$V(t)=\frac13\int_M\langle X_t,N_t\rangle dA_t,$$
where $N_t$ is a unit normal vector field to the immersion $X_t$. The value  $V(t)$ represents the enclosed volume by the cone whose
basis is $X_t(M)$ and whose vertex is the origin of coordinates.  Because $\e_1^3=T_pM\oplus T_pM^\bot$, the vector field $\xi$ decomposes as
$$\xi(p)=\xi(p)^T+f(p)N(p),\ \ f\in C^{\infty}(M).$$
For the tangent part  $\xi(p)^T$, there exists   $ \tilde{\xi}\in\mathfrak{X}(M)$ such that $(dx)_p \tilde{\xi}(p)=\xi(p)^T$.
The map $A(t)$ is differentiable at $t=0$ with
$$A'(0)=2\int_M (Hf) dA-\int_{\partial M}\langle \tilde{\xi},\nu\rangle\ ds,$$
where $H$ is the mean curvature of the immersion $x$ and $\nu$ is the unit conormal vector field along  $\partial M$.
Since the variation $X$  fixes the boundary $\partial M$, the second summand vanishes.

\begin{theorem} Let $x:M\rightarrow\e_1^3$ be a spacelike immersion of a compact surface $M$ and  let $X$ be a variation of $x$. The first variation of the area is
$$A'(0)=2\int_M (Hf)\ dA.$$
For the volume functional and for variations that preserve the boundary, we have
$$V'(0)=-\int_M f\ dA.$$
\end{theorem}
Hence it is immediate,

\begin{theorem}\label{41t1} Let $M$ be a compact surface and let  $x:M\rightarrow\e_1^3$ be a spacelike immersion in  $\e_1^3$. Then $x$ has constant mean curvature if and only if it is a critical point of the area functional for all volume preserving variation of $x$ that fixes the boundary.
\end{theorem}

First examples of cmc spacelike surfaces are the umbilical surfaces. Then a spacelike plane and a hyperbolic plane have constant mean curvature. More surfaces appeared in Ex. \ref{ex-cy}, where Lorentzian and hyperbolic cylinders have constant mean curvature. In the rest of this section, we give more examples with some added geometric assumption.

\subsection{Translation surfaces}

A translation surface is a surface that is the graph of a function of type  $z=f(x)+g(y)$, where $x\in I$, $y\in J$. The name of translation surface is motivated by its parametrization  $\Psi(x,y)=(x,y,f(x)+g(y))$, where the surface $S$ is viewed as the sum of two planar curves in orthogonal coordinate planes, namely, $C_x:x\longmapsto (x,0,f(x))$ and $C_y:y\longmapsto (0,y,g(y))$ and thus,
$$S=\bigcup_{x\in I} T_x(C_y), \ T_x(p)=p+(x,0,f(x)), p\in\e_1^3.$$
If $H=0$,   the only translation surface in $\e^3$ is a plane and the Scherk surface
$$z=\frac{1}{\lambda} \log{\left|\frac{\cos(\lambda y+a)}{\cos(\lambda x+b)}\right|},\ \lambda, a,b\in\r, \lambda\not=0.$$
If $H\not=0$ is constant, then the surface is a right circular cylinder \cite{liu}.

In order to ask a similar problem in $\e_1^3$, there are some differences.
\begin{enumerate}
\item First, the surface may be considered as the graph of a function defined in each one of the coordinate planes of $\r^3$, indicating that the problem depends on the choice of this plane.
    \item We may also assume that this plane is lightlike since the coordinate planes are not lightlike.
    \item The surface may be spacelike or timelike.
\end{enumerate}

\begin{theorem}\label{te-43}  Consider a spacelike translation surface in $\e_1^3$ given by $z=f(x)+g(y)$. If $H$ is constant, then the surface is
a plane, a Scherk surface or a cylinder.
\end{theorem}

\begin{proof} We do the proof only for $H=0$ (see \cite{liu} when $H$ is a non-zero constant). The mean curvature $H$ satisfies
$$-2H(1-f'^2-g'^2)^{3/2}=(1-f'^2)g''+(1-g'^2)f''.$$
Let $H=0$. This means
$$\frac{f''}{1-f'^2}=-\frac{g''}{1-g'^2}=\lambda$$
for some constant $\lambda\in\r$. If $\lambda=0$, then $f$ and $g$ are linear and the surface is a plane. If $\lambda\not=0$, an integration gives
$$f(x)=\frac{1}{\lambda}\log(\cosh(\lambda x+a)),\ \ g(y)=\frac{1}{\lambda}\log(\cosh(\lambda y+b)),\ \ a,b\in\r.$$
Thus the Lorentzian Scherk surface  is
$$z=\frac{1}{\lambda} \log{\left|\frac{\cosh(\lambda y+a)}{\cosh(\lambda x+b)}\right|},\ \lambda, a,b\in\r, \lambda\not=0.$$
\end{proof}

In Th. \ref{te-43}, the surface is a graph on a spacelike plane. We analyse the other two cases.
\begin{enumerate}
\item Suppose  that the surface is a graph on a timelike plane. After a rigid motion, we suppose that the plane is the $xz$-plane. A parametrization of the surface is $\Psi(x,z)=(x,f(x)+g(z),z)$. The computation of the mean curvature $H$ gives
$$2H(-1-f'^2+g'^2)^{3/2}=-(1+f'^2)g''+(g'^2-1)f''.$$
If $H=0$, then
$$\frac{f''}{1+f'^2}=\frac{g''}{g'^2-1}=\lambda,$$
for some constant $\lambda$. Again, if $\lambda=0$, the solution is a plane. If $\lambda\not=0$, then
$$f(x)=\frac{1}{\lambda}\log(\cos(\lambda x+a)).$$
The integration of $g$ is a bit different because  the condition  $g'^2>1$. So,  the integral $\int\frac{1}{x^2-1}\ dx$  depends on the sign of $x^2-1$. If it is negative, then the solution is $\mbox{arc}\tanh(x)$, but in our case is positive. This means that we can not do
$$\int\frac{1}{x^2-1}\ dx=-\int\frac{1}{1-x^2}\ dx=-\mbox{arc}\tanh(x)$$
because the domain of the integrand function is $x^2>1$ and the domain of $\mbox{arc}\tanh(x)$ is $x^2<1$. However, it is not difficult to find that $g(z)=-\frac{1}{\lambda}\log(\sinh( \lambda z+b)),\,b\in\r$. Thus
$$y=\frac{1}{\lambda}\log\left|\frac{\cos(\lambda x+a)}{\sinh(\lambda z+b)}\right|.$$
\item Suppose now that the surface is the graph of the functions $f(u)+g(v)$ over a lightlike plane. After an isometry, we assume that this plane is of equation $y+z=0$. Then a parametrization of the surface is
\begin{eqnarray*}
\Psi(u,v)&=&u(1,0,0)+v(0,1,-1)+(f(u)+g(v))(0,1,1)\\
&=&(u,v+f(u)+g(v),-v+f(u)+g(v)).
\end{eqnarray*}
Here the coefficients of the first fundamental form are $E=1$, $F=0$ and $G=4g'$ and the mean curvature $H$ satisfies:
$$-8Hg'^{3/2}=4f''g'+g''.$$
If $H=0$, then $g''/g'=-4f''=a$, for some constant $\lambda\in\r$. If $\lambda=0$, then $f$ and $g$ are linear and the surface is a plane. If $\lambda\not=0$, we obtain by integration:
$$f(u)=-\frac{\lambda}{4}u+a,\ \ g(v)=\frac{1}{\lambda}\ e^{\lambda v+b}+c,\ \ a,b,c\in\r.$$
We observe that the surface is non-degenerate if and only if $g'\not=0$ and the the sign of $g'$ gives the spacelike or timelike character to the surface.
\end{enumerate}

\subsection{Rotational surfaces}

An important class of surfaces are the rotational surfaces, called also surfaces of revolution.
\begin{definition}\label{def1}
A surface in Lorentz-Minkowski space $\e_1^3$ is a surface of revolution with respect to a straight-line $L$, called the rotational axis, if it invariant by the  uniparametric family of boosts $G_L=\{\phi_\theta: \theta\in\r\}$.
\end{definition}
 Therefore there are three types of rotational surface according to the causal character of $L$.

 Sometimes in the literature appears other  definition of a surface of revolution.

\begin{definition}[second option]\label{def2}
A surface of revolution of axis $L$ is  the surface obtained when we apply a uniparametric group $G_L$ of boosts to a plane curve contained in a plane through the axis $L$.
\end{definition}

Of course, a rotational surfaces according to  Def. \ref{def2} is  a surface of revolution in the sense of Def. \ref{def1} but in general, both definitions are not equivalent.  We analyse what happens in Euclidean space. Given a surface of revolution $S$ in $\e^3$, it is easy to prove that any plane containing the axis is $L$ intersects transversally $S$, obtaining a plane curve $C$. If one rotates $C$ about $L$ then one obtains the initial surface.

In Lorentz-Minkowski space this is partially true and it depends if the surface is spacelike or timelike. Exactly, both definitions are equivalente if the surface is spacelike or if the surface is timelike and the axis $L$ is timelike or lightlike. The problem appears for a timelike surface with spacelike axis.

In order to indicate where lies the problem, consider the  pseudosphere $\s_1^2$ of equation $x^2+y^2-z^2=1$, which is a timelike surface. This surface is invariant by the group of boosts with axis $\mbox{Span}\{E_1\}$. Indeed, if $(x,y,z)\in S$, $\phi_\theta(x,y,z)=(x,\cosh(\theta)y+\sinh(\theta)z,\sinh(\theta)y+\cosh(\theta)z)$, which satisfies the equation of the pseudosphere. On the other hand, the intersection of  $S$ with the $xz$-plane is a curve with two components. Let $\alpha(s)=(\cosh(s),0,\sinh(s))$ be the parametrization of the component corresponding with $x>0$. If we reverse the process by Def. \ref{def2}  and we rotate $\alpha$ with respect to the group $G_L$, we obtain
$$S_1= \left\{(\cosh(s),\sinh(\theta)\sinh(s),\cosh(\theta)\sinh(s)):\theta,s\in\r\right\}.$$
The surface $S_1$ is only a part of $S$ and it does not cover fully the pseudosphere. This contrast to with the Euclidean space. Even, if one considers the component of $S\cap \{y=0\}$ corresponding with $x<0$ and rotates it about $G_L$, we obtain a surface $S_2\subset S$, but the union $S_1\cup S_2$ does not coincide with $S$. In fact, we also need to intersect $\s_1^2$ with the plane of equation $x=0$, obtaining two curves and rotate with respect to $L$, giving two surfaces $S_3$ and $S_4$. It is now when $S=S_1\cup S_2\cup S_3\cup S_4$. A discussion can seen in \cite{bco}.

The  difference when one considers a spacelike surface is that  all curves in the surface are spacelike, imposing strong restrictions. Exactly, let $S$ be a spacelike surface of revolution with axis $L=\mbox{Span}\{E_1\}$ according to Def. \ref{def1}. Then $S$ can not intersect the set $A=\{(x,y,z)\in\r^3: y^2=z^2\}$ because in such a case, if $(x,y,\pm y)\in S\cap A$, the orbit of this point by the group $G_L$ is a lightlike curve. Then $S$ is included in one of the four components of $\r^3\setminus A$. In fact, $S$ does not intersect the set of $\r^3$ given by $y^2-z^2>0$. The argument is that the orbit of a point  $(x,y,z)$ with $y^2>z^2$ is a timelike curve. As conclusion,
$S$ must be included in the domain $\{(x,y,z)\in\r^3: y^2-z^2<0, z>0\}$ or in the domain $\{(x,y,z)\in\r^3: y^2-z^2<0, z<0\}$. Thus, when one intersects $S$ with the plane $y=0$, the curve generates the surface in the sense of Def. \ref{def2}.

We study rotational cmc spacelike surfaces. By simplicity, we only consider the case that the axis is timelike which, after a rigid motion, we assume that it is $\mbox{Span}\{E_3\}$. Then the surface is a surface of revolution in the Euclidean sense. Suppose that the generating curve $\alpha$ is contained in the $xz$-plane. Since the surface is spacelike, the curve $\alpha$ is spacelike and thus, it is a graph on the $x$-axis. Let $\alpha(s)=(s,0,f(s))$, $s\in I$. A parametrization of $S$ is $\Psi(s,\theta)=(s\cos(\theta),s\sin(\theta),f(s))$, $s\in\theta$.
Then the mean curvature $H$ satisfies:
$$s^3 f''+s^2 f'(1-f'^2)=-2Hs^3(1-f'^2)^{3/2},$$
which writes as
$$\left(\frac{sf'(s)}{\sqrt{1-f'(s)^2}}\right)'=(-Hs^2)'.$$
Then there exists $\lambda\in\r$ such that
$$ \frac{sf'(s)}{\sqrt{1-f'(s)^2}}=-Hs^2+\lambda.$$
This is similar as in Euclidean space.  See \cite{hn}.

Consider   $H=0$. Recall that in $\e^3$, besides planes, the only minimal rotational surface is the catenoid, whose generating curve is  $z=\lambda\cosh(x/\lambda+b)$, $\lambda,b\in\r, \lambda\not=0$. In Lorentz-Minkowski space, if $H=0$, there exists a constant $\lambda\in\r$ such that
$$\frac{sf'(s)}{\sqrt{1-f'(s)^2}}=\lambda.$$
If $\lambda=0$, then $f$ is constant, obtaining a horizontal plane. If $\lambda\not=0$, the solution is
$$f(s)=\lambda\mbox{ arc}\sinh(\frac{s}{\lambda})+b,\ \ b\in\r.$$
If we see the generating curve as a graph on the $z$-axis, then the curve writes as $x=\lambda\sinh(z/\lambda+b)$. This surface is called the Lorentzian catenoid. Observe that the curve intersects the $z$-axis and at this point, the surface is not regular. This is a difference with the Euclidean case.

Changing the axis to be spacelike or lightlike, we obtain all spacelike catenoids of $\e_1^3$ (\cite{ko}).

\begin{figure}[hbtp]
\begin{center}\includegraphics[width=.4\textwidth]{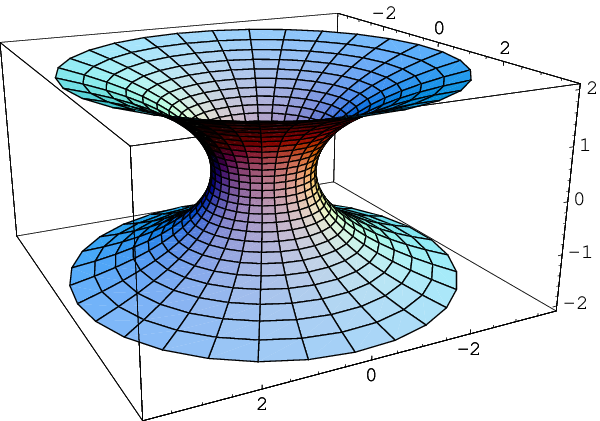},\ \includegraphics[width=.4\textwidth]{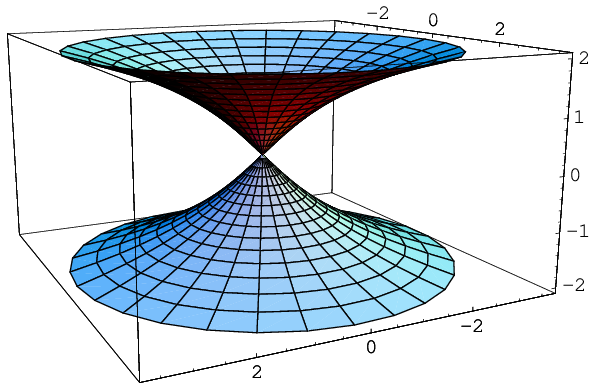}\end{center}
\caption{The catenoid in $\e^3$ (left) and in $\e_1^3$ (right)\label{catenoid}}
\end{figure}

\subsection{Riemann examples}

\begin{figure}[hbtp]
\begin{center}\includegraphics[width=.4\textwidth]{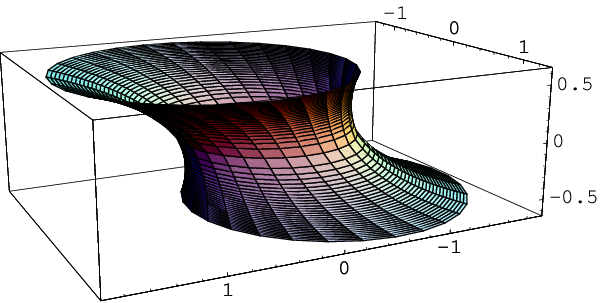},\ \includegraphics[width=.4\textwidth]{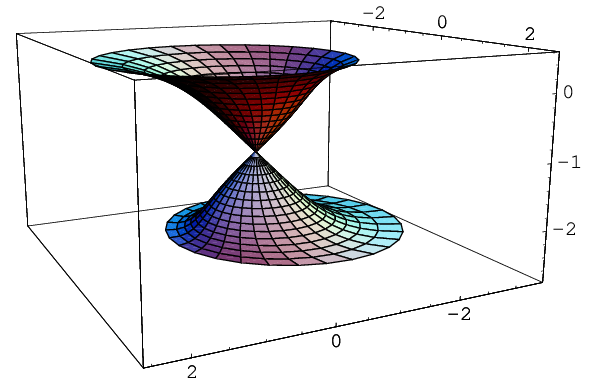}\end{center}
\caption{Minimal Riemann examples in $\e^3$ (left) and in $\e_1^3$ (right)\label{riemann}}
\end{figure}

A minimal Riemann example is  a spacelike  minimal surface foliated by circles in parallel planes. In order to motivate the problem, consider the general case that the surface is formed by a uniparametric family of circles not necessarily in parallel planes. Solving the equation $H=0$, one concludes  necessarily that the planes must be parallel. This allows to parametrize the surface as
\begin{equation}\label{rie3}
\Psi(t,s)=c(t)+r(t)(\cos(s) e_1 +\sin(s) e_2 ),
\end{equation}
where $c(t)\in\r^3$, $r(t)>0$ and $\{e_1,e_2\}$ is a couple of unit orthogonal vectors.  If $H=0$, then besides the catenoid, there are new surfaces which are called Riemann examples  (\cite{rie}). See Fig. \ref{riemann}, left. In a discrete set of heights, the surface is asymptotic to horizontal planes.

In Lorentz-Minkowski space one asks a similar problem.  First, we announce  that the results are similar and the only difference is that we have to distinguish the three cases of circles. In order to simplify the statements and proofs, we consider the case that the circles are included in parallel spacelike planes. After a rigid motion, we  suppose that they are horizontal, and so, the circles are indeed Euclidean circles. The surface parameterizes as in \eqref{rie3}.

\begin{theorem} Let $S$ be a spacelike surface foliated by circles in parallel spacelike planes. Suppose that   $H$ is constant.
\begin{enumerate}
\item If $H\not=0$, then the surface is rotational.
\item If $H=0$, then the surface is the catenoid or it belongs to a family of minimal Riemann examples in $\e_1^3$.
\end{enumerate}
\end{theorem}
\begin{proof}
We parametrize the surface as
$$\Psi(u,v)=(a(u),b(u),u)+r(u)\left(\cos(v) E_1+\sin(v) E_2\right), u\in I, v\in\r,$$
where $a,b,r>0$ are smooth functions of  $u$. The curve $c(u)=(a(u),b(u),u)$ describes the centers of the circles. The surface is rotational if and only if the functions  $a$ and $b$ are constant.

We compute the mean curvature using \eqref{media2} (for $\epsilon=-1$).  We distinguish two cases:
\begin{enumerate}
\item Case  $H\not=0$. After a homothety, we assume that  $H=1/2$. Squaring equation (\ref{media2}), we obtain
$W^3-P^2=0$. After a long computation (one may use a program as Mathematica),  this equation writes as a polynomial of type
\begin{equation}\label{6poli}
\sum_{n=0}^6 A_n(u) \cos(nv)+B_n(u) \sin(nv)=0.
\end{equation}
As the functions  $\{\cos(nv),\sin(nv):n\in{\mathbb Z}\}$ are linearly  independent, then $A_n(u)$ and $B_n(u)$ must vanish for $0\leq n\leq 6$. From $A_5=B_5=0$ we obtain
$$a'^4-10a'^2b'^2+5b'^4=0,\ \ 5a'^4-10 a'^2b'^2+b'^4=0.$$
Hence we deduce that  $a'=b'=0$, that is,  $a$ and $b$ are constant. This shows that the surface is rotational, proving the first part of the theorem.
\item Case $H=0$. Then $P=0$ gives a polynomial as  (\ref{6poli}), but only until  $n=1$. The coefficients  $A_1$ and $B_1$ imply
$$2a'r'-r a''=0,\ \ 2b'r'-rb''=0.$$
A first integration concludes that there exist constants  $\lambda$ and $\mu$ such that
$a'=\lambda r^2$ and $b'=\mu r^2$ In particular, the curve of centers lies in a plane.
\begin{enumerate}
\item If $\lambda=\mu=0$, the surface is rotational and this surface is the catenoid. See  Fig.  \ref{catenoid}, right.

\item If $\lambda\mu\not=0$, the function $r=r(u)$ satisfies
$$-1+(\lambda^2+\mu^2)r^4+r'^2-r r''=0.$$
The solution of this equation gives the Riemann examples (\cite{lo2}). See Fig. \ref{riemann}, right. Let us see that for initial conditions, the value of
  $W$ must be positive, since the surface is spacelike.
  \end{enumerate}
\end{enumerate}
\end{proof}
The same questions can be posed for timelike surfaces. See the survey \cite{lo3}  about the developments of this problem.

To end this part, we study what happens if the planes containing the circles are not parallel. Using the theory of Frenet equations for curves given in section \ref{section2}, we prove:

\begin{theorem} Consider a cmc spacelike surface in $\e_1^3$ foliated by circles contained in non-parallel spacelike planes. Then the surface is part of a hyperbolic plane.
\end{theorem}
The counterpart result in Euclidean space asserts that the only cmc surface foliated by circles in not parallel planes is the sphere.
\begin{proof} We only sketch the proof, showing the main ideas.
Let $\Pi=\Pi(u)$ be the planes of the foliation parametrized by the parameter $u$. Suppose that $\Pi(u)$ are not parallel. Consider a
curve   $\Gamma=\Gamma(u)$ orthogonal to each plane  $\Pi(u)$, where  $u$ denotes the arc length parameter of  $\Gamma$. Let $\kappa$ be the curvature of  $\Gamma$. Because $\Gamma$ is not a straight-line,   $\kappa\not=0$. Let $\{\t,\n,\b\}$ the Frenet frame of $\Gamma$. Then $\t$ is timelike. The surface    parametrizes as
$$\Psi(u,v)={\c}(u)+r(u)(\cos{v}\ \n(u)+ \sin{v}\ \b(u)),$$
where $r=r(u)>0$ and $\c={\c}(u)$ is the curve that describes the centers of circles. We write
\begin{equation}\label{alfa}
\c'=\alpha\t+\beta\n+\gamma\b,
\end{equation}
where $\alpha,\beta$ and $\gamma$ are smooth functions on  $u$. We compute again  $W$ and $P$ in \eqref{media2} and we distinguish the cases $H=0$ and $H\not=0$.
Assume that $H\not=0$. After a homothety, we assume that $H=1/2$ and that $P^2-W^3=0$.  Then
$$P^2-W^3=\sum_{n=0}^{12}\Big[A_n(u)\cos(nv)+B_n(u)\sin(nv)\Big]=0.$$

By the amount of computations and cases, we do not finish  the complete reasoning but we point out that a contradiction arrives when one
concludes that $\kappa=0$ or that $W=EG-F^2=0$. In order to show how one obtains a hyperbolic plane, we explicit how does this situation appear.
The coefficient $B_{12}$ is
$$B_{12}=\frac{r^{12}\tau^4}{2048}.$$
Thus  $\tau=0$.  The next coefficients vanish trivially until  $n=6$. For $n=6$, the equation $B_6=0$ gives some possibilities. We analyze one  of them corresponding with the case
$$\gamma=0\ \ \mbox{y}\ \ \beta^2=\kappa^2(4+r^2).$$
Now $B_5$ yields
$$B_5=\kappa^5 r^7\big(\alpha-\frac{rr'}{\sqrt{4+r^2}}\Big)=0.$$
Then  $\alpha=rr'/\sqrt{4+r^2}$. With this value of $\alpha$, the rest of coefficients $A_n$ and $B_n$ vanish trivially.  From  \eqref{alfa} we have
$$\c'=\frac{rr'}{\sqrt{4+r^2}}\t+\kappa\sqrt{4+r^2}\n=(\sqrt{4+r^2}\t)'.$$
Then  there is  $c_0\in\r^3$ such that
$$\c=c_0+\sqrt{4+r^2}\t.$$
The expression of the parametrization $\Psi(u,v)$ of the surface writes now as
$$\Psi(u,v)=c_0+\sqrt{4+r^2}\t(u)+r\cos(v)\n(u)+r\sin(v)\b(u).$$
Hence we deduce
$$\langle \Psi(u,v)-c_0,\Psi(u,v)-c_0\rangle=-4,$$
that is, the surface is the hyperbolic plane  $\h^2(2;c_0)$. Observe that $H=1/2$.
\end{proof}

\section{Elliptic equations on cmc spacelike surfaces}

The class of elliptic equations plays an important  role  in the PDE theory. The simplest case is the Laplace equation $\Delta u=0$ which appears in many areas of mathematics and physics.  The relation between elliptic equations and cmc spacelike surfaces follows by the next observation.

Assume that $S$ is a spacelike surface given by the graph $z=u(x,y)$ on the $xy$-plane, where $u\in C^\infty(\Omega)\cap C^0(\overline{\Omega})\subset\r^2$. Since the surface is spacelike, $u$ satisfies
$|Du|<1$ in $\Omega$. Consider the orientation $N$
$$N=\frac{-(Du,1)}{\sqrt{1-|Du|^2}}.$$
Then the mean curvature $H$ satisfies
\begin{equation}\label{mean2}
(1- u_y^2)u_{xx}+2 u_x u_y u_{xy}+(1-u_x^2)u_{yy}=-2H  (1-u_x^2-u_y^2)^{3/2}
\end{equation}
or
$$H=-\frac12\ \mbox{div}\Bigg(\frac{Du}{\sqrt{1-|Du|^2}}\Bigg),$$
We write this equation as
$$Qu:=\sum_{ij}a_{ij}(x,y,u_x,u_y)D_{ij}u+ b(x,y,u_x,u_y)=0, \ \ D_{ij}u=\frac{\partial^2 u}{\partial x_i\partial x_j}.$$
In order to classify this equation, we study the matrix $A=(a_{ij})$ of the coefficients of second order. This matrix  is
$$A=\Bigg(\begin{array}{cc}1- u_y^2&  u_x u_y\\ u_x u_y&1-  u_x^2\end{array}\Bigg).$$
The coefficient $a_{11}=1-u_y^2$ coincides with $|\Psi_y|^2=1-u_y^2>0$. The determinant of $A$ is $1-u_x^2-u_y^2$, which is the function $W=EG-F^2=1-u_x^2-u_y^2$.  This means that $A$ is a positive definite matrix.

\begin{proposition}
The mean curvature equation of a spacelike surface in $\e_1^3$ is elliptic.
\end{proposition}
Thus one can apply the machinery of elliptic equations to obtain geometric consequences.

It is possible to write in a unified equation the mean curvature equation in Euclidean space and in Lorentz-Minkowski space. Let
$(\r^3,(dx_1)^2+(dx_2)^2+\epsilon (dx_3)^2)$, with $\epsilon=\pm 1$. If $\epsilon=1$, then the space is $\e^3$, and if $\epsilon=-1$, then it is $\e_1^3$.
The mean curvature $H$ of a surface in $\e^3$ or a spacelike surface in $\e_1^3$ satisfies
$$(1+\epsilon u_y^2)-2\epsilon  u_x u_y u_{xy}+(1+\epsilon u_x^2)u_{yy}=2\epsilon  H  \left(1+\epsilon (u_x^2+u_y^2)\right)^{3/2}.$$

\subsection{The tangency principle}

 The maximum principle of elliptic theory is used  to compare two cmc spacelike surfaces that are tangent at a point.
Let $S_1$ and $S_2$ be two surfaces which are tangent at a common point $p\in S_1\cap S_2$.  Suppose that one of the surfaces, for example
$S_1$, lies below the other one around the point  $p$. Exactly we consider both surfaces as the graphs of smooth functions
$u_1$ and $u_2$ on a domain of the common tangent plane $T_p S_1=T_p S_2$. After an isometry, we assume that  $T_p S_i$ is the horizontal plane  $z=0$ (here we use the spacelike condition). Let us take orientations in both graphs in such way that they agree at $p$, that is,  $N_1(p)=N_2(p)$. In such a situation, we say that  $S_1$ lies below  $S_2$, and we write $S_1\leq S_2$, if $u_1\leq u_2$ in a neighbourhood of $p$ where the positive direction of the height coordinate is determined by $N_i(p)$.  If  $p\in\partial S_1\cap\partial S_2$, we add the condition  $T_p\partial S_1=T_p\partial S_2$.

An easy application of calculus proves that if  $S_1\leq S_2$, then $H_1(p)\leq H_2(p)$ (see \cite[p. 97]{mr} for the Euclidean case).
This result is known in the literature as the comparison principle.
\begin{proposition}[Comparison principle]  Let $S_1$ and $S_2$ be two spacelike surfaces tangent at $p$ and both are oriented so the Gauss maps coincide at $p$. If $S_1\leq S_2$ around $p$ then $H_1(p)\leq H_2(p)$.
\end{proposition}

\begin{proof} Take the usual parametrizations of a surface as a graph of a function. Suppose that at $p$, $N_i(p)=(0,0,-1)$. Thus, if $S_1\leq S_2$ around $p$, then $u_1\geq u_2$ around $p$. By \eqref{mean2},
\begin{equation}\label{hi}
(u_i)_{xx}(p)+(u_i)_{yy}(p)=-2H_i(p).
\end{equation}
As $u_1-u_2$ has a minimum at $p$, the Hessian of $u_1-u_2$ is positive semi-definite at $p$. In particular,
$(u_1)_{xx}(p)\geq (u_2)_{xx}(p)$ and $(u_1)_{yy}(p)\geq (u_2)_{yy}(p)$. By using \eqref{hi}, we conclude $H_1(p)\leq H_2(p)$.
\end{proof}

We ask what happens if $H_1=H_2$. Assume  that $S_1$ and $S_2$ are tangent at $p$,   $S_1\leq S_2$ at $p$ and they have the same constant mean curvature. Remark that the orientation is prescribed by the condition that the Gauss maps of $S_1$ and $S_2$ coincide  at $p$.

\begin{theorem}[Tangency principle]\label{41t2} Let $S_1$ and $S_2$ be two spacelike surfaces  with a common (interior or boundary) tangent point $p$. Suppose that $S_1\leq S_2$. If the mean curvatures agree and are constant, then $S_1=S_2$ in an open set around  $p$.
\end{theorem}
\begin{proof}
Equation  \eqref{mean2} is not linear, but if $u_1$ and $u_2$ satisfy \eqref{mean2}, the difference function $u=u_1-u_2$ satisfies a linear   elliptic equation, that is, $0=Q(u_1)-Q(u_2)=L(u)$, where $L$ is a  linear elliptic operator. Since $u_1-u_2\leq 0$ and  $u_1(p)=u_2(p)$,  the maximum principle  implies   $u=0$, showing that $u_1=u_2$.
\end{proof}
Some geometric consequences of the tangency and comparison principles are the following:

\begin{corollary}\label{c-side}  If $S$ is a   cmc spacelike compact surface with  $H\not=0$ and boundary included in a plane $P$, then the surface lies in one side of $P$.
\end{corollary}

\begin{proof}
Without loss of generality, let  $P$ be the  plane of equation $z=0$ and assume that the orientation $N$ points to the future. Suppose also that  $H>0$. We show that the surface lies below $P$. By contradiction, we assume that there exist points above $P$. We take a horizontal plane $\Pi$  in the highest point $p$ and tangent to the surface. Since $p\not\in \partial S$, then $N(p)$ is a vertical vector. As $\langle N(p),E_3\rangle<0$, then $N(p)=E_3$. Thus $S\leq\Pi$ around $p$ and the  comparison principle yields $H\leq 0$, a contradiction.
\end{proof}

\begin{corollary} Let $S$ be a minimal spacelike compact surface. Then  $S$ lies included in the spacelike convex hull  of  $\partial M$.
\end{corollary}
We point out that the spacelike convex hull of $C$ is the convex hull but using only spacelike planes.

\begin{proof} It  suffices to compare $S$ with spacelike planes. Take a spacelike plane $P$ disjoint from $S$. Next let us move $P$ parallelly towards $S$ until the first contact point $p$ with $S$. This point can not be   interior to the surface $S$ since in such a case, $S$ and $P$ are tangent at $p$ and both surfaces are minimal with any orientation. The tangency principle implies that $S\subset P$, a contradiction. As a conclusion, the point $p$ must be a boundary point of $S$. If we do the same arguments with any spacelike plane of $\e_1^3$, we obtain the result.
\end{proof}

The next result says that two graphs with the same mean curvature and the same boundary must coincide. This is a direct  consequence of the maximum principle of a solution of a quasilinear elliptic equation.  We do a proof by using the tangency principle.

\begin{corollary}\label{42t1} Let  $S_1, S_2\subset\e_1^3$ be two compact  graphs with the same constant mean curvature and the same  boundary. Then    $S_1=S_2$.
\end{corollary}

\begin{proof}
We suppose that the surfaces are graphs on a horizontal plane $P$. Let $C=\partial S_1=\partial S_2$. Let $S_1$ and we lift up it vertically  until it does not touch
 $S_2$ (this is possible because both surfaces are compact). Next, we descend $S_1$ until the first time that it touches $S_2$.
 If there is an interior tangent point  $p\in S_1\cap S_2$,   $S_2\leq S_1$ around $p$. As $S_1$ and $S_2$ have the same mean curvature,  the tangency principle says that both surfaces agree around  $p$. By connectedness,  $S_1=S_2$.

If $p$ is a boundary point, then  $S_1$ comes back to its original position and, furthermore, $S_2\leq S_1$. There are two possibilities:
\begin{enumerate}
\item If $p$ is a tangent point, the tangency principle concludes that  $S_1=S_2$ again.
\item On the contrary, the slope of $S_1$ along $C$ is strictly bigger than the one of  $S_2$. Now we descend  $S_1$ until it does not touch $S_2$. Next, we move it upwards. Since in its original position, $S_1$ lies strictly above $S_2$ along $C$, then there exists
    an interior tangent point between  $S_1$ and $S_2$ at some time $t$ before the original position. If we denote
by $S_1(t)$ the surface $S_1$ at this position, the tangency principle yields $S_1(t)=S_2$. But this is impossible because  $\partial S_1(t)\not=C$ and $\partial S_2=C$. This contradiction says that this situation can not occur.
\end{enumerate}

\end{proof}

For the next result, we need to introduce the concept of hyperbolic cap.   Denote by  $\h^2(r)$ the hyperbolic plane of radius
 $r$ and center the origin. If $R>0$,  a \emph{hyperbolic cap} is defined by
$$K(r;R)=\{(x,y,z)\in \h^2(r): z^2\leq r^2+ R^2\}.$$
The boundary of $K(r;R)$, namely,
$$\partial K(r;R)=\{(x,y,z)\in K(r;R): z^2=r^2+ R^2\},$$
is a circle of radius $R$  in the spacelike plane  of equation  $z=\sqrt{R^2+r^2}$.

\begin{corollary}The only cmc compact spacelike surfaces in $\e_1^3$ spanning a circle are a planar disc and a hyperbolic cap.
\end{corollary}

\begin{proof} Let $C\subset\e_1^3$ be a circle of radius $R>0$. Since the surface is compact, this circle is a spacelike curve contained in a spacelike plane $\Pi$  which we assume that it is the $xy$-plane. Because a circle  is a simple closed curve, the surface must be a graph on $\Pi$ by Prop. \ref{co-closed}.  Let $S$ be a graph with mean curvature $H$ and  $\partial S=C$. Then  $S$ is a graph on the round disc bounded by  $C$. If $H=0$, the planar disc that bounds $C$ is a graph with $H=0$. By Cor. \ref{42t1}, $S$ is a planar disc.

Let  $H\not=0$. We have only to show that there exists a hyperbolic cap with mean curvature  $H$ and boundary  $C$. For this, we take the hyperbolic cap  $K(1/H;R)$.
\end{proof}

As we observe, this result has a simple proof, but its counterpart in Euclidean space is completely different. In $\e^3$  there are compact cmc surfaces with circular boundary that are not umbilical. These surfaces are not embedded and have higher genus. Furthermore, it is an open problem to know if planar discs and spherical caps are the only cmc surfaces with circular boundary that are embedded or that are topological discs. See \cite{lo4} for an introduction to the problem.

\subsection{Estimates on cmc spacelike surfaces}

In this section we compute the Laplacian of the coordinate functions of  a spacelike surface as well as of   the Gauss map.  In order to compare the differences with the Euclidean case,  we work in both ambient spaces.  Let  $x:M\rightarrow\e_1^3$ or $\e^3$ be an
  immersion of a surface, which will be spacelike if the codomain if $\e_1^3$.   Let $a\in \r^3$. We  compute $\Delta\langle x,a\rangle$ and $\Delta\langle N,a\rangle$, where $\Delta$ is the Laplace operator on $M$ with the induced metric from $\e^3$ or $\e_1^3$ depending the case. The Laplacian of a  function  $f\in C^{\infty}(M)$ is defined as
$$\Delta f=\mbox{ div} (\nabla f)=\mbox{trace}\Big(v\longmapsto\nabla_v \nabla f\Big).$$
Let $p\in M$ and take an adapted orthonormal basis at $p$, that is, an orthonormal basis of tangent vector fields
 $\{E_1,E_2\}$ such that  $\nabla_{E_i(p)}E_j=0$, $i,j\in\{1,2\}$. Because the computations are done at $p$, we indicate
  $e_i=E_i(p)$. With respect to this basis, the Laplacian is
$$\Delta f(p)=\sum_{i=1}^2 e_i(E_i(f)).$$

\begin{enumerate}
\item Let $f=\langle x,a\rangle$. Then $E_i(f)=\langle E_i,a\rangle$ and at the point $p$ we have
\begin{eqnarray*}
e_i(E_i(f))&=&e_i\langle E_i,a\rangle=\langle\nabla^0_{e_i} E_i,a\rangle=\langle \nabla_{e_i}E_j+\sigma(e_i,e_j),a\rangle\\
&=&\langle \sigma (e_i,e_j),a\rangle=\epsilon\langle Ae_i,e_i\rangle\langle N,a\rangle.
\end{eqnarray*}
Thus
$$\Delta\langle x,a\rangle=\epsilon\Big(\sum_{i=1}^2\langle Ae_i,e_i\Big)\langle N,a\rangle=2H\langle N,a\rangle.$$
\item  Suppose  that $H$ is constant. Let $f=\langle N,a\rangle$. First we show that the vector $ \sum_{i=1}^2\nabla^0_{e_i}\nabla^0_{E_i} N$ is orthogonal to the surface, that is,
\begin{equation}\label{41claim}
\sum_{i=1}^2\langle \nabla^0_{e_i}\nabla^0_{E_i} N,e_k\rangle=0,\ \ k=1, 2.
\end{equation}
Since $\langle N,E_k\rangle=0$, then $\langle\nabla^0_{E_i}N,E_k\rangle=-\langle N,\nabla^0_{e_i}E_k\rangle$. Therefore
$$ \langle\nabla^0_{e_i}\nabla^0_{E_i} N,e_k\rangle+\langle\nabla^0_{e_i}N,\nabla^0_{e_i}E_k\rangle=-\langle\nabla^0_{e_i}N,\nabla^0_{e_i}E_k\rangle -
\langle N,\nabla^0_{e_i}\nabla^0_{e_i}E_k\rangle.$$
As $\nabla^0_{e_i}N$ is tangent to the surface and  $\nabla_{e_i}E_k=0$, then
$$\langle\nabla^0_{e_i}\nabla^0_{E_i}N,e_k\rangle=-\langle N,\nabla^0_{e_i}\nabla^0_{e_i}E_k\rangle=-\langle N,\nabla_{e_i}^0\nabla^0_{E_k}E_i\rangle.$$
The metric in  $\e^3$ and $\e_1^3$ is flat and this means  $\nabla^0_{e_i}\nabla^0_{E_k}E_i=\nabla^0_{e_k}\nabla^0_{E_i}E_i$. By substituting in the above equation, we obtain:
$$\langle\nabla^0_{e_i}\nabla^0_{E_i}N,e_k\rangle=-\langle\nabla^0_{e_k}\nabla^0_{E_i}E_i,N\rangle$$
and
\begin{equation}\label{41ayuda1}\sum_{i=1}^2\langle\nabla^0_{e_i}\nabla^0_{E_i}N,e_k\rangle=-\langle\nabla_{e_k}^0\Big(\sum_{i=1}^2 \nabla^0_{E_i}E_i\Big),N\rangle.
\end{equation}
The mean curvature is $2H=\epsilon\sum_{i=1}^2\langle\nabla^0_{e_i}E_i,N\rangle$. Because  $H$ is constant,
$$e_k\langle\sum_{i=1}^2\nabla^0_{e_i}E_i,N\rangle=0,$$
that is,
$$\langle\nabla_{e_k}^0\Big(\sum_{i=1}^2 \nabla^0_{E_i}E_i\Big),N\rangle+\langle\sum_{i=1}^2 \nabla^0_{E_i}E_i,\nabla^0_{e_k}N\rangle=0.$$
The second summand vanishes since that the left part is orthogonal to the surface and the right one is tangent. From  (\ref{41ayuda1}),
$$\langle\nabla_{e_k}^0\Big(\sum_{i=1}^2 \nabla^0_{E_i}E_i\Big),N\rangle=0,$$
and   the claim is proved.

Finally,  \eqref{41claim} yields
\begin{eqnarray}
\Delta\langle N,a\rangle&=& \epsilon \sum_{i=1}^2\langle\nabla_{e_i}^0\nabla^0_{E_i}N,N\rangle\langle N,a\rangle=-\epsilon \sum_{i=1}^2\langle\nabla^0_{e_i}N,\nabla^0_{e_i}N\rangle \langle N,a\rangle\\
&=&-\epsilon \sum_{i=1}^2\langle Ae_i,Ae_i\rangle\langle N,a\rangle=-\epsilon \sum_{i=1}^2\langle A^2e_i,e_i\rangle\langle N,a\rangle\\
&=&-\epsilon \mbox{ trace }(A^2)\langle N,a\rangle.\label{41ayuda2}
\end{eqnarray}
\end{enumerate}

\begin{theorem}\label{4th-two} Let $x:M\rightarrow\e_1^3$ or $\e^3$ be an  immersion and let $N$ be an orientation on $M$. Given $a\in\r^3$, we have
\begin{equation}\label{lap1}\Delta\langle x,a\rangle=2H\langle N,a\rangle.
\end{equation}
Furthermore, if  $H$ is constant,
$$\Delta\langle N,a\rangle=-\epsilon \mbox{ trace} (A^2)\langle N,a\rangle.$$
\end{theorem}

By \eqref{rho},
 $$\mbox{trace }(A^2)=4H^2-2\epsilon K=2H^2+2(H^2-\epsilon K),$$
and  $H^2-\epsilon K\geq 0$ on $M$.  Thus, the last equation writes as
\begin{equation}\label{lap2}
\Delta\langle N,a\rangle+\epsilon(4H^2-2\epsilon K)\langle N,a\rangle=0.
\end{equation}
Thanks to \eqref{lap1} and \eqref{lap2} we will obtain a priori estimates of the height of a cmc graph in $\e^3$ which are not possible to extend to  $\e_1^3$.

In both settings consider a cmc graph of a function $u\in C^{\infty}(\Omega)\cap C^{0}(\overline{\Omega})$, where $\Omega$ is  a bounded domain.  In particular, the graph is a compact set. Assume that the boundary of the graph is $\partial\Omega$, that is, $u=0$ on $\partial\Omega$ and the orientation of the graph points up.
Let $a=(0,0,1)$. Thus $\langle N,a\rangle>0$ if $\epsilon=1$ and $\langle N,a\rangle<0$ if $\epsilon=-1$. In fact, and since $N$ and $a$ are timelike vectors in the same timelike cone, $\langle N,a\rangle\leq -1$. On the other hand, Cor. \ref{c-side} implies that $\epsilon\langle x,a\rangle=u<0$ on $\Omega$.

A linear combination of
\eqref{lap1} and \eqref{lap2} gives
$$\Delta \Big(\epsilon H\langle x,a\rangle+\langle N,a\rangle\Big)=-2\epsilon(H^2-\epsilon K)\langle N,a\rangle\leq 0.$$
The maximum principle for elliptic equations asserts that
\begin{equation}\label{fails}
\epsilon H\langle x,a\rangle+\langle N,a\rangle \geq
\min_{\partial\Omega}(\epsilon H\langle x,a\rangle+\langle N,a\rangle)=  \min_{\partial\Omega}\langle N,a\rangle.
\end{equation}
In Euclidean space, $\langle N,a\rangle\geq 0$, and thus, we conclude
$$  H\langle x,a\rangle+\langle N,a\rangle \geq 0.$$
Hence
$$H x_3\geq -\langle N,a\rangle\geq -1\ \Rightarrow 0\geq x_3\geq\frac{-1}{H}.$$
This height estimate depends \emph{only} on $H$. Recall that $x_3=u$ measures the height of the graph with respect to $\Pi$.

\begin{proposition} \label{pr-he} Let $P$ be a plane and $H\not=0$. Given a compact graph  $S\subset\e^3$ on  $P$ with constant mean curvature  $H$ and $\partial S\subset P$, the height of   $S$ with respect to $P$ is less than  $1/|H|$.
\end{proposition}

For a cmc spacelike surface, the last part in \eqref{fails} can not follow because $\langle N,a\rangle\leq -1$. In fact, prescribing the value $H$ of the mean curvature,  there are spacelike graphs with planar boundary and constant mean curvature $H$ that have arbitrary  heights. It suffices to consider hyperbolic caps  $K(1/H;R)$ with arbitrary $R$, placing the boundaries of all them  on a fix plane.

\begin{proposition}\label{pr-hl} Let $\Omega$ be a bounded domain of a spacelike plane $P$ and  $H\in\r$. Let $S=\mbox{graph}(u)$ be  with mean curvature  $H$ spanning  $\partial \Omega$. Then there exists a constant  $c=c(H,\Omega)$ such that  $|u|\leq c(H,\Omega)$.
\end{proposition}

\begin{proof} Without loss of generality, we assume that $P$ is the plane $\Pi$ of equation $z=0$. Let $H>0$ be the mean curvature with the future orientation. We know by Cor. \ref{c-side}  that  $u\leq 0$. Consider a hyperbolic cap $K(1/H;R)$ whose mean curvature is $H$ with respect to the orientation pointing to the future. Let  $R$ sufficiently big so the circle $\partial K(1/H;R)$ lies included in  $\Pi$ and that $\Omega$ lies strictly in the bounded domain determined by
$\partial K(1/H;R)$. By a translation, we move downwards $K(1/H;R)$ until to be disjoint from  $S$. Now, we move it up. Recall that the mean curvatures agree but  the boundaries are different. By the tangency principle, there is not a contact point at least that the hyperbolic cap returns  its original position.
 This proves that the height of $S$ is less than the one of $K(1/H;R)$.
\end{proof}

We now study  $C^1$-estimates of a cmc graph. First we prove
\begin{equation}\label{es-c1}
\sup_{\Omega}|Du|=\max_{\partial\Omega}|Du|.
\end{equation}

From \eqref{lap2} and for $a=(0,0,1)$, we have
$$\Delta \langle N,a\rangle=-\epsilon(4H^2-2\epsilon K)\langle N,a\rangle\leq 0.$$
Then the maximum principle gives
$$\langle N,a\rangle\geq\min_{\partial\Omega}\langle N,a\rangle,\ \ \langle N,a\rangle=\epsilon \frac{1}{\sqrt{1+\epsilon|Du|^2}},$$
that is,
$$\epsilon \frac{1}{\sqrt{1+\epsilon |Du|^2}}\geq \epsilon \min_{\partial\Omega}\frac{1}{\sqrt{1+\epsilon |Du|^2}},$$
hence it follows \eqref{es-c1}.

We write \eqref{es-c1} in geometric terms. Let $\nu$ be
the interior conormal vector along $\partial S$. If $H>0$, the surface lies below the plane. Moreover,
$$\langle\nu,a\rangle^2+\epsilon \langle N,a\rangle^2=\epsilon.$$
$$\langle\nu,a\rangle^2=\frac{|Du|^2}{1+\epsilon|Du|^2}\Longrightarrow |Du|=\frac{|\langle \nu,a\rangle|}{\sqrt{1-\epsilon \langle\nu,a\rangle^2}}.$$
If $\epsilon=1$, let $\delta=\min_{\partial M}\langle\nu,a\rangle>-1$, then
$$\max_{\partial \Omega}|Du|=-\frac{\delta}{\sqrt{1-\delta^2}}.$$
Let us observe  that $\langle\nu,a\rangle$ measures the slope of $S$ with respect to $\Pi$ along $\partial S$.

\subsection{The Dirichlet problem in the Lorentzian case}

In Lorentz-Minkowski space, the Dirichlet problem ${\mathcal P}$ for the constant mean curvature equation is
$$
\left\{
\begin{array}{lll}
 & \mbox{div}\dfrac{Du}{\sqrt{1-|Du|^2}}=2H&\mbox{on $\Omega$,}\\
& u=0 &\mbox{in $\partial\Omega$.} \\
&|Du|<1  &\mbox{on $\partial\Omega$.}
\end{array}
\right.$$
The technique to solve this problem is the continuity method and that we explain now.   Denote by ${\mathcal P}_t$ the  Dirichlet problem which is the same than ${\mathcal P}$ except that we replace $H$ by $tH$, where $t\in [0,1]$. Define
$$A=\{t\in [0,1]: \mbox{ there exists a solution of ${\mathcal P}_t$}\}.$$
For $t=0$, the function $u=0$ is a solution. If we prove that the set $A$ is open and closed in $[0,1]$, by connectedness, $A=[0,1]$, in particular, $1\in A$, solving the initial Dirichlet problem ${\mathcal P}_1={\mathcal P}$. The set   $A$ is open as a consequence of the implicit function theorem in Banach space. We omit the details.

For the closeness, we need to get a priori $C^0$ and $C^1$ estimates of a solution of ${\mathcal P}_t$ for all $t\in [0,1]$,  that is, of $|u_t|$ and $|Du_t|$, where $u_t$ is the solution of ${\mathcal P}_t$. The estimates of $|u_t|$ are given by Props. \ref{pr-he} and \ref{pr-hl}. In order to estimate $|Du_t|$, it suffices to do it along $\partial \Omega$ thanks to \eqref{es-c1}.

Here there is a difference between the Euclidean and Lorentzian ambient space. To show these differences, suppose that $\Omega$ is a convex domain.  In Euclidean space, the domain $\Omega$ can not arbitrary large. In fact, if the Dirichlet problem has a solution, there do not exist a closed disc of radius $1/|H|$ included in $\Omega$. On the contrary, we can place a hemisphere $z=\sqrt{1/H^2-x^2-y^2}$ of radius $1/|H|$ over $\Omega$. Then moving up and next down, the hemisphere arrives until the graph  at a tangent point. But the hemisphere would lie above the graph at the tangent point, in contradiction with the tangency principle. Therefore, we have to impose restriction on the size of $\Omega$. An example of an existence result is the following (\cite{lot}):

\begin{theorem}\label{ex-e}
 Let $\Omega\subset\r^2$ be a strictly convex domain included in a plane  $P$  and let $\kappa$ be the curvature of  $\partial \Omega$.
If $H\in\r$ satisfies
 $\kappa>|H|>0$, then there exists a graph on  $\Omega$ with mean curvature  $H$ and boundary  $\partial \Omega$.
 \end{theorem}

In contrast, in Lorentz-Minkowski space we have:

\begin{theorem}
If $\Omega\subset\r^2$ is a bounded convex domain, then the Dirichlet problem has a solution (for any $H$).
\end{theorem}
\begin{proof}
Let $w$ be the diameter of $\overline{\Omega}$. This number has the property that given a horizontal direction $v\in \Pi$, there are two parallel straight-lines $\{L_v,L'_v\}$ contained in $\Pi$ and orthogonal to $v$ with the property that   $\overline{\Omega}$ is included in the strip determined by $\{L_v,L'_v\}$. Moreover, the distance between $L_v$ and $L'_v$ is less thant $w$. Consider the hyperbolic cylinder of equation $y^2-z^2=-1/(4H^2)$, $z>0$. This surface has (constant) mean curvature $H$ with the orientation pointing to the future (see Ex. \ref{ex-cy}). Take the   piece $C$  obtained when we cut   the cylinder by a horizontal plane in such way that the intersection are two  parallel straight-lines $w$ far apart. This surface $C$ has bounded height. Let us observe that  $C$ is a graph of a function $f$, which it is , up a constant,  $f(x,y)=\sqrt{1/(4H^2)+y^2}$) where the domain is a strip $\Omega_C\subset \Pi$ of width $w$.

We are in conditions to estimate $|Du_t|$ at any boundary point. Take $p\in\partial\Omega$. Move down $C$ and rotating with respect to $z$-axis, if necessary, until that we place $C$ below the graph $S_t$ of $u_t$, $\overline{\Omega}\subset \Omega_C$ and $p\in\partial\Omega_C$. Next move vertically down $C$ sufficiently until that $C$ does not intersect $S_t$. Next we move up. Because the mean curvature of $C$ is $H$, the tangency principle implies that there is not a contact point between $C$ and $S_t$ at least that  $\partial C$ touches $\Pi$. Then the surface $S_t$ lies above $C$ and $p$ is a common point of both surfaces. Thus $|Du_t|(p)$ is bounded by  $|Df|(p)$, the gradient of the cylinder at $p\in\partial\Omega$. However $|Df|(p)$ is a constant that does not depend on $S_t$ but only on $\Omega$ and $H$. This gives the a priori $C^1$-estimates of $u_t$ along $\partial\Omega$ that we are going to looking for.
\end{proof}

Similar arguments as above  proves Th. \ref{ex-e} by replacing hyperbolic cylinders by  spherical caps of $\e^3$.


\end{document}